\documentclass[11pt,onecoulme]{article}
\usepackage{amsfonts}
\usepackage{mathrsfs}
\usepackage{amssymb}
\usepackage{color, amsmath,amssymb, amsfonts, amstext,amsthm, latexsym}

\usepackage{amssymb, epsfig, amssymb, latexsym}
\usepackage{amsmath}
\usepackage{graphicx}
\usepackage{longtable}
\usepackage{appendix}
\DeclareMathOperator{\esssup}{esssup}

\numberwithin{equation}{section}

\allowdisplaybreaks

\newtheorem{definition}{Definition}[section]
\newtheorem{theorem}[definition]{Theorem}
\newtheorem{lemma}[definition]{Lemma}
\newtheorem{remark}[definition]{Remark}

\newtheorem{hyp}[definition]{Hypothesis}

  \renewcommand\appendix{\par
    \setcounter{section}{0}
    \setcounter{subsection}{0}
    \gdef\thesection{ Appendix \Alph{section}}}
\title{Viscosity Solutions to Second Order Path-Dependent Hamilton-Jacobi-Bellman Equations and Applications \thanks{This work was partially supported by  the National Natural Science Foundation of China  (Grant No. 11401474), Shaanxi Natural Science Foundation
               (Grant No. 2021JM-083) and the Fundamental Research Funds for the Central Universities (Grant No. 2452021063).}}

\author{  Jianjun Zhou   \\
       College of Science,
             Northwest A\&F University,\\ Yangling 712100, Shaanxi, P. R.
             China\\
      \emph{E-mail:zhoujianjun@nwsuaf.edu.cn} }
      \date{}

\begin{document}

\maketitle

\pagestyle{plain}

\begin{abstract}
In this article, a notion of viscosity solutions is introduced for second order path-dependent Hamilton-Jacobi-Bellman (PHJB) equations associated with optimal control problems for path-dependent stochastic differential equations. We identify the value functional of optimal control problems as unique viscosity solution to the associated PHJB equations. We also show that our notion of viscosity solutions is consistent with the corresponding notion of classical solutions, and satisfies a stability property. Applications to backward stochastic Hamilton-Jacobi-Bellman  equations are also given.

\medskip

\end{abstract}


\section{Introduction}

The notion of  viscosity solutions for  Hamilton-Jacobi-Bellman (HJB) equations, first introduced in 1983 by Crandall and Lions \cite{cra1}, has  become an indispensable tool in optimal control theory and numerous subjects
related to it. We refer to the survey paper of Crandall, Ishii and  Lions \cite{cran2} and  the monographs of Fleming
      and Soner \cite{fle1} and Yong and Zhou \cite{yong} for a detailed account
for the theory of viscosity solutions. For viscosity solutions in infinite dimensional Hilbert spaces, we refer to Fabbri, Gozzi and  \'{S}wi\c{e}ch \cite{fab1}, Gozzi,  Rouy and \'{S}wi\c{e}ch \cite{gozz3}, Lions \cite{lio1, lio2, lio3} and  \'{S}wi\c{e}ch \cite{swi}.
\par
   Dupire in his work \cite{dupire1} introduced horizontal and vertical derivatives for functionals defined in the space of c\`adl\`ag paths  and provided a functional It\^o formula
   (see Cont and  Fourni\'e \cite{cotn0, cotn1} for a more general and systematic research). Soon after, Dupire's functional It\^o formula was applied to second order path-dependent HJB (PHJB) equations.
Peng \cite{peng2.5} made the first attempt to extend Crandall-Lions framework to path-dependent case, by focusing on the uniqueness part. Tang and Zhang \cite{tang1} proposed a different notion of viscosity solutions for path-dependent Bellman equations. They verified that the value functional of the corresponding optimal control problem is a viscosity solution, but did not investigate the comparison principle.

\par
         At the same time,      Ekren, Keller, Touzi
and Zhang \cite{ekren1} introduced a new notion of viscosity solutions for semi-linear path-dependent partial differential equations (PPDEs) in the space of continuous paths in terms of a nonlinear expectation. 
 In the subsequent
works, Ekren \cite{ekren0}, Ekren,
Touzi and Zhang  
\cite{ekren3, ekren4}  and Ren \cite{ren}  extended  the notion   to  fully nonlinear case when the  Hamilton function $\mathbf{H}$ is uniformly nondegenerate.
Ren,  Touzi and  Zhang \cite{ren1}  considered the degenerate case and established  the comparison principle  when  the nonlinearity  $\mathbf{H}$ is $d_p$-uniformly continuous in the path.
  We refer to Barrasso and  Russo \cite{bar}, Cosso and Russo \cite{cosso}, Leao, Ohashi and Simas \cite{lea}, Peng and  Song \cite{penga} and Peng and  Wang \cite{pengb} for other notions of
solutions to path-dependent semi-linear equations, and  to Viens and  Zhang \cite{viens} and Wang, Yong and  Zhang \cite{wang} for some more general  PPDEs.
  We also mention that Lukoyanov \cite{luk} developed
a theory of viscosity solutions to fully non-linear path-dependent first order  Hamilton-Jacobi equations when Hamilton function $\mathbf{H}$ is $d_p$-locally Lipschitz continuous with respect to the path.
 \par

In this paper, we consider the following   controlled path-dependent stochastic  differential
                 equation (PSDE):
\begin{eqnarray}\label{state1}
\begin{cases}
            dX^{\gamma_t,u}(s)=
           b(X_s^{\gamma_t,u},u(s))ds+\sigma(X_s^{\gamma_t,u},u(s))dW(s),  \ \ s\in [t,T],\\
~~~~~X_t^{\gamma_t,u}=\gamma_t\in {\Lambda}_t.
\end{cases}
\end{eqnarray}
In the above equation,  $\Lambda_t$ denotes the set of all
 continuous $\mathbb{R}^{d}$-valued functions $\gamma$ defined over $[0,t]$, and  let ${\Lambda}^s=\bigcup_{l\in[s,T]}{\Lambda}_{l}$ and ${\Lambda}$ denote ${\Lambda}^0$; 
   the unknown $X^{\gamma_t,u}(s)$, representing the state of the system, is an $\mathbb{R}^{d}$-valued process; $X^{\gamma_t,u}_s$ is the whole history of  $X^{\gamma_t,u}(\cdot)$ from time 0 to $s$;
 $\{W(t),t\geq0\}$ is an $n$-dimensional
      standard Wiener process; 
      $u(\cdot)=(u(s))_{s\in [t,T]}$ is  progressively measurable with respect to the Wiener filtration  and takes values in some Polish space $(U,d_1)$ (we will say that $u(\cdot)\in {\cal{U}}[t,T]$).
We define a   norm on ${\Lambda}_t$  and a metric on ${\Lambda}$ as follows: for any  
  $(t,\gamma_t), (s,\eta_s)\in [0,T]\times{\Lambda}$,
\begin{eqnarray*}
   ||\gamma_t||_0:=\sup_{0\leq l\leq t}|\gamma_t(l)|,\ \ d_{\infty}(\gamma_t,\eta_s)
   :=|t-s|
               +\sup_{0\leq l\leq t\vee s}\left|\gamma_{t}(l\wedge t)-\eta_{s}(l\wedge s)\right|.
\end{eqnarray*}
 We assume that the coefficients $b:\Lambda\times U\rightarrow \mathbb{R}^{d}$ and  $\sigma:\Lambda\times U\rightarrow \mathbb{R}^{d\times n}$   satisfy  Lipschitz condition under $||\cdot||_0$
                         with respect to  the path.
\par
              We aim at maximizing  a cost functional of the form:
\begin{eqnarray}\label{cost1}
                     J(\gamma_t,u(\cdot)):=Y^{\gamma_t,u}(t),\ \ \ (t,\gamma_t)\in [0,T]\times {\Lambda},
\end{eqnarray}
 over  ${{\mathcal
                  {U}}}[t,T]$,
   where the process $Y^{\gamma_t,u}$ is defined via solution of backward stochastic differential equation (BSDE):
    \begin{eqnarray}\label{fbsde1}
Y^{\gamma_t,u}(s)&=&\phi(X_T^{\gamma_t,u})+\int^{T}_{s}q(X_l^{\gamma_t,u},Y^{\gamma_t,u}(l),Z^{\gamma_t,u}(l),u(l))dl\nonumber\\
                 &&-\int^{T}_{s}Z^{\gamma_t,u}(l)dW(l),\ \ \ a.s., \ \ \mbox{all}\ \ s\in [t,T].
\end{eqnarray}
Here  $q$ and $\phi$ are given real functionals on ${\Lambda}\times \mathbb{R}\times \mathbb{R}^n\times U$ and ${\Lambda}_T$, respectively, and satisfy  Lipschitz condition under $||\cdot||_0$
                         with respect to  the path.
             We define the value functional of the  optimal
                  control problem as follows:
\begin{eqnarray}\label{value1}
V(\gamma_t):=\mathop{\esssup}\limits_{u(\cdot)\in{\cal{U}}[t,T]} Y^{\gamma_t,u}(t),\ \  (t,\gamma_t)\in [0,T]\times {\Lambda}.
\end{eqnarray}
  We characterize this value functional  $V$ with the following
                     PHJB equation:
  \begin{eqnarray}\label{hjb1}
\begin{cases}
{\mathcal{L}}V(\gamma_t):=\partial_tV(\gamma_t)+{\mathbf{H}}(\gamma_t,V(\gamma_t),\partial_xV(\gamma_t),\partial_{xx}V(\gamma_t))= 0,\ \ \  (t,\gamma_t)\in
                               [0,T)\times {\Lambda},\\
 V(\gamma_T)=\phi(\gamma_T), \ \ \ \gamma_T\in {\Lambda}_T;
 \end{cases}
\end{eqnarray}
      where
\begin{eqnarray*}
                                {\mathbf{H}}(\gamma_t,r,p,\iota)&=&\sup_{u\in{
                                         {U}}}[
                        \langle p,b(\gamma_t,u)\rangle+\frac{1}{2}\mbox{tr}[ \iota \sigma(\gamma_t,u)\sigma^\top(\gamma_t,u)]\\
                        &&\ \ \ \ \   +q(\gamma_t,r,\sigma^\top(\gamma_t,u)p,u)], \ \ \ (t,\gamma_t,r,p,\iota)\in [0,T]\times{\Lambda}\times \mathbb{R}\times \mathbb{R}^{d}\times \mathcal{S}(\mathbb{R}^{d}).
\end{eqnarray*}
Here,    $\sigma^\top$ is the transpose of the matrix $\sigma$,  $\mathcal{S}(\mathbb{R}^{d})$  the set of all $(d\times d)$ symmetric matrices,  $\langle\cdot,\cdot\rangle$ the scalar product of $\mathbb{R}^{d}$,
  and
                         $\partial_t,\partial_x$ and $\partial_{xx}$  the so-called pathwise  (or functional or Dupire; see \cite{dupire1, cotn0, cotn1}) derivatives, where $ \partial_t$  is known as the horizontal derivative, while $\partial_x$ and $\partial_{xx}$  are the first and  second order vertical derivatives, respectively.
\par
                         The primary objective of this article is to
                         develop a concept  of  viscosity solutions
                         to
                         PHJB equations on the space of continuous paths (see Definition \ref{definition4.1} for details). We shall show that the value functional
                         $V$  defined in  (\ref{value1}) is   unique viscosity solution to the PHJB
                         equation  (\ref{hjb1}) when the coefficients ($b,\sigma,q,\phi$) are uniformly Lipschitz in the path  under $||\cdot||_0$.
\par
The main challenge for our path-dependent case comes from  both facts that the path space $\Lambda_T$ is an infinite dimensional Banach space, and that the maximal norm $||\cdot||_0$ is not G\^ateaux differentiable.
         Since $\Lambda_T$ is not a separable Hilbert space,
             the standard techniques for the comparison
principle in Hilbert space
introduced by Lions \cite{lio1, lio2, lio3}, which contain
a limiting procedure based on the existence of a countable basis, are not applicable in
                 our case. On the other hand,  noticing  that   the value functional is only $||\cdot||_0$-Lipschitz continuous  with respect to   the path,  the auxiliary
functional in the proof of uniqueness should include the term
$||\cdot||_0$ or a functional which is equivalent to $||\cdot||_0$.
 The lack
of smoothness of $||\cdot||_0$ makes it more difficult to define the viscosity solutions and to prove its uniqueness.
\par
In this paper we  want to extend the theory of viscosity solutions to the second order path-dependent case.
We adopt the natural generalization of the well-known Crandall-Lions definition  in terms of test functions.
Since we  assume the coefficients ($b,\sigma,q,\phi$) only satisfy $||\cdot||_0$-Lipschitz conditions with respect to  the path and do not impose  uniformly nondegenerate  requirement on the coefficients, none of these results in Ekren \cite{ekren0}, Ekren,
Touzi and Zhang
\cite{ekren3, ekren4}, Lukoyanov \cite{luk}, Peng \cite{peng2.5}, Ren \cite{ren} and
Ren,  Touzi and  Zhang \cite{ren1}    are directly applicable to our case.
\par
 The main contribution of this paper
is the introduction of an appropriate notion of viscosity solutions and the proof of uniqueness. The uniqueness property is derived from the comparison theorem.
For the proof of  the comparison theorem, we generalize the classical methodology in Crandall, Ishii  and  Lions \cite{cran2}  to the path-dependent
case.   We introduce functional $\Upsilon:\Lambda\times \Lambda\rightarrow {\mathbb{R}}$ defined by
          $$
          \Upsilon(\gamma_t,\eta_s)=S(\gamma_t,\eta_s)+3|\gamma_t(t)-\eta_s(s)|^6
          $$
          and
 \begin{eqnarray*}
S(\gamma_t,\eta_s)=\begin{cases}
            \frac{(||\gamma_{t}-\eta_s||_{0}^6-|\gamma_{t}(t)-\eta_s(s)|^6)^3}{||\gamma_{t}-\eta_s||^{12}_{0}}, \
         ~~ ||\gamma_{t}-\eta_s||_{0}\neq0; \\
0, \ ~~~~~~~~~~~~~~~~~~~ ~~~~~~~~~~~||\gamma_{t}-\eta_s||_{0}=0
\end{cases}
\end{eqnarray*}
for $\gamma_t,\eta_s\in \Lambda$. Here $||\gamma_{t}-\eta_s||_{0}=\sup_{l\in [0,t\vee s]}|\gamma_{t}(l\wedge t)-\eta_s(l\wedge s)|$.
\par
This key functional  is the starting point for the proof of  comparison theorem. 
 First, for every fixed $(\hat{t},a_{\hat{t}})\in [0,T)\times\Lambda$,   define $f: \Lambda^{\hat{t}}\rightarrow \mathbb{R}$ by
 $$
 f(\gamma_t):=\Upsilon(\gamma_t,a_{\hat{t}}), \ \  \gamma_t\in\Lambda^{\hat{t}}.
 $$
 We show that it is equivalent to $||\cdot||_0^6$ and  study its  regularity in the sense of  horizontal/vertical derivatives. 
   Then the test function in our definition of viscosity solutions and the auxiliary function $\Psi$ in the proof of  comparison theorem can include $f$ (see Step 1 in the proof of Theorem \ref{theoremhjbm}) as we show that it satisfies  a functional It\^o formula.
     By this,  the comparison theorem is  established
            when the coefficients only satisfy   Lipschitz assumption under $||\cdot||_0$.
\par
 Second, we use  $\Upsilon$ to  define a  smooth  gauge-type function $\overline{\Upsilon}:\Lambda\times \Lambda\rightarrow \mathbb{R}$ by
                 $$
                 \overline{\Upsilon}(\gamma_t,\eta_s)
                 :=\Upsilon(\gamma_t,\eta_s)+|s-t|^2, \ \  (t,\gamma_t), (s,\eta_s)\in [0,T]\times{\Lambda}.
                 $$
 Then we can apply
                 a modification of Borwein-Preiss variational principle (see    Borwein and Zhu  \cite[Theorem 2.5.2]{bor1}) to get a maximum of a perturbation of the auxiliary function $\Psi$.
                 \par
                 Unfortunately, the second spatial derivative $\partial_{xx}S(\cdot,a_{\hat{t}})$ is not equal to $\mathbf{0}$ (see Lemma \ref{theoremS}), 
  in order to apply  Crandall-Ishii lemma (see Theorem \ref{theorem0513})  to prove comparison  theorem (see Theorem \ref{theoremhjbm}), a stronger convergence property of auxiliary functional is needed. 
 Thanks to the Step 2 in the proof of Theorem \ref{theoremhjbm}, we can find the expected convergence property of auxiliary functional  and prove the comparison theorem.
\par
%
%
  Regarding existence, we  prove that the value functional  $V$  defined in  (\ref{value1}) is  a viscosity solution to the PHJB
                         equation given in  (\ref{hjb1})   under our definition by functional It\^o  formula and dynamic programming principle.
                         \par
                           Cosso and  Russo \cite{cosso1} define a smooth functional and apply the Borwein-Preiss  variational principle to
          study
 the comparison theorem (which imply the uniqueness) of viscosity solutions for the path-dependent heat equation.  The construction of their smooth functional seems to be  
   more complicated than ours, and they only study the path-dependent heat equation.

   \par
Backward stochastic partial differential equation (BSPDE) is another interesting topic. Peng \cite{peng0} obtained the existence and uniqueness theorem for the solution to  backward stochastic HJB (BSHJB) equations in a triple. The relationship
between forward backward stochastic differential equations  and a class of semi-linear BSPDEs was established in Ma and Yong \cite{ma2}.
 For viscosity
solutions of BSPDEs, we refer to  Ekren and Zhang \cite{ekren5} and Qiu \cite{qiu, qiu1}.
As an application of our results, we give a definition of viscosity solutions to BSHJB equations, and characterize the value functional of the optimal stochastic control problem as
unique viscosity solution to the associated BSHJB equation.

\par
                 The outline of this article is  as follows. In the following
              section, we introduce the framework of \cite{cotn1} and \cite{dupire1},  preliminary results on  path-dependent stochastic optimal control problems, and a modification of Borwein-Preiss variational principle.
           In Section 3, we  present  the smooth functionals $\Upsilon^{m,M}$ which are the key to proving the stability and uniqueness results of viscosity solutions.
             In Section 4, we define classical and viscosity solutions to our
             PHJB equations and  prove  that the value functional $V$ defined by (\ref{value1}) is a viscosity solution
               to  PHJB equation  (\ref{hjb1}).   We also show
             the consistency with the notion of classical solutions and the stability result.
                 A Crandall-Ishii lemma for path-dependent case is given in Section 5.  The uniqueness of viscosity solutions for  (\ref{hjb1}) is proved in Section 6 and  Section 7 is devoted to applications to BSHJB equations.
\section{Preliminaries}  \label{RDS}
\vbox{}
2.1. \emph{Pathwise derivatives}.
 For the vectors $x,y\in \mathbb{R}^{d}$, the scalar product is denoted by $\langle x,y\rangle$ and the
         Euclidean norm $\langle x,x\rangle^{\frac{1}{2}}$ is denoted by $|x|$ (we use the same symbol $|\cdot|$ to denote the Euclidean norm on ${\mathbb{{R}}}^k$, for any $k\in {\mathbf{N}}^+$). If $A$ is a vector or matrix, its transpose is denoted by $A^\top$; For  a matrix $A$, denote its operator norm and
          Hilbert-Schmidt norm    by $|A|$ and $|A|_2$, respectively. Denote by $\mathcal{S}(\mathbb{R}^{d})$  the set of all $(d\times d)$ symmetric matrices.
Let 
$T>0$ be a fixed number.  For each  $t\in[0,T]$,
          define
         $\hat{\Lambda}_t:=D([0,t];\mathbb{R}^{d})$ as  the set of c\`adl\`ag  $\mathbb{R}^{d}$-valued
         functions on $[0,t]$.
       We denote $\hat{\Lambda}^t=\bigcup_{s\in[t,T]}\hat{\Lambda}_{s}$  and let  $\hat{\Lambda}$ denote $\hat{\Lambda}^0$.
       \par
A very important remark on the notations: as in  \cite{dupire1}, we will denote elements of $\hat{\Lambda}$ by lower
case letters and often the final time of its domain will be subscripted, e.g. $\gamma\in \hat{\Lambda}_t\subset \hat{\Lambda}$ will be
denoted by $\gamma_t$. Note that, for any $\gamma\in  \hat{\Lambda}$, there exists only one $t$ such that $\gamma\in  \hat{\Lambda}_t$. For any $0\leq s\leq t$, the value of
$\gamma_t$ at time $s$ will be denoted by $\gamma_t(s)$. Moreover, if
a path $\gamma_t$ is fixed, the path $\gamma_t|_{[0,s]}$, for $0\leq s\leq t$, will denote the restriction of the path  $\gamma_t$ to the interval
$[0,s]$.  
\par
 For convenience, define for $ x
 \in \mathbb{R}^{d},\gamma_t\in \hat{\Lambda}$, $0\leq t\leq \bar{t}\leq T$,
\begin{eqnarray*}
  \gamma^x_{t}(s)&:=&\gamma_t(s){\mathbf{1}}_{[0,t)}(s)+(\gamma_t(t)+x){\mathbf{1}}_{\{t\}}(s), \ \ s\in [0, {t}];\\
  \gamma_{t,\bar{t}}(s)&:=&\gamma_t(s){\mathbf{1}}_{[0,t)}(s)+\gamma_t(t){\mathbf{1}}_{[t,\bar{t}]}(s), \ \ s\in [0,\bar{t}].
\end{eqnarray*}
 We define a norm on  $\hat{\Lambda}_t$  and a metric on $\hat{\Lambda}$ as follows:
  for any $0\leq t\leq \bar{t}\leq T$ and $\gamma_t,\bar{\gamma}_{\bar{t}}\in \hat{\Lambda}$,
\begin{eqnarray}\label{2.1}
   ||\gamma_t||_{0}:=\sup_{0\leq s\leq t}|\gamma_t(s)|,\ \ d_{\infty}(\gamma_t,\bar{\gamma}_{\bar{t}})= d_{\infty}(\bar{\gamma}_{\bar{t}},\gamma_t):=|t-\bar{t}|
               +||\gamma_{t,\bar{t}}-\bar{\gamma}_{\bar{t}}||_{0}.
\end{eqnarray}
In the sequel,
 for notational simplicity,
we use $||\gamma_{t}-\bar{\gamma}_{\bar{t}}||_{0}$ to denote $||\gamma_{t,\bar{t}}-\bar{\gamma}_{\bar{t}}||_{0}$.
%
Then 
$(\hat{\Lambda}_t, ||\cdot||_{0})$ is a Banach space and $(\hat{\Lambda}^t, d_{\infty})$ is a complete metric space.
Following Dupire \cite{dupire1}, we define spatial derivatives of $f:\hat{\Lambda}\rightarrow \mathbb{R}$, if exist, in the standard sense: for the basis $e_i$ of $\mathbb{R}^{d}$, $i, j=1,2,\ldots,d$,
\begin{eqnarray}\label{2.2}
 \partial_{x_i}f(\gamma_s):=\lim_{h\rightarrow0}\frac{1}{h}\bigg{[}f(\gamma_s^{he_i})-f(\gamma_s)\bigg{]},\
 \partial_{x_ix_j}f:=\partial_{x_i}(\partial_{x_j}f),
    \ (s, \gamma_s)\in [0,T]\times\hat{\Lambda}, \ \
\end{eqnarray}
and the right time-derivative of $f$, if exists, as:
\begin{eqnarray}\label{2.3}
               \partial_tf(\gamma_s):=\lim_{h\rightarrow0,h>0}\frac{1}{h}\bigg{[}f(\gamma_{s,s+h})-f(\gamma_s)\bigg{]},\ (s,\gamma_s)\in [0,T)\times\hat{\Lambda}.
\end{eqnarray}
For the final time $T$, we define
$$
\partial_tf(\gamma_T):=\lim_{t<T,t\uparrow T}\partial_tf(\gamma_T|_{[0,t]}), \ \gamma_T\in \hat{\Lambda}.
$$
 We take the convention that $\gamma_s$ and  $\partial_{x}f$ denote  column vector  
    and $\partial_{xx}f$ denotes $d\times d$-matrix.
\begin{definition}\label{definitionc}
       Let $t\in[0,T)$ and $f:\hat{\Lambda}^t\rightarrow \mathbb{R}$ be given.
\begin{description}
        \item{(i)}
                 We say $f\in C^0(\hat{\Lambda}^t)$ if $f$ is continuous in $\gamma_s$  on $\hat{\Lambda}^t$ under $d_{\infty}$.
\par
       \item{(ii)}  We say $f\in C^{1,2}(\hat{\Lambda}^t)\subset C^0(\hat{\Lambda}^t)$ if  $\partial_tf$, $\partial_{x_i}f, \partial_{x_ix_j}f$ exist in $\hat{\Lambda}^t$   and are in $C^0(\hat{\Lambda}^t)$
        for all $i,j=1,2,\ldots,d$.
       \par
       \item{(iii)} We say $f\in C^{1,2}_p(\hat{\Lambda}^t)\subset C^{1,2}(\hat{\Lambda}^t)$ if $f$ and all of its derivatives grow  in a polynomial way.
\end{description}
\end{definition}
\par
 For each  $t\in[0,T]$, let $\Lambda_t:= C([0,t],\mathbb{R}^{d})$ be the set of all continuous $\mathbb{R}^{d}$-valued functions defined over $[0,t]$. We denote ${\Lambda}^t=\bigcup_{s\in[t,T]}{\Lambda}_{s}$  and let  ${\Lambda}$ denote ${\Lambda}^0$.
 Clearly, $\Lambda:=\bigcup_{t\in[0,T]}{\Lambda}_{t}\subset\hat{\Lambda}$, and each $\gamma\in \Lambda$ can also be viewed as an element of $\hat{\Lambda}$. 
 $(\Lambda_t, ||\cdot||_{0})$ is a
 Banach space, and $(\Lambda^t,d_{\infty})$ is a complete metric space.
 $f:\Lambda^t\rightarrow \mathbb{R}$ and $\hat{f}:\hat{\Lambda}^t\rightarrow \mathbb{R}$ are called consistent
  on $\Lambda^t$ if $f$ is the restriction of $\hat{f}$ on $\Lambda^t$.
\begin{definition}\label{definitionc2}
       Let  $t\in [0,T)$ and  $f:{\Lambda}^t\rightarrow \mathbb{R}$  be given.
\begin{description}
        \item{(i)}
                 We say $f\in C^0({\Lambda}^t)$ (resp., $f\in USC^0({\Lambda}^t)$, $f\in LSC^0({\Lambda}^t)$) if $f$ is continuous (resp., upper semicontinuous,  lower semicontinuous) in $\gamma_s$  on $\Lambda^t$ under $d_{\infty}$. 
\par
       \item{(ii)}  
    We say $f\in C^{1,2}_p({\Lambda}^t)$ if
        there exists $\hat{f}\in C_p^{1,2}(\hat{{\Lambda}}^t)$ which is consistent with $f$ on $\Lambda^t$.
\end{description}
\end{definition}
By Dupire \cite{dupire1} and Cont and Fournie \cite[Theorem 4.1]{cotn1}, we have the following functional It\^o formula.
\begin{theorem}\label{theoremito}
\ \ Suppose $X$ is a continuous semi-martingale and $f\in C^{1,2}_p({\hat{\Lambda}}^{\hat{t}})$ for some fixed $\hat{t}\in [0,T)$. Then for any $t\in [\hat{t},T]$:
\begin{eqnarray}\label{statesop}
                f(X_t)=f(X_{\hat{t}})+\int_{{\hat{t}}}^{t}\partial_tf(X_s)ds+\frac{1}{2}\int^{t}_{{\hat{t}}}\partial_{xx}f(X_s)d\langle X\rangle(s)
                +\int^{t}_{{\hat{t}}}\partial_xf(X_s)dX(s),\ \   P\mbox{-}a.s.
\end{eqnarray}
 Here and in the following, for every $s\in [0,T]$, $X(s)$ denotes  the value of $X$  at
 time $s$, and $X_s$ the whole history path of $X$ from time 0 to $s$.
\end{theorem}
\par
{\bf  Proof}. \ \ Define a family of functionals $F:=(F_t)_{t\in [0,T]}$ where
 $$F_t(\gamma_t,\xi_t)=f(\gamma_t),\ \ (\gamma_t,\xi_t)\in \hat{\Lambda}_t\times D([0,t], S^+_d),$$
 and $D([0,t], S^+_d)$ denotes the space of c\`adl\`ag functions on $[0,t]$, taking  values in the set $S^+_d$ of positive $d\times d$ matrices. It is clear that $F$ satisfies condtions (10) and (14) in Cont and Fournie \cite{cotn1}.
    By $f\in C^{1,2}_p({\hat{\Lambda}}^{\hat{t}})$, we have $F\in \mathbb{C}^{1,2}([\hat{t},T))$ (see Definition 3.6 in Cont and Fournie \cite{cotn1}). Thus $F$ satisfies conditions in Theorem 4.1 of Cont and Fournie \cite{cotn1}. Then, by functional It\^o formula (30) in Theorem 4.1 of Cont and Fournie \cite{cotn1}, we get (\ref{statesop}) holds true. \ \ $\Box$
\par
By the above Theorem, we have the following important result.
\begin{lemma}\label{0815lemma}
             Let $f\in C_p^{1,2}(\Lambda^t)$ and $\hat{f}\in C_p^{1,2}(\hat{\Lambda}^t)$ such that $\hat{f}$ is consistent with $f$ on $\Lambda^t$, then the following definition
             $$
             \partial_tf:=\partial_t\hat{f}, \ \ \ \partial_xf:=\partial_x\hat{f}, \ \ \ \partial_{xx}f:=\partial_{xx}\hat{f} \ \ \mbox{on} \ \Lambda^t
             $$
              is independent of the choice of $\hat{f}$. Namely, if there is another $\hat{f}'\in C_p^{1,2}(\hat{\Lambda}^t)$ such that $\hat{f}'$ is consistent with $f$ on $\Lambda^t$, then the derivatives of $\hat{f}'$
              coincide with those of $\hat{f}$ on $\Lambda^t$.
\end{lemma}
\par
{\bf  Proof}. \ \
 By the definition of the horizontal derivative, it is clear that $\partial_t\hat{f}(\gamma_l)=\partial_t\hat{f}'(\gamma_l)$ for every $\gamma_l\in \Lambda^t$.
  Next,  for every $l\in [0,T]$, define
  $$
                    \Theta_l:=\{\gamma_l\in \Lambda_l: \ \mbox{there exists a constant}\ L>0 \ \mbox{such that}\ \sup_{r_1,r_2\in [0,l]}|\gamma_l(r_1)-\gamma_l(r_2)|\leq L|r_1-r_2|\}.
  $$
  Then $\Theta_l$ is dense in $(\Lambda_l,||\cdot||_0)$.
  For every $(l,\gamma_l)\in [t,T)\times\Theta_l$ and $h\in \mathbb{R}^d$, let
  $$
  X(r)=\gamma_l(r\wedge l)+(r\vee l-l)h,\ \ r\in [0,T].
  $$
  It is clear that $X$ is a continuous semi-martingale. Then  by  Lemma \ref{theoremito},
$$
                        \int^{s}_{l}\langle\partial_x\hat{f}(X_r),h\rangle dr=\int^{s}_{l}\langle\partial_x\hat{f}'(X_r),h\rangle dr, \ \ s\in [l,T].
$$
 By the continuity of $\partial_x\hat{f}$ and $\partial_x\hat{f}'$  and the arbitrariness of $h\in {\mathbb{R}^{d}}$, we  have  $\partial_x\hat{f}(\gamma_l)=\partial_x\hat{f}'(\gamma_l)$. Notice that $\Theta_l$ is dense in $(\Lambda_l,||\cdot||_0)$,  by also the continuity of $\partial_x\hat{f}$ and $\partial_x\hat{f}'$,  we obtain that $\partial_x\hat{f}(\gamma_l)=\partial_x\hat{f}'(\gamma_l)$ for every $\gamma_l\in \Lambda^t$.
 Finally,  for every $(l,\gamma_l)\in [t,T)\times\Theta_l$ and $a\in \mathbb{R}^{d\times n}$, let
  $$
  X(r)=\gamma_l(r\wedge l)+a(W(r\vee l)-W(l)),\ \ r\in [0,T].
  $$
  It is clear that $X$ is a continuous semi-martingale. Then by  Lemma \ref{theoremito},
 $$
                        \int^{s}_{l}\mbox{tr}(\partial_{xx}\hat{f}(X_r)aa^\top)dr =\int^{s}_{l}\mbox{tr}(\partial_{xx}\hat{f}'(X_r)aa^\top)dr, \ \ s\in [l,T].
$$
 By the continuity of $\partial_{xx}\hat{f}$ and $\partial_{xx}\hat{f}'$, the arbitrariness of $a\in \mathbb{R}^{d\times n}$ and the denseness of $\Theta_l$, we also have  $\partial_{xx}\hat{f}(\gamma_l)=\partial_{xx}\hat{f}'(\gamma_l)$ for every $\gamma_l\in \Lambda^t$. \ \ $\Box$
\par
  \vbox{}
2.2. \emph{Value functional}.
           Let $\Omega:=\{\omega\in C([0,T],\mathbb{R}^n):\omega(0)={\mathbf{0}}\}$, the set of continuous functions with initial value ${\mathbf{0}}$,
           $W$ the canonical process, $P$ the Wiener measure, ${\cal {F}}$ the  Borel $\sigma$-field over $(\Omega,||\cdot||_0)$, completed with
            respect to the Wiener measure $P$ on this space.  Then  $(\Omega,{\cal {F}},P)$ is a complete   space. Here and in the sequel,
 for notational simplicity,
we use $\mathbf{0}$ to denote vectors or matrices with appropriate dimensions whose components are all equal to 0.  By $\{{\cal {F}}_t\}_{0\leq t\leq T}$ we denote  the filtration generated by $\{W(t),0\leq t\leq T\}$, augmented
       with the family $\mathcal {N}$ of $P$-null of ${\cal {F}}$.
 The filtration $\{{\cal {F}}_t\}_{0\leq t\leq T}$ satisfies the
       usual conditions.
       \par
We introduce the admissible control. Let $t,s$ be two deterministic times, $0\leq t\leq s\leq T$.
\begin{definition}
                An admissible control process $u(\cdot)=\{u(r),  r\in [t,s]\}$ on $[t,s]$  is an $\{{\cal{F}}_r\}_{t\leq r\leq s}$-progressively measurable process taking values in some Polish space $(U,d_1)$. The set of all admissible controls on $[t,s]$ is denoted by ${\cal{U}}[t,s]$. We identify two processes $u(\cdot)$ and $\tilde{u}(\cdot)$ in ${\cal{U}}[t,s]$
                and write $u(\cdot)\equiv\tilde{u}(\cdot)$ on $[t,s]$, if $P(u(\cdot)=\tilde{u}(\cdot) \ a.e. \ \mbox{in}\ [t,s])=1$.
\end{definition}
 \par
Now we consider the controlled  state equation
(\ref{state1}) and cost equation (\ref{fbsde1}). First we make the following assumption.
\begin{hyp}\label{hypstate}
$b:{\Lambda}\times U\rightarrow \mathbb{R}^{d}$, $\sigma:{\Lambda}\times U\rightarrow \mathbb{R}^{d\times n}$, $
        q: {\Lambda}\times \mathbb{R}\times \mathbb{R}^n\times U\rightarrow \mathbb{R}$ and $\phi: {\Lambda}_T\rightarrow \mathbb{{R}}$ are continuous, and
                 there exists   $L>0$
                 such that, for all $(t,\gamma_t,\eta_T,y,z,u)$,  $ (t, \gamma'_t,\eta'_T,y',z',u)
                 \in [0,T]\times\Lambda\times {\Lambda}_T\times \mathbb{R}\times \mathbb{R}^n\times U$,
      \begin{eqnarray*}
               && |b(\gamma_t,u)|^2\vee|\sigma(\gamma_t,u)|_2^2\leq
                 L^2(1+||\gamma_t||_0^2),\\
                 &&|b(\gamma_t,u)-b(\gamma'_{t},u)|\vee|\sigma(\gamma_t,u)-\sigma(\gamma'_{t},u)|_2\leq
                 L||\gamma_t-\gamma'_{t}||_0,\\
                  &&|q(\gamma_t,y,z,u)|\leq L(1+||\gamma_t||_0+|y|+|z|),
                 \\
               &&|q(\gamma_t,y,z,u)-q(\gamma'_{t},y',z',u)|\leq L(||\gamma_t-\gamma'_{t}||_0+|y-y'|+|z-z'|),
                 \\
                 &&|\phi(\eta_T)-\phi(\eta'_T)|\leq L||\eta_T-\eta'_{T}||_0.
\end{eqnarray*}
\end{hyp}
\par
The following lemma is standard; see, for example,  \cite[Lemmas 3.1 and 3.2]{tang1} (see also  El Karoui, Peng and Quenez \cite[Theorem 2.1 and Proposition 2.1]{el1}
and  Karatzas and  Shreve \cite[Theorem 3.28 in Chapter 3 and Theorem 2.9 in Chapter 5]{iakara}  for details).
\begin{lemma}\label{lemmaexist}
\ \ Assume that Hypothesis \ref{hypstate}  holds. Then for every $u(\cdot)\in {\cal{U}}[0,T]$,
$(t,\gamma_t)\in [0,T]\times {\Lambda}$ and $p\geq2$, PSDE (\ref{state1}) admits a
unique strong  solution $X^{\gamma_t,u}$, and BSDE (\ref{fbsde1}) admits a unique
 pair of solutions $(Y^{\gamma_t,u}, Z^{\gamma_t,u})$.  Furthermore, let  $X^{\gamma'_t,v}$ and $(Y^{\gamma'_t,v}, Z^{\gamma'_t,v})$ be the solutions of PSDE (\ref{state1}) and BSDE (\ref{fbsde1})
 corresponding $(t,\gamma'_t)\in [0,T]\times \Lambda$ and $v(\cdot)\in {\cal{U}}[0,T]$. Then the following estimates hold:
\begin{eqnarray}\label{fbjia1}
                \mathbb{E}[||X^{\gamma_t,u}_T-X^{\gamma_t',v}_T||_0^p|{\cal{F}}_t]&\leq& C_p||\gamma_t-\gamma_t'||_0^p
               +C_p\int^{T}_{t}{\mathbb{E}}[|b(X^{\gamma'_t,v}_l,u(l))-b(X^{\gamma_t',v}_l,v(l))|^{p}|{\cal{F}}_t]dl\nonumber\\
&&+C_p\int^{T}_{t}{\mathbb{E}}[|\sigma(X^{\gamma'_t,v}_l,u(l))-\sigma(X^{\gamma_t',v}_l,v(l))|_2^{p}|{\cal{F}}_t]dl;
\end{eqnarray}
\begin{eqnarray}\label{fbjia2}
                \mathbb{E}[||X_T^{\gamma_t,u}||_0^p|{\cal{F}}_t]\leq C_p(1+||\gamma_t||_0^p);
                \end{eqnarray}
\begin{eqnarray}\label{fbjia3}
                \mathbb{E}[||X_{r}^{\gamma_t,u}-\gamma_{t}||_0^p|{\cal{F}}_t]\leq C_p{(}1+||\gamma_t||_0^p)(r-t)^{\frac{p}{2}}, \  \ r\in [t,T];
\end{eqnarray}
and
                \begin{eqnarray}\label{fbjia4}
                 &&\mathbb{E}[||Y^{\gamma_t,u}_T-Y^{\gamma'_t,v}_T||_0^p|{\cal{F}}_t]\nonumber\\
                &\leq& C_p||\gamma_t-\gamma'_t||_0^p+C_p\int^{T}_{t}\mathbb{E}{[}{|}b(X^{\gamma'_{t},v}_l,u(l))-b(X^{\gamma'_{t},v}_l,v(l)){|}^p\nonumber\\
                &&+|\sigma(X^{\gamma'_{t},v}_l,u(l))-\sigma(X^{\gamma'_{t},v}_l,v(l))|_2^p\nonumber\\
                                                             &&+|q(X^{\gamma'_{t},v}_l,Y^{\gamma'_t,v}(l),Z^{\gamma'_t,v}(l),u(l))-q(X^{\gamma'_{t},v}_l,Y^{\gamma'_t,v}(l),Z^{\gamma'_t,v}(l),v(l))|^p|{\cal{F}}_t{]}dl;
                \end{eqnarray}
                \begin{eqnarray}\label{fbjia5}
                \mathbb{E}\left[\sup_{t\leq s\leq T}|Y^{\gamma_t,u}(s)|^p\bigg{|}{\cal{F}}_t\right]
                +\mathbb{E}\left[\bigg{(}\int^{T}_{t}|Z^{\gamma_t,u}(l)|^2dl\bigg{)}^{\frac{p}{2}}\bigg{|}{\cal{F}}_t\right]\leq C_p(1+||\gamma_t||_0^{p}).
                \end{eqnarray}
              The constant $C_p$ depending only on  $p$, $T$ and $L$.
\end{lemma}
\par
Formally,  under 
  Hypothesis \ref{hypstate}, the  value functional $V(\gamma_t)$
 defined by (\ref{value1})
 is  ${\cal{F}}_t$-measurable. 
However, according to the similar proof procedure of Proposition 3.3 in  Buckdahn and Li \cite{buck1}, we can prove the following.
\begin{theorem}\label{valuedet} \ \  Suppose the Hypothesis \ref{hypstate} holds true.
                      Then $V$ is a deterministic functional.
\end{theorem}
\par
The following  property of the value functional $V$ 
 is an immediate consequence of Lemma \ref{lemmaexist}.
\begin{lemma}\label{lemmavaluev}
\ \ Assume that Hypothesis  \ref{hypstate}     holds, then $V\in {{C}}^0(\Lambda)$ and there exists a  constant $C>0$ such that, for all  
$(t,\gamma_t,\gamma'_t)\in[0,T]\times\Lambda\times \Lambda$,
\begin{eqnarray}\label{valuelip}
|V(\gamma_t)-V(\gamma'_t)|\leq C||\gamma_t-\gamma_t'||_0; \ \ \ \ \
               |V(\gamma_t)|\leq C(1+||\gamma_t||_0).
\end{eqnarray}
\end{lemma}
\par
We now discuss a dynamic programming principle  (DPP) for the optimal control problem (\ref{state1}), (\ref{fbsde1}) and (\ref{value1}).
For this purpose, we define the family of backward semigroups associated with BSDE (\ref{fbsde1}), following the
idea of Peng \cite{peng1}.
\par
Given the initial condition $(t,\gamma_t)\in [0,T)\times{\Lambda}$, a positive number $\delta\leq T-t$, an admissible control $u(\cdot)\in {\cal{U}}[t,t+\delta]$ and
a real-valued random variable $\eta\in L^2(\Omega,{\cal{F}}_{t+\delta},P;\mathbb{R})$, we put
\begin{eqnarray}\label{gdpp}
                        G^{\gamma_t,u}_{s,t+\delta}[\eta]:=\tilde{Y}^{\gamma_t,u}(s),\ \
                        \ \ \ \ s\in[t,t+\delta],
\end{eqnarray}
                        where $(\tilde{Y}^{\gamma_t,u}(s),\tilde{Z}^{\gamma_t,u}(s))_{t\leq s\leq
                        t+\delta}$ is the solution of the following
                        BSDE with the time horizon $t+\delta$:
\begin{eqnarray}\label{bsdegdpp}
\begin{cases}
d\tilde{Y}^{\gamma_t,u}(s) =-q(X^{\gamma_t,u}_s,\tilde{Y}^{\gamma_t,u}(s),\tilde{Z}^{\gamma_t,u}(s),u(s))ds+\tilde{Z}^{\gamma_t,u}(s)dW(s),\\
 ~\tilde{Y}^{\gamma_t,u}(t+\delta)=\eta,
 \end{cases}
\end{eqnarray}
                  and $X^{\gamma_t,u}(\cdot)$ is the solution of PSDE (\ref{state1}).
                  %
%
%
\begin{theorem}\label{theoremddp} (see  \cite[Theorem 3.4]{tang1})
    Assume Hypothesis \ref{hypstate}  holds true, the value functional
                              $V$ obeys the following DPP: for
                              any $(t,\gamma_t)\in [0,T)\times{\Lambda}$ and $0<\delta\leq T-t$,
\begin{eqnarray}\label{ddpG}
                              V(\gamma_t)=\mathop{\esssup}\limits_{u(\cdot)\in{\mathcal
                              {U}}[t,t+\delta]}G^{\gamma_t,u}_{t,t+\delta}[V(X^{\gamma_t,u}_{t+\delta})].
\end{eqnarray}
\end{theorem}
 Lemma \ref{lemmavaluev} says that the value functional $V$ is
Lipschitz continuous in $\Lambda_t$. From Theorem \ref{theoremddp}, we have
\begin{theorem}\label{theorem3.9} (see    \cite[Theorem 3.7]{tang1})
                          Under Hypothesis \ref{hypstate},
                           there is a constant $C>0$ such that, for every  $0\leq t\leq t'\leq T$ and $\gamma_t, \gamma'_{t'}\in {\Lambda}$,
\begin{eqnarray}\label{hold}
                 |V(\gamma_t)-V(\gamma'_{t'})|\leq
                                    C[||\gamma_{t}-\gamma'_{t'}||_0+(1+||\gamma_t||_0)(t'-t)^{\frac{1}{2}}].
\end{eqnarray}
\end{theorem}

\vbox{}
2.3. \emph{ Borwein-Preiss variational principle}. In this subsection, we introduce a modification of Borwein-Preiss variational principle (see    Borwein and Zhu  \cite[Theorem 2.5.2]{bor1}) which
 plays a crucial role in the proof of the comparison Theorem.
We
firstly recall the definition of gauge-type function for compete  metric  space $(H,d)$.
\begin{definition}\label{gaupe}
Let $(H,d)$ be a compete metric space.
               We say that a  continuous  functional $\rho:H\times H\rightarrow [0,\infty)$ is a {gauge-type function} on the compete metric space $(H,d)$
              provided that:
             \begin{description}
        \item{(i)} $\rho(x,x)=0$ for all $x\in H$,
        \item{(ii)} for any $\varepsilon>0$, there exists $\delta>0$ such that, for all $x,y\in H$, we have $\rho(x,y)\leq \delta$ implies that
        $d(x,y)<\varepsilon$.
        \end{description}
\end{definition}
%
\begin{lemma}\label{theoremleft} 
Let $t\in [0,T]$ be fixed and let $f:\Lambda^t\rightarrow \mathbb{R}$ be an upper semicontinuous functional  bounded from above. Suppose that $\rho$ is a gauge-type function on $(\Lambda^t,d_\infty)$
 and $\{\delta_i\}_{i\geq0}$ is a sequence of positive number, and suppose that $\varepsilon>0$ and $(t_0,\gamma^0_{t_0})\in [t,T]\times \Lambda^t$ satisfy
 $$
f(\gamma^0_{t_0})\geq \sup_{(s,\gamma_s)\in [t,T]\times \Lambda^t}f(\gamma_s)-\varepsilon.
 $$
 Then there exist $(\hat{t},\hat{\gamma}_{\hat{t}})\in [t,T]\times \Lambda^t$ and a sequence $\{(t_i,\gamma^i_{t_i})\}_{i\geq1}\subset [t_0,T]\times \Lambda^t$ such that
  \begin{description}
        \item{(i)} $\rho(\gamma^0_{t_0},\hat{\gamma}_{\hat{t}})\leq \frac{\varepsilon}{\delta_0}$,  $\rho(\gamma^i_{t_i},\hat{\gamma}_{\hat{t}})\leq \frac{\varepsilon}{2^i\delta_0}$ and $t_i\uparrow \hat{t}$ as $i\rightarrow\infty$,
        \item{(ii)}  $f(\hat{\gamma}_{\hat{t}})-\sum_{i=0}^{\infty}\delta_i\rho(\gamma^i_{t_i},\hat{\gamma}_{\hat{t}})\geq f(\gamma^0_{t_0})$, and
        \item{(iii)}  $f(\gamma_s)-\sum_{i=0}^{\infty}\delta_i\rho(\gamma^i_{t_i},\gamma_s)
            <f(\hat{\gamma}_{\hat{t}})-\sum_{i=0}^{\infty}\delta_i\rho(\gamma^i_{t_i},\hat{\gamma}_{\hat{t}})$ for all $(s,\gamma_s)\in [\hat{t},T]\times \Lambda^{\hat{t}}\setminus \{(\hat{t},\hat{\gamma}_{\hat{t}})\}$.

        \end{description}
\end{lemma}
This lemma is  similar to    Borwein and Zhu  \cite[Theorem 2.5.2]{bor1}, the only difference is that here contains one more result than the latter, that is, the time series $\{t_i\}_{i\geq1}$ is monotonically increasing. For the convenience of readers, we give its proof in  Appendix A.
\par
For every $t\in [0,T]$, define $\Lambda^t\otimes\Lambda^t:=\{(\gamma_s,\eta_s)|\gamma_s,\eta_s\in \Lambda^t\}$.
It is clear that $(\Lambda^t\otimes\Lambda^t,d_{1,\infty})$ is a compete metric space,
     where $d_{1,\infty}((\gamma_t^1,\gamma_t^2),(\eta_s^1,\eta_s^2))=d_\infty(\gamma_t^1,\eta_s^1)+d_\infty(\gamma_t^2,\eta_s^2)$
       for all $(t,(\gamma_t^1,\gamma_t^2))$, $(s,(\eta_s^1,\eta_s^2))\in [0,T]\times(\Lambda\otimes\Lambda)$. Similar to Lemma \ref{theoremleft}, we have the following.

       \begin{lemma}\label{theoremleft1} 
Let $t\in [0,T]$ be fixed and let $f:\Lambda^t\otimes\Lambda^t\rightarrow \mathbb{R}$ be an upper semicontinuous functional  bounded from above. Suppose that $\rho$ is a gauge-type function on $(\Lambda^t\otimes\Lambda^t, d_{1,\infty})$
 and $\{\delta_i\}_{i\geq0}$ is a sequence of positive number, and suppose that $\varepsilon>0$ and $(t_0,(\gamma^0_{t_0},\eta^0_{t_0}))\in [t,T]\times (\Lambda^t\otimes\Lambda^t)$ satisfy
 $$
f(\gamma^0_{t_0},\eta^0_{t_0})\geq \sup_{(s,(\gamma_s,\eta_s))\in [t,T]\times (\Lambda^t\otimes\Lambda^t)}f(\gamma_s,\eta_s)-\varepsilon.
 $$
 Then there exist $(\hat{t},(\hat{\gamma}_{\hat{t}},\hat{\eta}_{\hat{t}}))\in [t,T]\times (\Lambda^t\otimes\Lambda^t)$ and a sequence $\{(t_i,(\gamma^i_{t_i},\eta^i_{t_i}))\}_{i\geq1}\subset [t_0,T]\times (\Lambda^t\otimes\Lambda^t)$ such that
  \begin{description}
        \item{(i)} $\rho((\gamma^0_{t_0},\eta^0_{t_0}),(\hat{\gamma}_{\hat{t}},\hat{\eta}_{\hat{t}}))\leq \frac{\varepsilon}{\delta_0}$,
         $\rho((\gamma^i_{t_i},\eta^i_{t_i}),(\hat{\gamma}_{\hat{t}},\hat{\eta}_{\hat{t}}))\leq \frac{\varepsilon}{2^i\delta_0}$ and $t_i\uparrow \hat{t}$ as $i\rightarrow\infty$,
        \item{(ii)}  $f(\hat{\gamma}_{\hat{t}},\hat{\eta}_{\hat{t}})-\sum_{i=0}^{\infty}\delta_i\rho((\gamma^i_{t_i},\eta^i_{t_i}),(\hat{\gamma}_{\hat{t}},\hat{\eta}_{\hat{t}}))\geq f(\gamma^0_{t_0},\eta^0_{t_0})$, and
        \item{(iii)}  $f(\gamma_s,\eta_s)-\sum_{i=0}^{\infty}\delta_i\rho((\gamma^i_{t_i},\eta^i_{t_i}),(\gamma_s,\eta_s))
            <f(\hat{\gamma}_{\hat{t}},\hat{\eta}_{\hat{t}})-\sum_{i=0}^{\infty}\delta_i\rho((\gamma^i_{t_i},\eta^i_{t_i}),(\hat{\gamma}_{\hat{t}},\hat{\eta}_{\hat{t}}))$ for all $(s,(\gamma_s,\eta_s)) $ $\in [\hat{t},T]\times (\Lambda^{\hat{t}}\otimes\Lambda^{\hat{t}})\setminus \{(\hat{t},(\hat{\gamma}_{\hat{t}},\hat{\eta}_{\hat{t}}))\}$.

        \end{description}
\end{lemma}
\section{Smooth  gauge-type functions.}
In this section we introduce the  functionals $\Upsilon^{m,M}$, which are the key to proving  the uniqueness and  stability  of viscosity solutions.
\par
For every $m\in {\mathbf{N}}^+$, define  $S_m:\hat{\Lambda}\times\hat{\Lambda}\rightarrow \mathbb{R}$ by, for every $(t,\gamma_t),(s,\eta_s)\in [0,T]\times{\hat{\Lambda}}$,
 \begin{eqnarray}\label{220817c}
S_m(\gamma_t,\eta_s)=\begin{cases}
            \frac{(||\gamma_{t}-\eta_s||_{0}^{2m}-|\gamma_{t}(t)-\eta_s(s)|^{2m})^3}{||\gamma_{t}-\eta_s||^{4m}_{0}}, \
         ~~ ||\gamma_{t}-\eta_s||_{0}\neq0; \\
0, \ ~~~~~~~~~~~~~~~~~~~~~~~~~~~~~~~~~~ ||\gamma_{t}-\eta_s||_{0}=0.
\end{cases}
\end{eqnarray}
We recall that $||\gamma_{t}-\eta_s||_{0}=||\gamma_{t,t\vee s}-\eta_{s,t\vee s}||_{0}$.
For every $m\in {\mathbf{N}}^+$ and $M\in {\mathbb{R}}$, define $\Upsilon^{m,M}$ and $\overline{\Upsilon}^{m,M}$ by
\begin{eqnarray}\label{220817c1}
          \Upsilon^{m,M}(\gamma_t,\eta_s)
          :=S_m(\gamma_t,\eta_s)+M|\gamma_t(t)-\eta_s(s)|^{2m}, \ \ \ (t,\gamma_t), (s,\eta_s)\in [0,T]\times\hat{\Lambda},
\end{eqnarray}
and
\begin{eqnarray}\label{220817c2}
        \overline{\Upsilon}^{m,M}(\gamma_t,\eta_s)
        := \Upsilon^{m,M}(\gamma_t,\eta_s)+|s-t|^2, \ \ \ (t,\gamma_t), (s,\eta_s)\in [0,T]\times\hat{\Lambda}.
\end{eqnarray}
For simplicity, we let $\Upsilon^{m,M}(\gamma_t)$ denote $\Upsilon^{m,M}(\gamma_t,\eta_t)$ when $\eta_t(l)\equiv{\mathbf{0}}$ for all $l\in [0,t]$. It is clear that $\Upsilon^{m,M}(\gamma_t,\eta_t)=\Upsilon^{m,M}(\gamma_t-\eta_t)$ for all $\gamma_t,\eta_t\in \hat{\Lambda}$. We also let $S$, $\Upsilon$ and $\overline{\Upsilon}$ denote $S_3$, $\Upsilon^{3,3}$ and $\overline{\Upsilon}^{3,3}$, respectively.\\
Now we study the   regularity  of $\Upsilon^{m,M}$ in the sense of  horizontal/vertical derivatives.
\begin{lemma}\label{theoremS}
For every fixed $m\in {\mathbf{N}}^+$, $M\in \mathbb{R}$ and $(\hat{t}, a_{\hat{t}}) \in [0,T)\times \hat{\Lambda}_{\hat{t}}$, define   $\Upsilon^{m,M}_{a_{\hat{t}}}:{\hat{\Lambda}}^{\hat{t}}\rightarrow \mathbb{R}$ by
 $$\Upsilon^{m,M}_{a_{\hat{t}}}(\gamma_t):=\Upsilon^{m,M}(\gamma_{t},a_{\hat{t}}),\ \ (t,\gamma_t)\in [\hat{t}, T]\times{\hat{\Lambda}}^{\hat{t}}.$$
 Then,  for all $(t,\gamma_t)\in [\hat{t}, T]\times{\hat{\Lambda}}^{\hat{t}}$,  $\partial_{t}\Upsilon^{m,M}_{{a_{\hat{t}}}}(\gamma_t),\partial_{x}\Upsilon^{m,M}_{{a_{\hat{t}}}}(\gamma_t),\partial_{xx}\Upsilon^{m,M}_{{a_{\hat{t}}}}(\gamma_t)$ exist and
 \begin{eqnarray}\label{220817a0}
              \partial_{t}\Upsilon^{m,M}_{{a_{\hat{t}}}}(\gamma_t)=0;
\end{eqnarray}
 \begin{eqnarray}\label{220817a}
              |\partial_{x}\Upsilon^{m,M}_{{a_{\hat{t}}}}(\gamma_t)|\leq 2m(3+|M-3|)|\gamma_t(t)-a_{\hat{t}}(\hat{t})|^{2m-1};
\end{eqnarray}
\begin{eqnarray}\label{220817a1}
              |\partial_{xx}\Upsilon^{m,M}_{{a_{\hat{t}}}}(\gamma_t)|\leq 2m[3(6m-1)+(2m-1)|M-3|]|\gamma_t(t)-a_{\hat{t}}(\hat{t})|^{2m-2}.
\end{eqnarray}
 If we also assume $m\geq2$, $\Upsilon^{m,M}_{a_{\hat{t}}}(\cdot)\in C^{1,2}_p(\hat{\Lambda}^{\hat{t}})$. 
\end{lemma}
\par
  {\bf  Proof}. \ \  First, by the definition of $\Upsilon^{m,M}_{a_{\hat{t}}}$, it is clear that $\Upsilon^{m,M}_{a_{\hat{t}}}\in C^0(\hat{\Lambda}^{\hat{t}})$ and $\partial_t\Upsilon^{m,M}_{a_{\hat{t}}}(\gamma_{t})=0$ for
   $(t,\gamma_t)\in [\hat{t},T]\times\hat{\Lambda}^{\hat{t}}$.
  Second, we consider $\partial_{x}\Upsilon^{m,M}_{a_{\hat{t}}}$.
 Define   $S_m^{a_{\hat{t}}}:{\hat{\Lambda}}^{\hat{t}}\rightarrow \mathbb{R}$ by
 $$S_m^{a_{\hat{t}}}(\gamma_t):=S_m(\gamma_{t},a_{\hat{t}}),\ \ (t,\gamma_t)\in [\hat{t}, T]\times{\hat{\Lambda}}^{\hat{t}},$$
 and define $g:{\hat{\Lambda}}^{\hat{t}}\rightarrow \mathbb{R}$ by
 $$g(\gamma_t):=|\gamma_t(t)-a_{\hat{t}}(\hat{t})|^{2m},\ \ (t,\gamma_t)\in [\hat{t}, T]\times{\hat{\Lambda}}^{\hat{t}}.$$
 It is clear that
 $$
                  \Upsilon^{m,M}_{a_{\hat{t}}}(\gamma_t)=S_m^{a_{\hat{t}}}(\gamma_t)+Mg(\gamma_t),\ \ (t,\gamma_t)\in [\hat{t}, T]\times{\hat{\Lambda}}^{\hat{t}}.
 $$
 Clearly, if $\hat{t}=0$,
   \begin{eqnarray}\label{0622a}
   \partial_{x}S_m^{a_{\hat{t}}}(\gamma_{0})=\mathbf{0},\  \ \gamma_0\in \hat{\Lambda}_0.
\end{eqnarray}
For every $(t,\gamma_t)\in (0,T]\times {\hat{\Lambda}}$, let $||\gamma_t||^{2m}_{0^-}=\sup_{0\leq s<t}|\gamma_t(s)|^{2m}$  and $(\gamma_t)_i(t)=\langle\gamma_t(t),e_i\rangle,\ i=1,2,\ldots, d$.
Then, for every $t\in[\hat{t},T]$ and $t>0$, if $|\gamma_t(t)-a_{\hat{t}}(\hat{t})|<||\gamma_t-a_{\hat{t}}||_{0-}$,
\begin{eqnarray}\label{s1}
   &&\partial_{x_i}S_m^{a_{\hat{t}}}(\gamma_{t})
   =\lim_{h\rightarrow0}\frac{S_m^{{a_{\hat{t}}}}(\gamma_t^{he_i})-S_m^{{a_{\hat{t}}}}(\gamma_{t})}{h}\nonumber\\
   &=&\lim_{h\rightarrow0}\frac{{(||\gamma_t-a_{\hat{t}}||^{2m}_{0}-|\gamma_t(t)+{he_i}-a_{\hat{t}}(\hat{t})|^{2m})^3}
   -{(||\gamma_t-a_{\hat{t}}||^{2m}_{0}-|\gamma_t(t)-a_{\hat{t}}(\hat{t})|^{2m})^3}}{h||\gamma_t-a_{\hat{t}}||^{4m}_{0}}\nonumber\\
   &=&-\frac{6m(||\gamma_t-a_{\hat{t}}||^{2m}_{0}-|\gamma_t(t)-a_{\hat{t}}(\hat{t})|^{2m})^2|\gamma_t(t)-a_{\hat{t}}(\hat{t})|^{2m-2}((\gamma_t)_i(t)-(a_{\hat{t}})_i(\hat{t}))}
   {||\gamma_t-a_{\hat{t}}||^{4m}_{0}};\ \ \ \ \
   \end{eqnarray}
 if  $|\gamma_t(t)-a_{\hat{t}}(\hat{t})|>||\gamma_t-a_{\hat{t}}||_{0^-}$,
\begin{eqnarray}\label{s2}
   \partial_{x_i}S_m^{{a_{\hat{t}}}}(\gamma_{t})=0;
   \end{eqnarray}
if  $|\gamma_t(t)-a_{\hat{t}}(\hat{t})|=||\gamma_t-a_{\hat{t}}||_{0^-}\neq0$,
since
\begin{eqnarray}\label{jiaxis}
&&||\gamma^{he_i}_t-a_{\hat{t}}||_{0}^{2m}-|\gamma_t(t)+{he_i}-a_{\hat{t}}(\hat{t})|^{2m}\nonumber
\\
&=&
\begin{cases}
0,\ \ \ \ \ \ \ \ \ \  \ \ \ \ \ \ \ \ \ \ \  \ \ \ \ \ \ \ \ \ \  \  \ \  \ \ \  \ \ \ \ ~~ ~~~~  \ \ \ ~~~~~
|\gamma_t(t)+he_i-a_{\hat{t}}(\hat{t})|\geq |\gamma_t(t)-a_{\hat{t}}(\hat{t})|;\\
  |\gamma_t(t)-a_{\hat{t}}(\hat{t})|^{2m}-|\gamma_t(t)+{h}e_i-a_{\hat{t}}(\hat{t})|^{2m},\  \ \ \ |\gamma_t(t)+he_i-a_{\hat{t}}(\hat{t})|<|\gamma_t(t)-a_{\hat{t}}(\hat{t})|,
  \end{cases}
\end{eqnarray}
we have
\begin{eqnarray}\label{s3}
 0&\leq&\lim_{h\rightarrow0}\frac{|S_m^{{a_{\hat{t}}}}(\gamma_t^{he_i})-S_m^{{a_{\hat{t}}}}(\gamma_t)|}{|h|}\nonumber\\
  &\leq&\lim_{h\rightarrow0}\frac{{
  ||\gamma_t(t)-a_{\hat{t}}(\hat{t})|^{2m}-|\gamma_t(t)+{h}e_i-a_{\hat{t}}(\hat{t})|^{2m}|^3}}{|h|\times||\gamma^{he_i}_t-a_{\hat{t}}||_{0}^{4m}} =0; \ \ \
   \end{eqnarray}
    if $|\gamma_t(t)-a_{\hat{t}}(\hat{t})|=||\gamma_t-a_{\hat{t}}||_{0^-}=0$,
\begin{eqnarray}\label{ss4}
 \partial_{x_i}S_m^{{a_{\hat{t}}}}(\gamma_{t})=0.
   \end{eqnarray}
Notice that
\begin{eqnarray}\label{2208171}
\partial_xg(\gamma_t)=2m|\gamma_t(t)-a_{\hat{t}}(\hat{t})|^{2m-2}(\gamma_t(t)-a_{\hat{t}}(\hat{t})),\ \ (t,\gamma_t)\in [\hat{t}, T]\times{\hat{\Lambda}}^{\hat{t}}.
\end{eqnarray}
From (\ref{0622a}), (\ref{s1}), (\ref{s2}), (\ref{s3}),(\ref{ss4}) and (\ref{2208171}) we obtain that,
   for all $(t,\gamma_t)\in [\hat{t}, T]\times{\hat{\Lambda}}^{\hat{t}}$,
\begin{eqnarray}\label{0528a}
    &&\partial_{x}\Upsilon^{m,M}_{{a_{\hat{t}}}}(\gamma_t)=\partial_{x}S_m^{{a_{\hat{t}}}}(\gamma_t)+M\partial_{x}g(\gamma_t)\nonumber\\
    &=&
    6m\frac{||\gamma_t-a_{\hat{t}}||^{4m}_{0}-(||\gamma_t-a_{\hat{t}}||^{2m}_{0}-|\gamma_t(t)-a_{\hat{t}}(\hat{t})|^{2m})^2}
   {||\gamma_t-a_{\hat{t}}||^{4m}_{0}}|\gamma_t(t)-a_{\hat{t}}(\hat{t})|^{2m-2}\nonumber\\
   &&\times(\gamma_t(t)-a_{\hat{t}}(\hat{t}))\mathbf{1}_{\{||\gamma_t-a_{\hat{t}}||_{0}\neq0\}}+2m(M-3)|\gamma_t(t)-a_{\hat{t}}(\hat{t})|^{2m-2}(\gamma_t(t)-a_{\hat{t}}(\hat{t})).
\end{eqnarray}
From (\ref{0528a}) it follows  that
\begin{eqnarray*}
              |\partial_{x}\Upsilon^{m,M}_{{a_{\hat{t}}}}(\gamma_t)|&\leq& 6m|\gamma_t(t)-a_{\hat{t}}(\hat{t})|^{2m-1}+2m|M-3||\gamma_t(t)-a_{\hat{t}}(\hat{t})|^{2m-1}\\
              &=&2m(3+|M-3|)|\gamma_t(t)-a_{\hat{t}}(\hat{t})|^{2m-1}.
\end{eqnarray*}
That is (\ref{220817a}).
\par
We now consider $\partial_{xx}\Upsilon^{m,M}_{a_{\hat{t}}}$.  Clearly, if $\hat{t}=0$,
   \begin{eqnarray}\label{0622b}
  \partial_{xx}S_m^{a_{\hat{t}}}(\gamma_{0})=\mathbf{0},\  \ \gamma_0\in \Lambda_0.
\end{eqnarray}
For every $(t,\gamma_t)\in [\hat{t},T]\times \hat{\Lambda}$ and $t>0$, if $|\gamma_t(t)-a_{\hat{t}}(\hat{t})|<||\gamma_t-a_{\hat{t}}||_{0^-}$,
\begin{eqnarray}\label{s5}
  && \partial_{x_jx_i}S_m^{{a_{\hat{t}}}}(\gamma_t)\nonumber\\
   &=&\lim_{h\rightarrow0}\bigg{[}\frac{{-6m(||\gamma_t-a_{\hat{t}}||_{0}^{2m}-|\gamma_t(t)+he_j-a_{\hat{t}}(\hat{t})|^{2m})^2|\gamma_t(t)+he_j-a_{\hat{t}}(\hat{t})|^{2m-2}
   }
   }{h{||\gamma_t-a_{\hat{t}}||_{0}^{4m}}}\nonumber\\
   &&~~~~~~\times ((\gamma_t)_i(t)-(a_{\hat{t}})_i(\hat{t})+h{\mathbf{1}}_{\{i=j\}})\nonumber\\
   &&~~~~~~+{\frac{6m(||\gamma_t-a_{\hat{t}}||^{2m}_{0}-|\gamma_t(t)-a_{\hat{t}}(\hat{t})|^{2m})^2|\gamma_t(t)-a_{\hat{t}}(\hat{t})|^{2m-2}((\gamma_t)_i(t)-(a_{\hat{t}})_i(\hat{t}))}
   {h||\gamma_t-a_{\hat{t}}||^{4m}_{0}}}\bigg{]}\nonumber\\
   &=&\frac{24m^2(||\gamma_t-a_{\hat{t}}||_{0}^{2m}-|\gamma_t(t)-a_{\hat{t}}(\hat{t})|^{2m})|\gamma_t(t)-a_{\hat{t}}(\hat{t})|^{4m-4}((\gamma_t)_i(t)-(a_{\hat{t}})_i(\hat{t}))}
   {{||\gamma_t-a_{\hat{t}}||_{0}^{4m}}}\nonumber\\
   &&\times ((\gamma_t)_j(t)-(a_{\hat{t}})_j(\hat{t}))-\frac{12m(m-1)(||\gamma_t-a_{\hat{t}}||_{0}^{2m}-|\gamma_t(t)-a_{\hat{t}}(\hat{t})|^{2m})^2|\gamma_t(t)-a_{\hat{t}}(\hat{t})|^{2m-4}}
   {{||\gamma_t-a_{\hat{t}}||_{0}^{4m}}}\nonumber\\
   &&\times ((\gamma_t)_i(t)-(a_{\hat{t}})_i(\hat{t}))((\gamma_t)_j(t)-(a_{\hat{t}})_j(\hat{t}))\mathbf{1}_{\{m>1\}}\nonumber\\
   &&-\frac{6m(||\gamma_t-a_{\hat{t}}||_{0}^{2m}-|\gamma_t(t)-a_{\hat{t}}(\hat{t})|^{2m})^2|\gamma_t(t)-a_{\hat{t}}(\hat{t})|^{2m-2}
   {\mathbf{1}}_{\{i=j\}}}{{||\gamma_t-a_{\hat{t}}||_{0}^{4m}}};
   \end{eqnarray}
   if $|\gamma_t(t)-a_{\hat{t}}(\hat{t})|>||\gamma_t-a_{\hat{t}}||_{0^-}$,
\begin{eqnarray}\label{s4}
   \partial_{x_jx_i}S_m^{{a_{\hat{t}}}}(\gamma_t)=0;
   \end{eqnarray}
 if $|\gamma_t(t)-a_{\hat{t}}(\hat{t})|=||\gamma_t-a_{\hat{t}}||_{0^-}\neq0$, by (\ref{jiaxis}),
we have
\begin{eqnarray}\label{s6666}
 0&\leq&\lim_{h\rightarrow0}\frac{|\partial_{x_i}S_m^{{a_{\hat{t}}}}(\gamma_t^{he_j})-\partial_{x_i}S_m^{{a_{\hat{t}}}}(\gamma_t)|}{|h|}\nonumber\\
  &\leq&\lim_{h\rightarrow0}\frac{6m{(|\gamma_t(t)-a_{\hat{t}}(\hat{t})|^{2m}-|\gamma_t(t)-a_{\hat{t}}(\hat{t})+he_j|^{2m})^2|\gamma_t(t)-a_{\hat{t}}(\hat{t})+he_j|^{2m-2}}}
  {|h|\times||\gamma^{he_j}_t-a_{\hat{t}}||_{0}^{4m}} \nonumber\\
  &&\times |(\gamma_t)_i(t)-(a_{\hat{t}})_i(\hat{t})+h{\mathbf{1}}_{\{i=j\}}|=0;
   \end{eqnarray}
    if $|\gamma_t(t)-a_{\hat{t}}(\hat{t})|=||\gamma_t-a_{\hat{t}}||_{0^-}=0$,
\begin{eqnarray}\label{ss42}
\partial_{x_jx_i}S_m^{{a_{\hat{t}}}}(\gamma_t)=0.
   \end{eqnarray}
%
%
Notice that
\begin{eqnarray}\label{2208172}
\partial_{xx}g(\gamma_t)&=&4m(m-1)|\gamma_t(t)-a_{\hat{t}}(\hat{t})|^{2m-4}(\gamma_t(t)-a_{\hat{t}}(\hat{t}))(\gamma_t(t)-a_{\hat{t}}(\hat{t}))^\top\mathbf{1}_{\{m>1\}}\nonumber\\
&&+2m|\gamma_t(t)-a_{\hat{t}}(\hat{t})|^{2m-2}I,\ \ (t,\gamma_t)\in [\hat{t}, T]\times{\hat{\Lambda}}^{\hat{t}}.
\end{eqnarray}
   Combining (\ref{0622b}), (\ref{s5}), (\ref{s4}), (\ref{s6666}), (\ref{ss42}) and (\ref{2208172}), we obtain, for all $(t,\gamma_t)\in [\hat{t}, T]\times{\hat{\Lambda}}^{\hat{t}}$,
   \begin{eqnarray}\label{0528b}
    &&\partial_{xx}\Upsilon^{m,M}_{{a_{\hat{t}}}}(\gamma_t)=\partial_{xx}S_m^{{a_{\hat{t}}}}(\gamma_t)+M\partial_{xx}g(\gamma_t)\nonumber\\
    &=&\bigg{[}24m^2\frac{(||\gamma_t-a_{\hat{t}}||_{0}^{2m}-|\gamma_t(t)-a_{\hat{t}}(\hat{t})|^{2m})|\gamma_t(t)-a_{\hat{t}}(\hat{t})|^{4m-4}
    (\gamma_t(t)-a_{\hat{t}}(\hat{t}))(\gamma_t(t)-a_{\hat{t}}(\hat{t}))^\top}
   {{||\gamma_t-a_{\hat{t}}||_{0}^{4m}}}\nonumber\\
  && +12m(m-1)\frac{||\gamma_t-a_{\hat{t}}||_{0}^{4m}-(||\gamma_t-a_{\hat{t}}||_{0}^{2m}-|\gamma_t(t)-a_{\hat{t}}(\hat{t})|^{2m})^2
 }
   {{||\gamma_t-a_{\hat{t}}||_{0}^{4m}}}|\gamma_t(t)-a_{\hat{t}}(\hat{t})|^{2m-4}\nonumber\\
   &&\times  (\gamma_t(t)-a_{\hat{t}}(\hat{t}))(\gamma_t(t)-a_{\hat{t}}(\hat{t}))^\top\mathbf{1}_{\{m>1\}}\nonumber\\
&&+6m\frac{||\gamma_t-a_{\hat{t}}||_{0}^{4m}-(||\gamma_t-a_{\hat{t}}||_{0}^{2m}-|\gamma_t(t)-a_{\hat{t}}(\hat{t})|^{2m})^2}{{||\gamma_t-a_{\hat{t}}||_{0}^{4m}}}|\gamma_t(t)-a_{\hat{t}}(\hat{t})|^{2m-2}
I\bigg{]}\mathbf{1}_{\{||\gamma_t-a_{\hat{t}}||_{0}\neq0\}}\nonumber\\
&&+4m(m-1)(M-3)|\gamma_t(t)-a_{\hat{t}}(\hat{t})|^{2m-4}(\gamma_t(t)-a_{\hat{t}}(\hat{t}))(\gamma_t(t)-a_{\hat{t}}(\hat{t}))^\top\mathbf{1}_{\{m>1\}}\nonumber\\
&&+2m(M-3)|\gamma_t(t)-a_{\hat{t}}(\hat{t})|^{2m-2}I+6mI\mathbf{1}_{\{||\gamma_t-a_{\hat{t}}||_{0}=0,m=1\}}.
\end{eqnarray}
From (\ref{0528b}) it follows that
\begin{eqnarray*}
                         |\partial_{xx}\Upsilon^{m,M}_{{a_{\hat{t}}}}(\gamma_t)|&\leq&  (24m^2+12m(m-1)+6m)|\gamma_t(t)-a_{\hat{t}}(\hat{t})|^{2m-2}\\
                            &&+2m|M-3|(1+2(m-1))|\gamma_t(t)-a_{\hat{t}}(\hat{t})|^{2m-2}\\
                            &=&2m[3(6m-1)+(2m-1)|M-3|]|\gamma_t(t)-a_{\hat{t}}(\hat{t})|^{2m-2}.
\end{eqnarray*}
That is (\ref{220817a1}).
By  (\ref{0528a}), it is clear that  $\partial_{x_i}\Upsilon^{m,M}_{{a_{\hat{t}}}}\in C^0(\hat{\Lambda}^{\hat{t}})$ for all $i=1,2,\ldots,d$.
By  (\ref{0528b}), it is clear that $\partial_{x_jx_i}\Upsilon^{m,M}_{{a_{\hat{t}}}}\in C^0(\hat{\Lambda}^{\hat{t}})$ for all $i,j=1,2,\ldots,d$ and $m\geq 2$.
 Notice that
 \begin{eqnarray}\label{220817b}
                                   |\Upsilon^{m,M}_{{a_{\hat{t}}}}(\gamma_t)|\leq (|M|+1)||\gamma_t-a_{\hat{t}}||_0^{2m},\ \ (t,\gamma_t)\in [\hat{t}, T]\times{\hat{\Lambda}}^{\hat{t}}.
 \end{eqnarray}
 Then, from (\ref{220817a0}), (\ref{220817a}) and (\ref{220817a1}) we have  $\Upsilon^{m,M}_{{a_{\hat{t}}}}(\cdot)\in C^{1,2}_{p}(\hat{\Lambda}^{\hat{t}})$ for $m\geq 2$.
 The proof is now complete. \ \ $\Box$

 For every $m\in \mathbf{N}^+$ and $M\in \mathbb{R}$, define $\overline{\Upsilon}^{m,M,2}:(\Lambda\otimes \Lambda)^2\rightarrow \mathbb{R}$ as follows: for  every $(\gamma_t,\gamma_t'), (\eta_s,\eta_s')  \in \Lambda\otimes \Lambda$,
 $$
                \overline{\Upsilon}^{m,M,2}((\gamma_t,\gamma_t'), (\eta_s,\eta_s')) :=   {\Upsilon}^{m,M}(\gamma_t,\eta_s) +{\Upsilon}^{m,M}(\gamma_t',\eta_s')+|t-s|^2.
 $$
 We now collect some properties of $\Upsilon^{m,M}$, $\overline{\Upsilon}^{m,M}$ and $\overline{\Upsilon}^{m,M,2}$.
 \begin{lemma}\label{theoremS0719}
For every $m\in \mathbf{N}^+$,  $M\in \mathbb{R}$,
\begin{eqnarray}\label{s0}
                ||\gamma_t||_{0}^{2m}+(M-3)|\gamma_t(t)|^{2m}\leq   \Upsilon^{m,M}(\gamma_t)
                \leq 3||\gamma_t||_{0}^{2m}+(M-3)|\gamma_t(t)|^{2m},  (t,\gamma_t)\in [0, T]\times{\hat{\Lambda}}.
\end{eqnarray}
If we also assume $M\geq 3$, then, for every $t\in [0,T]$,  $\overline{\Upsilon}^{m,M}$ (resp., $\overline{\Upsilon}^{m,M,2}$) is  a gauge-type function on compete metric space $(\Lambda^t,d_{\infty})$ (resp., $(\Lambda^t\otimes \Lambda^t,d_{1,\infty})$).
\end{lemma}
\par
   {\bf  Proof  }. \ \  If $||\gamma_t||_0=0$, it is clear that (\ref{s0}) holds. Then we may assume that
       $||\gamma_t||_0\neq0$.
Letting $\alpha:=|\gamma_t(t)|^{2m}$ and $M=3$, we have
 \begin{eqnarray*}
     \Upsilon^{m,3}(\gamma_t)=\frac{(||\gamma_t||_0^{2m}-|\gamma_t(t)|^{2m})^3}{||\gamma_t||_0^{4m}}+3|\gamma_t(t)|^{2m}:=f(\alpha)=\frac{(||\gamma_t||_0^{2m}-\alpha)^3}{||\gamma_t||_0^{4m}}+3\alpha.
\end{eqnarray*}
By
\begin{eqnarray}\label{03130}
              f'(\alpha)=-3\frac{(||\gamma_t||_0^{2m}-\alpha)^2}{||\gamma_t||_0^{4m}}+3\geq0, \ \ \mbox{for all}\   0\leq \alpha\leq ||\gamma_t||_0^{2m},
\end{eqnarray}
we get that
$$
                                 ||\gamma_t||_0^{2m} =f(0)\leq  \Upsilon^{m,3}(\gamma_t)=f(\alpha)\leq f(||\gamma_t||_0^{2m})=3||\gamma_t||_0^{2m},\ \ (t,\gamma_t)\in [0,T]\times\hat{\Lambda}.
$$
Notice that ${\Upsilon}^{m,M}(\gamma_t)={\Upsilon}^{m,3}(\gamma_t)+(M-3)|\gamma_t(t)|^{2m}$, we get (\ref{s0}).
\par
 Let $M\geq3$. It follows from (\ref{s0}) that, for every    $(t,(\gamma_t,\gamma_t')), (s,(\eta_s,\eta_s'))\in  [0,T]\times (\Lambda\otimes \Lambda)$,
 \begin{eqnarray}\label{0612d}
         \overline{\Upsilon}^{m,M}(\gamma_t,\eta_s)={\Upsilon}^{m,M}(\gamma_{t,t\vee s}-\eta_{s, t\vee s})+|s-t|^2
         \geq ||\gamma_t-\eta_s||_0^{2m}+|s-t|^2,
\end{eqnarray}
\begin{eqnarray}\label{0612d0719}
         \overline{\Upsilon}^{m,M,2}((\gamma_t,\gamma_t'),(\eta_s,\eta_s'))&=&{\Upsilon}^{m,M}(\gamma_t,\eta_s) +{\Upsilon}^{m,M}(\gamma_t',\eta_s')+|t-s|^2\nonumber\\
         &\geq& ||\gamma_t-\eta_s||_0^{2m}+||\gamma_t'-\eta_s'||_0^{2m}+|s-t|^2.
\end{eqnarray}
Recalling that $d_{\infty}(\gamma_t,\eta_s)=|t-s|
               +||\gamma_{t}-\eta_s||_{0}$ and $d_{1,\infty}((\gamma_t,\gamma_t'),(\eta_s,\eta_s'))=d_{\infty}(\gamma_t,\eta_s)+d_{\infty}(\gamma_t',\eta_s')$,
we get that  $\overline{\Upsilon}^{m,M}$ (resp., $\overline{\Upsilon}^{m,M,2}$) is  a gauge-type function on compete metric space $(\Lambda^t,d_{\infty})$ (resp., $(\Lambda^t\otimes \Lambda^t,d_{1,\infty})$).
  \ \ $\Box$
 \par
\begin{remark}\label{remarks}
 For every fixed  $(\hat{t}, a_{\hat{t}}) \in [0,T)\times \hat{\Lambda}_{\hat{t}}$, since $||\gamma_t-a_{\hat{t}}||_{0}^6$ does not belong to $C_p^{1,2}(\hat{\Lambda}^{\hat{t}})$, it cannot  appear  as an auxiliary functional in the proof of the
    uniqueness and stability of viscosity solutions. However,     by the above Lemmas \ref{theoremS} and \ref{theoremS0719}, we can replace $||\gamma_t-a_{\hat{t}}||_{0}^6$
   with the equivalent functional  ${\Upsilon}(\gamma_t,a_{\hat{t}})$.
\end{remark}

\par
In the proof of uniqueness of viscosity solutions,  
  we also need the following lemma.
\begin{lemma}\label{theoremS000} For $m\in {\mathbf{N}}^+$ and $M\geq3$,  we have
\begin{eqnarray}\label{up}
\left(\Upsilon^{m,M}(\gamma_t+\eta_t)\right)^{\frac{1}{2m}}\leq \left(\Upsilon^{m,M}(\gamma_t)\right)^{\frac{1}{2m}}+ \left(\Upsilon^{m,M}(\eta_t)\right)^{\frac{1}{2m}}, \ \ (t,\gamma_t,\eta_t)\in [0,T]\times \hat{\Lambda}\times \hat{\Lambda}.
\end{eqnarray}
\end{lemma}
\par
{\bf  Proof}. \ \  
    If one of $||\gamma_t||_{0}$, $||\eta_t||_{0}$ and $||\gamma_t+\eta_t||_{0}$ is equal to $0$, it is clear that (\ref{up}) holds. Then we may assume that
      all of $||\gamma_t||_{0}$, $||\eta_t||_{0}$ and $||\gamma_t+\eta_t||_{0}$ are not equal to $0$.
      By the definition of $\Upsilon^{m,M}$, we get, for every $(t, \gamma_t, \eta_t)\in [0,T]\times {\hat{\Lambda}}\times {\hat{\Lambda}}$,
   \begin{eqnarray*}
   \Upsilon^{m,M}(\gamma_t+\eta_t)&=&\frac{(||\gamma_t+\eta_t||_{0}^{2m}-|\gamma_t(t)+\eta_t(t)|^{2m})^3}
   {||\gamma_t+\eta_t||_{0}^{4m}}+M|\gamma_t(t)+\eta_t(t)|^{2m}\\
   &=&||\gamma_t+\eta_t||_{0}^{2m}-\frac{|\gamma_t(t)+\eta_t(t)|^{6m}}{||\gamma_t+\eta_t||_{0}^{4m}}+3\frac{|\gamma_t(t)+\eta_t(t)|^{4m}}{||\gamma_t+\eta_t||_{0}^{2m}}+(M-3)|\gamma_t(t)+\eta_t(t)|^{2m}.
\end{eqnarray*}
Letting $x:=||\gamma_t+\eta_t||_{0}^{2m}$ and $y:=|\gamma_t(t)+\eta_t(t)|^{2m}$, we have
 \begin{eqnarray*}
   \Upsilon^{m,M}(\gamma_t+\eta_t)=f(x,y):=x-\frac{y^3}{x^2}+3\frac{y^2}{x}+(M-3)y.
\end{eqnarray*}
By $$
              f_x(x,y)=1+2\bigg{(}\frac{y}{x}\bigg{)}^3-3\bigg{(}\frac{y}{x}\bigg{)}^2=\bigg{(}\frac{2y}{x}+1\bigg{)}\bigg{(}\frac{y}{x}-1\bigg{)}^2\geq0, \ \ 0\leq y\leq x, \ x>0,
$$
$$
              f_y(x,y)=-3\frac{y^2}{x^2}+6\frac{y}{x}+(M-3)\geq0, \ \ 0\leq y\leq x,\ x>0,
$$
 $$||\gamma_t+\eta_t||_{0}\leq ||\gamma_t||_{0}+||\eta_t||_{0},
 $$
and $$|\gamma_t(t)+\eta_t(t)|\leq |\gamma_t(t)|+|\eta_t(t)|,
$$
we have
\begin{eqnarray*}
  \Upsilon^{m,M}(\gamma_t+\eta_t)&\leq&(||\gamma_t||_{0}+||\eta_t||_{0})^{2m}-\frac{(|\gamma_t(t)|+|\eta_t(t)|)^{6m}}{(||\gamma_t||_{0}+||\eta_t||_{0})^{4m}}
   +3\frac{(|\gamma_t(t)|+|\eta(t)|)^{4m}}{(||\gamma_t||_{0}+||\eta_t||_{0})^{2m}}\\
   &&+(M-3)(|\gamma_t(t)|+|\eta_t(t)|)^{2m}.
\end{eqnarray*}
Let $a=\frac{|\gamma_t(t)|}{||\gamma_t||_{0}}$, $b=\frac{|\eta_t(t)|}{||\eta_t||_{0}}$, $\alpha=||\gamma_t||_{0}$ and $\beta=||\eta_t||_{0}$,
   we get that
\begin{eqnarray*}
                      \Upsilon^{m,M}(\gamma_t+\eta_t)&\leq& (\alpha+\beta)^{2m}-\frac{(a\alpha+b\beta)^{6m}}{(\alpha+\beta)^{4m}}+3\frac{(a\alpha+b\beta)^{4m}}{(\alpha+\beta)^{2m}}+(M-3)(a\alpha+b\beta)^{2m}\\
                      &=& (\alpha+\beta)^{2m}g^{2m}\left(\frac{a\alpha+b\beta}{\alpha+\beta}\right),    \\
                                        \Upsilon^{m,M}(\gamma_t)&=& \alpha^{2m}-{\alpha^{2m}}a^{6m}+3{\alpha^{2m}}a^{4m}+(M-3)a^{2m}\alpha^{2m}
                                        =\alpha^{2m}g^{2m}(a),\\
                      \Upsilon^{m,M}(\eta_t)&=& \beta^{2m}-{\beta^{2m}}b^{6m}+3{\beta^{2m}}b^{4m}+(M-3)b^{2m}\beta^{2m}
                      =\beta^{2m}g^{2m}(b),
\end{eqnarray*}
where
\begin{eqnarray}\label{0404}
                      g(x)=\left(1-{x^{6m}}+3{x^{4m}}+(M-3)x^{2m}
                      \right)^{\frac{1}{2m}},\ \  x\in [0,1].
\end{eqnarray}
By the following Lemma \ref{theoremS0404} , we have that $g$ is a  convex function on $[0,1]$, then 
\begin{eqnarray*}
                     && \left(\Upsilon^{m,M}(\gamma_t+\eta_t)\right)^{\frac{1}{2m}}-\left(\Upsilon^{m,M}(\gamma_t)\right)^{\frac{1}{2m}}-\left(\Upsilon^{m,M}(\eta_t)\right)^{\frac{1}{2m}} \\
                     &\leq&(\alpha+\beta)g\left(\frac{a\alpha+b\beta}{\alpha+\beta}\right)-\alpha g(a)-\beta g(b)\\
                    & =&(\alpha+\beta)\left(g\left(\frac{\alpha}{\alpha+\beta}a+\frac{\beta}{\alpha+\beta}b\right)-\frac{\alpha}{\alpha+\beta} g(a)-\frac{\beta}{\alpha+\beta} g(b)\right)\leq0.
\end{eqnarray*}
Thus we obtain (\ref{up}).
The proof is now complete. \ \ $\Box$
\par
 To complete the previous proof, it remains to state and prove the following lemma.
\begin{lemma}\label{theoremS0404} For $m\in {\mathbf{N}}^+$ and $M\geq3$, the function $g$ defined by (\ref{0404}) is a convex function on $[0,1]$.
\end{lemma}
\par
{\bf  Proof}. \ \
   By the definition of $g$, for all $x\in [0,1]$,
   \begin{eqnarray*}
                       g'(x)=\frac{1}{2m}g^{1-2m}(x)(-6mx^{6m-1}+12mx^{4m-1}+2m(M-3)x^{2m-1}
                       ),
   \end{eqnarray*}
   and
   \begin{eqnarray*}
                       g''(x)&=&\frac{1}{2m}g^{1-2m}(x)(-6m(6m-1)x^{6m-2}+12m(4m-1)x^{4m-2}+2m(2m-1)(M-3)x^{2m-2}
                       )\\
                       &&-\frac{1}{2m}(1-\frac{1}{2m})g^{1-4m}(x)(-6mx^{6m-1}+12mx^{4m-1}+2m(M-3)x^{2m-1}
                       )^2\\
                       &=& 3g^{1-4m}(x)x^{4m-2}(2x^{8m}-(2m+7)x^{6m}+6x^{4m}-(6m-1)x^{2m}+8m-2)\\
                        &&+g^{1-4m}(x)x^{2m-2}(M-3)(-(8m+2)x^{6m}+(6m+3)x^{4m}+2m-1).
   \end{eqnarray*}
   Notice that
    \begin{eqnarray*}
   &&2x^{8m}-(2m+7)x^{6m}+6x^{4m}-(6m-1)x^{2m}+8m-2\\
   &\geq& 2x^{8m}-(2m+7)x^{6m}+6x^{4m}-(6m-1)x^{2m}+(2m-1)x^{4m}+6m-1\\
   &=&x^{4m}(1-x^{2m})(2m+5-2x^{2m})+(6m-1)(1-x^{2m})\geq0, \  \ x\in [0,1],
     \end{eqnarray*}
     and
   \begin{eqnarray*}
                  &&-(8m+2)x^{6m}+(6m+3)x^{4m}+2m-1\geq -(8m+2)x^{6m}+(6m+3)x^{4m}+(2m-1)x^{4m}\\
                  &=&(8m+2)x^{4m}(1-x^{2m})\geq0, \  \ x\in [0,1],
   \end{eqnarray*}
 we get that   $g''(x)\geq0$ for all $x\in [0,1]$ and  the function $g$  is a convex function on $[0,1]$.
The proof is now complete. \ \ $\Box$
\section{Viscosity solutions to  PHJB equations: Existence theorem.}

\par
In this section, we consider the  second order path-dependent  Hamilton-Jacobi-Bellman
                   (PHJB) equation (\ref{hjb1}). As usual, we start with classical solutions.
\par
\begin{definition}\label{definitionccc}     (Classical solution)
              A functional $v\in C_p^{1,2}({\Lambda})$       is called a classical solution (resp., subsolution, supersolution) to the PHJB equation (\ref{hjb1}) if the terminal condition,  $v(\gamma_T)=(\mbox{resp.},  \ \leq,  \geq) \phi(\gamma_T)$ for all
                             $\gamma_T\in {\Lambda}_T$ is satisfied, and   ${\mathcal{L}}v(\gamma_t)=(\mbox{resp.}, \ \geq, \ \leq) 0$ for all $(t,\gamma_t)\in [0,T)\times \Lambda$.
 \end{definition}

        \par
        We will prove that the value functional $V$ defined by (\ref{value1}) is a viscosity solution of PHJB equation (\ref{hjb1}). We  give the following definition for the viscosity solutions.
                      For every $(t,\gamma_t)\in [0,T]\times \Lambda$ and $w\in USC^0(\Lambda)$, define
$$
             \mathcal{A}^+(\gamma_t,w):=\bigg{\{}\varphi\in C^{1,2}_p({\Lambda}^{t}):  0=({w}-{{\varphi}})({\gamma_t})=\sup_{(s,\eta_s)\in [t,T]\times{\Lambda}}
                         ({w}- {{\varphi}})(
                         \eta_s)\bigg{\}},
$$
and, for every $(t,\gamma_t)\in [0,T]\times \Lambda$ and $w\in LSC^0(\Lambda)$, define
$$
             \mathcal{A}^-(\gamma_t,w):=\bigg{\{}\varphi\in C^{1,2}_p({\Lambda}^{t}):  0=({w}-{{\varphi}})({\gamma_t})=\inf_{(s,\eta_s)\in [t,T]\times{\Lambda}}
                         ({w}-{{\varphi}})(
                         \eta_s)\bigg{\}}.
$$
\begin{definition}\label{definition4.1} \ \
 $w\in USC^0({\Lambda})$ (resp., $w\in LSC^0({\Lambda})$) is called a
                             viscosity subsolution (resp.,  supersolution)
                             to  (\ref{hjb1}) if the terminal condition,  $w(\gamma_T)\leq \phi(\gamma_T)$(resp.,  $w(\gamma_T)\geq \phi(\gamma_T)$) for all
                             $\gamma_T\in {\Lambda}_T$ is satisfied,
                                and whenever  ${\varphi}\in \mathcal{A}^+(\gamma_s,w)$ (resp.,  ${\varphi}\in \mathcal{A}^-(\gamma_s,w)$)  with $(s,{\gamma}_{s})\in [0,T)\times{\Lambda}$,  we have
\begin{eqnarray*}
                     {\mathcal{L}}{\varphi}(\gamma_s) \geq0 \ \ (\mbox{resp.},\ {\mathcal{L}}{\varphi}(\gamma_s) \leq0).
\end{eqnarray*}
                                $w\in C^0({\Lambda})$ is said to be a
                             viscosity solution to equation (\ref{hjb1}) if it is
                             both a viscosity subsolution and a viscosity
                             supersolution.
\end{definition}
\begin{remark}\label{remarkv}
  Assume that  the coefficients $b(\gamma_t,u)=\overline{b}(t,\gamma_t(t),u),
                     \ \sigma(\gamma_t,u)=\overline{\sigma}(t,\gamma_t(t),u)$,  $ q(\gamma_t,$ $y,z,u)=\overline{q}(t,\gamma_t(t),y,z,u),
                     \phi(\eta_T)=\overline{\phi}(\eta_T(T))$ for all
                     $(t,\gamma_t,y,z,u) \in [0,T]\times{\Lambda}\times \mathbb{R}\times \mathbb{R}^n\times U$ and $\eta_T\in \Lambda_T$. Then there exists a function  $\overline{V}: [0,T]\times \mathbb{R}^{d}\rightarrow \mathbb{R}$ such that  $V(\gamma_t)=\overline{V}(t,\gamma_t(t))$ for all
                     $(t,\gamma_t) \in [0,T]\times \Lambda$, and
                       PHJB equation (\ref{hjb1}) reduces to the following HJB equation:
 \begin{eqnarray}\label{hjb3}
\begin{cases}
\overline{V}_{t^+}(t,x)+\overline{{\mathbf{H}}}(t,x,\overline{V}(t,x),\nabla_x\overline{V}(t,x),\nabla^2_{x}\overline{V}(t,x))= 0,\ \ \  (t, x)\in
                               [0,T)\times \mathbb{R}^{d},\\
 \overline{V}(T,x)=\overline{\phi}(x), \ \ \ x\in \mathbb{R}^{d};
 \end{cases}
\end{eqnarray}
where
 \begin{eqnarray*}
                                \overline{{\mathbf{H}}}(t,x,r,p,\iota)&=&\sup_{u\in{
                                         {U}}}[\langle p,\overline{b}(t,x,u)\rangle +\frac{1}{2}\mbox{tr}[ \iota \overline{\sigma}(t,x,u)\overline{\sigma}^\top(t,x,u)]\\
                        &&\ \ \ \ \   +\overline{q}(t,x,r,\overline{\sigma}^\top(t,x,u)p,u)], \ \ \ (t,x,r,p,\iota)\in [0,T]\times \mathbb{R}^{d}\times \mathbb{R}\times \mathbb{R}^{d}\times \mathcal{S}(\mathbb{R}^{d}).
\end{eqnarray*}
Here  and in the sequel, $\nabla_x$ and $\nabla^2_{x}$ denote the standard  first and second order derivatives
with respect to $x$. However, slightly different from the HJB literature, $\overline{V}_{t^+}$ denotes the right time-derivative of $\overline{V}$.
 %
%
%
%
\end{remark}
 The following theorem show that our definition of viscosity solutions to PHJB equation (\ref{hjb1}) is a natural  extension of classical viscosity solution to HJB equation  (\ref{hjb3}).
 \begin{theorem}\label{theoremnatural} \ \ Consider the setting in Remark \ref{remarkv}.
                          Assume that $V$ is a viscosity solution (resp., subsolution, supersolution) of PHJB equation (\ref{hjb1}) in the sense of Definition \ref{definition4.1}. Then $\overline{V}$ is a
                          viscosity solution (resp., subsolution, supersolution) of HJB equation (\ref{hjb3}) in the standard sense (see Definition 5.1 on  page 190 of \cite{yong}).
\end{theorem}
 {\bf  Proof}. \ \ Without loss of generality, we shall only prove the viscosity subsolution property.
First,  from $V$ is a viscosity subsolution of equation (\ref{hjb1}), it follows that, for every $x\in \mathbb{R}^{d}$,
 $$
                            \overline{V}(T,x)=V(\gamma_T)\leq \phi(\gamma_T)=\overline{\phi}(x),
$$
where $\gamma_T\in \Lambda$ with $\gamma_T(T)= x$.
\\
Next, let $\overline{\varphi}\in C^{1,2}([0,T]\times \mathbb{R}^{d})$ and $(t,x)\in [0,T)\times \mathbb{R}^{d}$
 such that $$
                         0=(\overline{V}- {\overline{\varphi}})(t,x)=\sup_{(s,y)\in [0,T]\times \mathbb{R}^{d}}
                         (\overline{V}- {\overline{\varphi}})(
                         s,y).
$$
We can modify $\overline{\varphi}$  such that $\overline{\varphi}$, $\overline{\varphi}_t$, $\nabla_{x}\overline{\varphi}$ and $\nabla^2_{x}\overline{\varphi}$  grow  in a polynomial way. Here  and in the sequel,
$\overline{\varphi}_t$ denotes the  time-derivative of $\overline{\varphi}$.
%
%
 Define $\varphi:\hat{\Lambda}\rightarrow \mathbb{R}$ by
 $$
 \varphi(\gamma_s)=\overline{\varphi}(s,\gamma_s(s)),\  (s, \gamma_s)\in [0,T]\times\hat{\Lambda},
 $$
 and define $\hat{\gamma}_{t}\in \Lambda_t$ by
 $$
   \hat{\gamma}_{t}(s)= x,\ \ s\in [0,t].$$
 It is clear that,
 $$ \partial_x\varphi(\gamma_s)=\nabla_x\overline{\varphi}(s,\gamma_s(s)),\ \
 \partial_{xx}\varphi(\gamma_s)=\nabla^2_{x}\overline{\varphi}(s,\gamma_s(s)),\ \ (s,\gamma_s)\in [0,T]\times\hat{\Lambda},
 $$
 $$\partial_t\varphi(\gamma_s)=\overline{\varphi}_{t^+}(s,\gamma_s(s)),\ \ (s,\gamma_s)\in [0,T)\times\hat{\Lambda},
 $$
 and
 $$
 \partial_t\varphi(\gamma_T)=\lim_{s<T,s\uparrow T}\partial_t\varphi(\gamma_T|_{[0,s]})=\lim_{s<T,s\uparrow T}\overline{\varphi}_{t^+}(s,\gamma_T(s))=\overline{\varphi}_{t^+}(T,\lim_{s<T,s\uparrow T}\gamma_T(s)),\
 \gamma_T\in \hat{\Lambda}_T.
 $$
 Thus we have $\varphi\in C^{1,2}_p(\Lambda)\subset C^{1,2}_p(\Lambda^{t})$. Moreover, by the definitions of $V$ and $\varphi$,
 $$
 0=(V-{{\varphi}})(\hat{\gamma}_{t})=(\overline{V}-\overline{\varphi})(t,x)=\sup_{(s,y)\in [0,T]\times \mathbb{R}^{d}}
                         (\overline{V}- {\overline{\varphi}})(
                         s,y)=\sup_{(s,\gamma_s)\in [t,T]\times\Lambda}
                         (V-{{\varphi}})(
                         \gamma_s).
 $$
 Therefore,    $ \varphi \in \mathcal{A}^+(\hat{\gamma}_{t},V)$ with $(t,\hat{\gamma}_{t})\in [0,T)\times \Lambda$.  
 Since $V$ is a viscosity subsolution of PHJB equation (\ref{hjb1}), we have
 $$
{\mathcal{L}}{\varphi}(\hat{\gamma}_{t}) \geq0.
$$
Thus,
$$\overline{\varphi}_{t^+}(t,x)+\overline{{\mathbf{H}}}(t,x,\overline{\varphi}(t,x),\nabla_x\overline{\varphi}(t,x),\nabla^2_{x}\overline{\varphi}(t,x))\geq 0.$$
 By the arbitrariness of   $\overline{\varphi}\in C^{1,2}([0,T]\times \mathbb{R}^{d})$, we see that $\overline{V}$ is a viscosity subsolution of HJB equation (\ref{hjb3}), and thus completes the proof. \ \ $\Box$
 \par
  We note that  the viceversa of Theorem \ref{theoremnatural} does not hold true. Indeed, let $ \varphi \in \mathcal{A}^+(\hat{\gamma}_{t},V)$ with $(t,\hat{\gamma}_{t})\in [0,T)\times \Lambda$. Since $\varphi$
   belongs only to $C^{1,2}_p(\Lambda^t)$ and not to $C^{1,2}_p(\Lambda)$, we cannot construct $\overline{\varphi}\in C^{1,2}([0,T]\times \mathbb{R}^{d})$ from $\varphi$
 such that $\overline{\varphi}_{t}(t,\hat{\gamma}_t(t))=\partial_t\varphi(\hat{\gamma}_t),\ \nabla_x\overline{\varphi}(t,\hat{\gamma}_t(t))=\partial_x\varphi(\hat{\gamma}_t), \
 \nabla^2_{x}\overline{\varphi}(t,\hat{\gamma}_t(t))=\partial_{xx}\varphi(\hat{\gamma}_t)
 $ and
 $$
                         0=(\overline{V}- {\overline{\varphi}})(t,\hat{\gamma}_{t}(t))=\sup_{(s,y)\in [0,T]\times \mathbb{R}^{d}}
                         (\overline{V}- {\overline{\varphi}})(
                         s,y).
$$
\par
We are now in a  position  to give  the existence and consistency results 
for the viscosity solutions.
\begin{theorem}\label{theoremvexist} \ \
                          Suppose that Hypothesis \ref{hypstate}  holds. Then the value
                          functional $V$ defined by (\ref{value1}) is a
                          viscosity solution to equation  (\ref{hjb1}).
\end{theorem}
\begin{theorem}\label{theorem3.2}
                      Let Hypothesis \ref{hypstate} hold true, $v\in C_p^{1,2}({\Lambda})$. Then
                 $v$ is a classical solution (resp., subsolution, supersolution) of  equation (\ref{hjb1}) if and only if it is a viscosity solution (resp., subsolution, supersolution).
\end{theorem}
The proof of Theorems  \ref{theoremvexist} and \ref{theorem3.2} is rather standard. Moreover, note that a viscosity solution in the sense of \cite{ekren3} is a viscosity solution  in our sense, then these results  can be  implied by
\cite{ekren3} directly. However, our conditions are weaker than those in  \cite{ekren3}. For the sake of the completeness of the article and the convenience of readers, we give their proof in the appendix B.
               \par
We conclude this section with   the stability of viscosity solutions.
\begin{theorem}\label{theoremstability}
                      Let $b,\sigma,q,\phi$ satisfy  Hypothesis \ref{hypstate}, and $v\in C^0(\Lambda)$ (resp., $v\in USC^0(\Lambda)$, $v\in LSC^0(\Lambda)$).  Assume
                       \item{(i)}      for any $\varepsilon>0$, there exist $b^\varepsilon, \sigma^\varepsilon, q^\varepsilon, \phi^\varepsilon$ and $v^\varepsilon\in C^0(\Lambda)$ (resp., $v^\varepsilon\in USC^0(\Lambda)$, $v^\varepsilon\in LSC^0(\Lambda)$) such that  $b^\varepsilon, \sigma^\varepsilon, q^\varepsilon, \phi^\varepsilon$ satisfy  Hypothesis \ref{hypstate} and $v^\varepsilon$ is a viscosity solution (resp., subsolution,  supersolution) of PHJB equation (\ref{hjb1}) with generators $b^\varepsilon, \sigma^\varepsilon, q^\varepsilon, \phi^\varepsilon$;
                           \item{(ii)} as $\varepsilon\rightarrow0$, $(b^\varepsilon, \sigma^\varepsilon, q^\varepsilon, \phi^\varepsilon,v^\varepsilon)$ converge to
                           $(b, \sigma, q, \phi, v)$  uniformly in the following sense: 
\begin{eqnarray}\label{sss}
                          \lim_{\varepsilon\rightarrow0}&&\sup_{(t,\gamma_t,x,y,u)\in [0,T]\times \Lambda \times \mathbb{R}\times \mathbb{R}^{d}\times U}\sup_{\eta_T\in  \Lambda_T}[(|b^\varepsilon-b|+|\sigma^\varepsilon-\sigma|_2)({\gamma}_t,u)
                         \nonumber\\
                         &&~~~~~~~~~~~~+|q^\varepsilon-q|({\gamma}_t,{x},\sigma^\top({\gamma}_{{t}},u){y},u)+|\phi^\varepsilon-\phi|(\eta_T)+|v^\varepsilon-v|({\gamma}_t)]=0.
\end{eqnarray}
                   Then $v$ is a viscosity solution (resp., subsolution, supersolution) of PHJB equation (\ref{hjb1}) with generators $b,\sigma,q,\phi$.
\end{theorem}
{\bf  Proof}. \ \ Without loss of generality, we shall only prove the viscosity subsolution property.
First,  from $v^{\varepsilon}$ is a viscosity subsolution of  equation (\ref{hjb1}) with generators $b^{\varepsilon}, \sigma^{\varepsilon},  q^{\varepsilon}, \phi^{\varepsilon}$, it follows that
 $$
                            v^{\varepsilon}(\gamma_T)\leq \phi^{\varepsilon}(\gamma_T),\ \ \gamma_T\in \Lambda_T.
$$
Letting $\varepsilon\rightarrow0$, we  have
 $$
                            v(\gamma_T)\leq \phi(\gamma_T),\ \ \gamma_T\in \Lambda_T.
$$
Next,     we let   $\varphi\in \mathcal{A}^+(\hat{\gamma}_{\hat{t}}, v)$ with
  $(\hat{t},\hat{\gamma}_{\hat{t}})\in [0,T)\times\Lambda$. By (\ref{sss}), there exists a constant $\delta>0$ such that for all $\varepsilon\in (0,\delta)$,
  $$\sup_{(t,\gamma_t)\in [\hat{t},T]\times\Lambda^{\hat{t}}}(v^\varepsilon(\gamma_t)-\varphi(\gamma_t))\leq 1.$$
 Denote $\varphi_{1}(\gamma_t):=\varphi(\gamma_t)+\overline{\Upsilon}(\gamma_t,\hat{\gamma}_{{\hat{t}}})$
 for all
 $(t,\gamma_t)\in [0,T]\times\Lambda$. By Lemma \ref{theoremS}, we have  $\varphi_{1}\in C^{1,2}_{p}(\Lambda^{\hat{t}})$.
 For every  $\varepsilon\in (0,\delta)$, it is clear that $v^{\varepsilon}-{{\varphi_{1}}}$ is an  upper semicontinuous functional
  and bounded from above on $\Lambda^{\hat{t}}$.  Define a sequence of positive numbers $\{\delta_i\}_{i\geq0}$  by 
        $\delta_i=\frac{1}{2^i}$ for all $i\geq0$. 
        Since 
         $\overline{\Upsilon}(\cdot,\cdot)$ is a gauge-type function on $(\Lambda^{\hat{t}}, d_{\infty})$, from Lemma \ref{theoremleft} it follows that,
 for every  $(t_0,\gamma^0_{t_0})\in [\hat{t},T]\times \Lambda^{\hat{t}}$ satisfying
 $$
(v^{\varepsilon}-{{\varphi_{1}}})(\gamma^0_{t_0})\geq \sup_{(s,\gamma_s)\in [\hat{t},T]\times \Lambda^{\hat{t}}}(v^{\varepsilon}-{{\varphi_{1}}})(\gamma_s)-\varepsilon,\
\    \mbox{and} \ \ (v^{\varepsilon}-{{\varphi_{1}}})(\gamma^0_{t_0})\geq (v^{\varepsilon}-{{\varphi_{1}}})(\hat{\gamma}_{\hat{t}}),
 $$
  there exist $(t_{\varepsilon},{\gamma}^{\varepsilon}_{t_{\varepsilon}})\in [\hat{t},T]\times \Lambda^{\hat{t}}$ and a sequence $\{(t_i,\gamma^i_{t_i})\}_{i\geq1}\subset [t_0,T]\times \Lambda^{\hat{t}}$ such that
  \begin{description}
        \item{(i)} $\overline{\Upsilon}(\gamma^0_{t_0},{\gamma}^{\varepsilon}_{t_{\varepsilon}})\leq {\varepsilon}$,  $\overline{\Upsilon}(\gamma^i_{t_i},{\gamma}^{\varepsilon}_{t_{\varepsilon}})\leq \frac{\varepsilon}{2^i}$ and $t_i\uparrow t_{\varepsilon}$ as $i\rightarrow\infty$,
        \item{(ii)}  $(v^{\varepsilon}-{{\varphi_{1}}})({\gamma}^{\varepsilon}_{t_{\varepsilon}})-\sum_{i=0}^{\infty}\frac{1}{2^i}\overline{\Upsilon}(\gamma^i_{t_i},{\gamma}^{\varepsilon}_{t_{\varepsilon}})\geq (v^{\varepsilon}-{{\varphi_{1}}})(\gamma^0_{t_0})$, and
        \item{(iii)}  $(v^{\varepsilon}-{{\varphi_{1}}})(\gamma_s)-\sum_{i=0}^{\infty}\frac{1}{2^i}\overline{\Upsilon}(\gamma^i_{t_i},\gamma_s)
            <(v^{\varepsilon}-{{\varphi_{1}}})({\gamma}^{\varepsilon}_{t_{\varepsilon}})-\sum_{i=0}^{\infty}\frac{1}{2^i}\overline{\Upsilon}(\gamma^i_{t_i},{\gamma}^{\varepsilon}_{t_{\varepsilon}})$ for all $(s,\gamma_s)\in [t_{\varepsilon},T]\times \Lambda^{t_{\varepsilon}}\setminus \{(t_{\varepsilon},{\gamma}^{\varepsilon}_{t_{\varepsilon}})\}$.

        \end{description}
  %
%
%
%
We claim that
\begin{eqnarray}\label{gamma}
d_\infty({\gamma}^{\varepsilon}_{{t}_{\varepsilon}},\hat{\gamma}_{\hat{t}})\rightarrow0  \ \ \mbox{as} \ \ \varepsilon\rightarrow0.
\end{eqnarray}
 Indeed, if not,  by (\ref{0612d}) and the definition of $d_{\infty}$, we can assume that there exists a constant  $\nu_0>0$
 such
                    that
$$
  \overline{\Upsilon}({\gamma}^{\varepsilon}_{{t}_{\varepsilon}},\hat{\gamma}_{{\hat{t}}})
  \geq\nu_0.
$$
 Thus, by  the property (ii) of $(t_{\varepsilon},{\gamma}^{\varepsilon}_{t_{\varepsilon}})$,  we obtain that
\begin{eqnarray*}
   &&0=(v- {{\varphi}})(\hat{\gamma}_{\hat{t}})= \lim_{\varepsilon\rightarrow0}(v^\varepsilon-{{\varphi_{1}}})(\hat{\gamma}_{\hat{t}})
   \leq {\limsup_{\varepsilon\rightarrow0}}\bigg{[}(v^{\varepsilon}-{{\varphi_{1}}})({\gamma}^{\varepsilon}_{t_{\varepsilon}})-\sum_{i=0}^{\infty}\frac{1}{2^i}\overline{\Upsilon}(\gamma^i_{t_i},{\gamma}^{\varepsilon}_{t_{\varepsilon}})\bigg{]}\\
   &=&{\limsup_{\varepsilon\rightarrow0}}\bigg{[}(v^\varepsilon-{{\varphi}})({\gamma}^{\varepsilon}_{{t}_{\varepsilon}})
  -\overline{\Upsilon}({\gamma}^{\varepsilon}_{{t}_{\varepsilon}},\hat{\gamma}_{{\hat{t}}})-
  \sum_{i=0}^{\infty}\frac{1}{2^i}\overline{\Upsilon}(\gamma^i_{t_i},{\gamma}^{\varepsilon}_{t_{\varepsilon}})\bigg{]}\\
   &\leq&{\limsup_{\varepsilon\rightarrow0}}\bigg{[}(v-{{\varphi}})({\gamma}^{\varepsilon}_{{t}_{\varepsilon}})+(v^\varepsilon-v)({\gamma}^{\varepsilon}_{{t}_{\varepsilon}})
     -\sum_{i=0}^{\infty}\frac{1}{2^i}\overline{\Upsilon}(\gamma^i_{t_i},{\gamma}^{\varepsilon}_{t_{\varepsilon}})\bigg{]}-\nu_0\leq (v- {{\varphi}})(\hat{\gamma}_{\hat{t}})-\nu_0=-\nu_0,
\end{eqnarray*}
 contradicting $\nu_0>0$.  We notice that, by (\ref{220817a}), (\ref{220817a1})  
   and the property (i) of $(t_{\varepsilon},{\gamma}^{\varepsilon}_{t_{\varepsilon}})$, there exists a generic constant $C>0$ such that
 \begin{eqnarray*}
  2\sum_{i=0}^{\infty}\frac{1}{2^i}({t_{\varepsilon}}-{t}_{i})
  \leq2\sum_{i=0}^{\infty}\frac{1}{2^i}\bigg{(}\frac{\varepsilon}{2^i}\bigg{)}^{\frac{1}{2}}\leq C\varepsilon^{\frac{1}{2}};
    \end{eqnarray*}
    \begin{eqnarray*}
    |{\partial_x{\Upsilon}}({\gamma}^{\varepsilon}_{{t}_{\varepsilon}}-\hat{\gamma}_{{\hat{t}},{t}_{\varepsilon}})|\leq C|\hat{\gamma}_{{\hat{t}}}(\hat{t})-{\gamma}^{\varepsilon}_{{t}_{\varepsilon}}({t}_{\varepsilon})|^5;\ \
   |{\partial_{xx}{\Upsilon}}({\gamma}^{\varepsilon}_{{t}_{\varepsilon}}-\hat{\gamma}_{{\hat{t}},{t}_{\varepsilon}})|\leq C|\hat{\gamma}_{{\hat{t}}}(\hat{t})-{\gamma}^{\varepsilon}_{{t}_{\varepsilon}}({t}_{\varepsilon})|^4;
    \end{eqnarray*}
    \begin{eqnarray*}
  \bigg{|}\partial_x\sum_{i=0}^{\infty}\frac{1}{2^i}
                      \Upsilon({\gamma}^{\varepsilon}_{t_{\varepsilon}}-\gamma^i_{t_i,t_{\varepsilon}})
                      \bigg{|}
                      \leq18\sum_{i=0}^{\infty}\frac{1}{2^i}|\gamma^i_{t_i}({t}_{i})-{\gamma}^{\varepsilon}_{t_{\varepsilon}}(t_{\varepsilon})|^5
                     \leq18\sum_{i=0}^{\infty}\frac{1}{2^i}\bigg{(}\frac{\varepsilon}{2^i}\bigg{)}^{\frac{5}{6}}
                      \leq C{\varepsilon}^{\frac{5}{6}};
                        \end{eqnarray*}
    and
     \begin{eqnarray*}
  \bigg{|}\partial_{xx}\sum_{i=0}^{\infty}\frac{1}{2^i}
                      \Upsilon({\gamma}^{\varepsilon}_{t_{\varepsilon}}-\gamma^i_{t_i,t_{\varepsilon}})
                      \bigg{|}
                      \leq306\sum_{i=0}^{\infty}\frac{1}{2^i}|\gamma^i_{t_i}({t}_{i})-{\gamma}^{\varepsilon}_{t_{\varepsilon}}(t_{\varepsilon})|^4
                     \leq306\sum_{i=0}^{\infty}\frac{1}{2^i}\bigg{(}\frac{\varepsilon}{2^i}\bigg{)}^{\frac{2}{3}}
                      \leq C{\varepsilon}^{\frac{2}{3}}.
                        \end{eqnarray*}
 Then for any $\varrho>0$, by (\ref{sss}) and (\ref{gamma}), there exists $\varepsilon>0$ small enough such that
$$
            \hat{t}\leq {t}_{\varepsilon}< T,  \        \  
             2|{t}_{\varepsilon}-\hat{t}|+2\sum_{i=0}^{\infty}\frac{1}{2^i}({t_{\varepsilon}}-{t}_{i})\leq \frac{\varrho}{4}, $$
             and
             $$
|\partial_t{\varphi}({\gamma}^{\varepsilon}_{{t}_{\varepsilon}})-\partial_t{\varphi}(\hat{\gamma}_{\hat{t}})|\leq \frac{\varrho}{4}, \ |I|\leq \frac{\varrho}{4}, \ |II|\leq \frac{\varrho}{4},
$$
where
\begin{eqnarray*}
I&=&{\mathbf{H}}^{\varepsilon}({\gamma}^{\varepsilon}_{{t}_{\varepsilon}}, v^{\varepsilon}({\gamma}^{\varepsilon}_{{t}_{\varepsilon}}),
                           \partial_x{\varphi_{2}}({\gamma}^{\varepsilon}_{{t}_{\varepsilon}}),\partial_{xx}{\varphi_{2}}({\gamma}^{\varepsilon}_{{t}_{\varepsilon}}))
                           -{\mathbf{H}}({\gamma}^{\varepsilon}_{{t}_{\varepsilon}}, v^{\varepsilon}({\gamma}^{\varepsilon}_{{t}_{\varepsilon}}),
                           \partial_x{\varphi_{2}}({\gamma}^{\varepsilon}_{{t}_{\varepsilon}}),\partial_{xx}{\varphi_{2}}({\gamma}^{\varepsilon}_{{t}_{\varepsilon}})),\\
II&=&{\mathbf{H}}({\gamma}^{\varepsilon}_{{t}_{\varepsilon}}, v^{\varepsilon}({\gamma}^{\varepsilon}_{{t}_{\varepsilon}}),
                           \partial_x{\varphi_{2}}({\gamma}^{\varepsilon}_{{t}_{\varepsilon}}),\partial_{xx}{\varphi_{2}}({\gamma}^{\varepsilon}_{{t}_{\varepsilon}}))
                       -{\mathbf{H}}(\hat{\gamma}_{\hat{t}},{\varphi}(\hat{\gamma}_{\hat{t}}),\partial_x{\varphi}(\hat{\gamma}_{\hat{t}}),\partial_{xx}{\varphi}(\hat{\gamma}_{\hat{t}})),\\
{\varphi_{2}}({\gamma}^{\varepsilon}_{{t}_{\varepsilon}})&=&{\varphi_{1}}({\gamma}^{\varepsilon}_{{t}_{\varepsilon}})+\sum_{i=0}^{\infty}\frac{1}{2^i}\overline{\Upsilon}(\gamma^i_{t_i},{\gamma}^{\varepsilon}_{t_{\varepsilon}}),
\end{eqnarray*}
and
\begin{eqnarray*}
                                {\mathbf{H}}^{\varepsilon}(\gamma_t,r,p,\iota)&=&\sup_{u\in{
                                         {U}}}[
                        \langle p,b^{\varepsilon}(\gamma_t,u)\rangle+\frac{1}{2}\mbox{tr}[ \iota \sigma^{\varepsilon}(\gamma_t,u){\sigma^{\varepsilon}}^\top(\gamma_t,u)]\\
                        &&\ \ \ \ \ +q^{\varepsilon}(\gamma_t,r,{\sigma^{\varepsilon}}^\top(\gamma_t,u)p,u)],  \ \ (t,\gamma_t,r,p,\iota)\in [0,T]\times {\Lambda}\times \mathbb{R}\times \mathbb{R}^{d}\times \mathcal{S}(\mathbb{R}^{d}).
\end{eqnarray*}
 Since $v^{\varepsilon}$ is a viscosity subsolution of PHJB equation (\ref{hjb1}) with generators $b^{\varepsilon}, \sigma^{\varepsilon},  q^{\varepsilon}, \phi^{\varepsilon}$, we have
$$
                          \partial_t{\varphi_{2}}({\gamma}^{\varepsilon}_{{t}_{\varepsilon}})
                           +{\mathbf{H}}^{\varepsilon}({\gamma}^{\varepsilon}_{{t}_{\varepsilon}}, v^{\varepsilon}({\gamma}^{\varepsilon}_{{t}_{\varepsilon}}),
                           \partial_x{\varphi_{2}}({\gamma}^{\varepsilon}_{{t}_{\varepsilon}}),\partial_{xx}{\varphi_{2}}({\gamma}^{\varepsilon}_{{t}_{\varepsilon}}))\geq0.
$$
Thus
\begin{eqnarray*}
                       0\leq  \partial_t{\varphi}({\gamma}^{\varepsilon}_{{t}_{\varepsilon}})
                       +2({t}_{\varepsilon}-\hat{t})+2\sum_{i=0}^{\infty}\frac{1}{2^i}({t_{\varepsilon}}-{t}_{i})
                       +{\mathbf{H}}(\hat{\gamma}_{\hat{t}},{\varphi}(\hat{\gamma}_{\hat{t}}),\partial_x{\varphi}(\hat{\gamma}_{\hat{t}}),
                       \partial_{xx}{\varphi}(\hat{\gamma}_{\hat{t}}))+I+II
\leq{\mathcal{L}}{\varphi}(\hat{\gamma}_{\hat{t}})+\varrho.
\end{eqnarray*}
Letting $\varrho\downarrow 0$, we show that ${\mathcal{L}}{\varphi}(\hat{\gamma}_{\hat{t}})\geq0$.
Since ${\varphi}\in C^{1,2}_{p}(\Lambda^{\hat{t}})$ is arbitrary, we see that $v$ is a viscosity subsolution of PHJB equation (\ref{hjb1}) with generators $b,\sigma,q,\phi$, and thus completes the proof.
\ \ $\Box$

\section{Viscosity solutions to  PHJB equations: Crandall-Ishii lemma.}
\par
 In this section we extend Crandall-Ishii lemma to the path-dependent case. 
It is the cornerstone of the theory of viscosity solutions, and is the key result in
the comparison proof that will be given in the next section.
\begin{definition}\label{definition0513}  Let $(\hat{t},\hat{x})\in (0,T)\times \mathbb{R}^d$ and $f:[0,T]\times \mathbb{R}^d\rightarrow \mathbb{R}$ be an upper semicontinuous function bounded from above. We say $f\in \Phi(\hat{t},\hat{x})$ if there is a constant $r>0$ such that,
 for every $L>0$, there is  a constant $C_0\geq0$ depending only on $L$ such that, for every  function $\varphi\in C^{1,2}([0,T]\times \mathbb{R}^{d})$ such that
            $f(s,y)-\varphi(s,y)$  has a  maximum over $[0,T]\times \mathbb{R}^d$ at a point $(t,x)\in (0,T)\times \mathbb{R}^d$, 
          if   
\begin{eqnarray*}
                    |t-\hat{t}|+|x-\hat{x}|<r,\ \ \ \
                    |f(t,x)|+|\nabla_x\varphi(t,x)|
                    +|\nabla^2_x\varphi(t,x)|\leq L,
\end{eqnarray*}
then
\begin{eqnarray}\label{05281}
{\varphi}_{t}(t,x) \geq -C_0.
\end{eqnarray}
\end{definition}

\begin{definition}\label{definition0607}
 Let $\hat{t}\in [0,T)$ be fixed and  $w:\Lambda\rightarrow \mathbb{R}$ be an upper semicontinuous function bounded from above.
Define,  for $(t,x)\in [0,T]\times \mathbb{R}^{d}$,
\begin{eqnarray*}
                             &&\tilde{w}^{\hat{t}}(t,x):=\sup_{\xi_t\in \Lambda^{\hat{t}},\xi_t(t)=x}
                             w(\xi_t), \ \   t\in [\hat{t},T];\ \
                             \tilde{w}^{\hat{t}}(t,x):=\tilde{w}^{\hat{t}}(\hat{t},x)-(\hat{t}-t)^{\frac{1}{2}}, \ \   t\in[0,\hat{t}).
\end{eqnarray*}
Let $\tilde{w}^{\hat{t},*}$ be the upper
semicontinuous envelope of $\tilde{w}^{\hat{t}}$ (see    
 \cite[Definition D.10]{fab1}), i.e.,
$$
\tilde{w}^{\hat{t},*}(t,x)=\limsup_{(s,y)\in [0,T]\times \mathbb{R}^d, (s,y)\rightarrow(t,x)}\tilde{w}^{\hat{t}}(s,y).
$$
\end{definition}
In what follows, by a $modulus \ of \ continuity$, we mean a continuous function $\rho_1:[0,\infty)\rightarrow[0,\infty)$, with $\rho_1(0)=0$ and subadditivity: $\rho_1(t+s)\leq \rho_1(t)+\rho_1(s)$, for all $t,s>0$; by a $local\ modulus\ of $ $ continuity$, we mean  a continuous function $\rho_1:[0,\infty)\times[0,\infty) \rightarrow[0,\infty)$, with the properties that for each $r\geq0$, $t\rightarrow \rho_1(t,r)$ is a modulus of continuity and $\rho_1$ is non-decreasing in the  second variable.
\begin{theorem}\label{theorem0513} (Crandall-Ishii lemma)\ \  Let $w_1,w_2:\Lambda\rightarrow \mathbb{R}$ be upper semicontinuous functionals bounded from above and such that
\begin{eqnarray}\label{05131}
                     \limsup_{||\gamma_t||_0\rightarrow\infty}\frac{w_1(\gamma_t)}{||\gamma_t||_0}<0;\ \ \  \limsup_{||\gamma_t||_0\rightarrow\infty}\frac{w_2(\gamma_t)}{||\gamma_t||_0}<0.
\end{eqnarray}
Recalling  $\Lambda^t\otimes\Lambda^t:=\{(\gamma_s,\eta_s)|\gamma_s,\eta_s\in \Lambda^t\}$ for all $ t\in [0,T]$.
Let $\varphi\in C^2( \mathbb{R}^{d}\times \mathbb{R}^{d})$ be such that
$$
                w_1(\gamma_t)+w_2(\eta_t)-\varphi(\gamma_t(t),\eta_t(t))
$$
has a 
 maximum over $\Lambda^{\hat{t}}\otimes \Lambda^{\hat{t}}$ at a point $(\hat{\gamma}_{\hat{t}},\hat{\eta}_{\hat{t}})$ with $\hat{t}\in (0,T)$. 
Assume, moreover, $\tilde{w}_{1}^{\hat{t},*}\in \Phi(\hat{t},\hat{\gamma}_{\hat{t}}(\hat{t}))$ and $\tilde{w}_{2}^{\hat{t},*}\in \Phi(\hat{t},\hat{\eta}_{\hat{t}}(\hat{t}))$, and there exists a local modulus of continuity  $\rho_1$  such that, for all  $\hat{t}\leq t\leq s\leq T, \ \gamma_t\in \Lambda $,
\begin{eqnarray}\label{0608a}
w_1(\gamma_t)-w_1(\gamma_{t,s})\leq \rho_1(|s-t|,||\gamma_t||_0),\ \  w_2(\gamma_t)-w_2(\gamma_{t,s})\leq \rho_1(|s-t|,||\gamma_t||_0).
\end{eqnarray}
  Then for every $\kappa>0$, there exist
sequences  $(t_{k},\gamma^{k}_{t_k}), (s_{k},\eta^{k}_{s_k})\in [\hat{t},T]\times \Lambda^{\hat{t}}$ and
  sequences of functionals $\varphi_k\in C_p^{1,2}(\Lambda^{t_k}),\psi_k\in C_p^{1,2}(\Lambda^{s_k})$ bounded from below   
   such that 
$$
 w_{1}(\gamma_t)-\varphi_k(\gamma_t)
$$
has a strict  maximum $0$ at  $\gamma^{k}_{t_k}$ over $\Lambda^{t_k}$,
$$
w_{2}(\eta_t)-\psi_k(\eta_{t})
$$
has a strict  maximum $0$ at $\eta^{k}_{s_k}$ over $\Lambda^{s_k}$, and
\begin{eqnarray}\label{0608v}
       &&\left(t_{k}, {\gamma}^{k}_{t_{k}}(t_{k}), w_1({\gamma}^{k}_{t_{k}}),\partial_t\varphi_k({\gamma}^{k}_{t_{k}}),\partial_x\varphi_k({\gamma}^{k}_{t_{k}}),\partial_{xx}\varphi_k({\gamma}^{k}_{t_{k}})\right)\nonumber\\
       &&\underrightarrow{k\rightarrow\infty}\left({\hat{t}},\hat{{\gamma}}_{{\hat{t}}}(\hat{t}), w_1(\hat{{\gamma}}_{{\hat{t}}}),b_1, \nabla_{x_1}\varphi(\hat{{\gamma}}_{{\hat{t}}}({{\hat{t}}}),\hat{{\eta}}_{{\hat{t}}}({{\hat{t}}})), X\right),
\end{eqnarray}
\begin{eqnarray}\label{0608vw}
       &&\left(s_{k}, {\eta}^{k}_{s_{k}}(s_{k}), w_2({\eta}^{k}_{s_{k}}),\partial_t\psi_k({\eta}^{k}_{s_{k}}),\partial_x\psi_k({\eta}^{k}_{s_{k}}),\partial_{xx}\psi_k({\eta}^{k}_{s_{k}})\right)\nonumber\\
       &&\underrightarrow{k\rightarrow\infty}\left({\hat{t}},\hat{{\eta}}_{{\hat{t}}}(\hat{t}), w_2(\hat{{\eta}}_{{\hat{t}}}),b_2, \nabla_{x_2}\varphi(\hat{{\gamma}}_{{\hat{t}}}({{\hat{t}}}),\hat{{\eta}}_{{\hat{t}}}({{\hat{t}}})), Y\right),
\end{eqnarray}
 where $b_{1}+b_{2}=0$ and $X,Y\in \mathcal{S}(\mathbb{R}^{d})$ satisfy the following inequality:
\begin{eqnarray}\label{II0615}
                              -\left(\frac{1}{\kappa}+|A|\right)I\leq \left(\begin{array}{cc}
                                    X&0\\
                                    0&Y
                                    \end{array}\right)\leq A+\kappa A^2,
\end{eqnarray}
and  $A=\nabla^2_{x}\varphi(\hat{{\gamma}}_{{\hat{t}}}({{\hat{t}}}),\hat{{\eta}}_{{\hat{t}}}({{\hat{t}}}))$. Here  $\nabla^2_{x}\varphi$  denotes  the standard  second  order derivative of $\varphi$
with respect to  the  variable $x=(x_1,x_2)\in \mathbb{R}^{2d}$, and $\nabla_{x_1}\varphi$ and $\nabla_{x_2}\varphi$ denote  the standard  first order derivative of $\varphi$
with respect to  the first variable and  the second variable, respectively.  
\end{theorem}
\par
  {\bf  Proof}. \ \ 
By the following Lemma \ref{lemma4.30615}, we have that
\begin{eqnarray}\label{05211}
\tilde{w}_{1}^{\hat{t},*}(t,x)+ \tilde{w}_{2}^{\hat{t},*}(t,y)-\varphi(x,y) \ \mbox{has a  maximum over}\  [0,T]\times \mathbb{R}^{d}\times \mathbb{R}^{d} \
\mbox{ at}\ (\hat{t},\hat{{\gamma}}_{{\hat{t}}}({{\hat{t}}}),\hat{{\eta}}_{{\hat{t}}}({{\hat{t}}})).
\end{eqnarray}
 Moreover, we have
  $\tilde{w}_{1}^{\hat{t},*}(\hat{t},\hat{{\gamma}}_{{\hat{t}}}({{\hat{t}}}))={w}_{1}(\hat{{\gamma}}_{{\hat{t}}})$, $\tilde{w}_{2}^{\hat{t},*}(\hat{t},\hat{{\eta}}_{{\hat{t}}}({{\hat{t}}}))={w}_{2}(\hat{{\eta}}_{{\hat{t}}})$.
  Then, by $\tilde{w}_{1}^{\hat{t},*}\in \Phi(\hat{t},\hat{\gamma}_{\hat{t}}(\hat{t}))$, $\tilde{w}_{2}^{\hat{t},*}\in \Phi(\hat{t},\hat{\eta}_{\hat{t}}(\hat{t}))$ and Remark \ref{remarkv0528}, the Theorem 8.3 in \cite{cran2} and
  Lemma 5.4 of Chapter 4 in \cite{yong} can be used  to obtain  sequences of   functions
  $\tilde{\varphi}_k,\tilde{\psi}_k\in C^{1,2}([0,T]\times \mathbb{R}^{d})$  bounded from below such that
$\tilde{w}_{1}^{\hat{t},*}(t, x)-\tilde{\varphi}_k(t,x)$ has a strict  maximum $0$ at some point $(t_k,x^{k})\in (0,T)\times \mathbb{R}^d$ over $[0,T]\times \mathbb{R}^d$, $\tilde{w}_{2}^{\hat{t},*}(s,y)-\tilde{\psi}_k(s,y)$ has a strict  maximum $0$ at some point $(s_k,y^{k})\in (0,T)\times \mathbb{R}^d$ over $[0,T]\times \mathbb{R}^d$, and such that
\begin{eqnarray}\label{0607a}
       &&\left(t_{k}, x^k, \tilde{w}_{1}^{\hat{t},*}(t_{k}, x^k), (\tilde{\varphi}_k)_t(t_{k}, x^k),\nabla_x\tilde{\varphi}_k(t_{k}, x^k),\nabla^2_x\tilde{\varphi}_k(t_{k}, x^k)\right)\nonumber\\
       &&\underrightarrow{k\rightarrow\infty}\left({\hat{t}},\hat{{\gamma}}_{{\hat{t}}}(\hat{t}), w_1(\hat{{\gamma}}_{{\hat{t}}}),b_1, \nabla_{x_1}\varphi(\hat{{\gamma}}_{{\hat{t}}}({{\hat{t}}}),\hat{{\eta}}_{{\hat{t}}}({{\hat{t}}})), X\right),
\end{eqnarray}
\begin{eqnarray}\label{0607b}
       &&\left(s_{k}, y^k, \tilde{w}_{2}^{\hat{t},*}(s_{k}, y^k),(\tilde{\psi}_k)_t(s_{k}, y^k),\nabla_x\tilde{\psi}_k(s_{k}, y^k),\nabla^2_x\tilde{\psi}_k(s_{k}, y^k)\right)\nonumber\\
       &&\underrightarrow{k\rightarrow\infty}\left({\hat{t}},\hat{{\eta}}_{{\hat{t}}}(\hat{t}), w_2(\hat{{\eta}}_{{\hat{t}}}),b_2, \nabla_{x_2}\varphi(\hat{{\gamma}}_{{\hat{t}}}({{\hat{t}}}),\hat{{\eta}}_{{\hat{t}}}({{\hat{t}}})), Y\right),
\end{eqnarray}
%
%
%
  where $b_{1}+b_{2}=0$ and (\ref{II0615}) is satisfied.
  \\
  We claim that we can assume the  sequences $\{t_{k}\}_{k\geq1}\in [\hat{t},T)$ and $\{s_{k}\}_{k\geq1}\in [\hat{t},T)$. Indeed, if not, for example, there exists a subsequence of $\{t_{k}\}_{k\geq1}$ still denoted by itself such that $t_{k}<\hat{t}$ for all $k\geq1$.
Since $\tilde{w}_{1}^{\hat{t},*}(t, x)-\tilde{\varphi}_k(t,x)$
            has a maximum at $(t_{k},x^{k})$ on $[0,T]\times \mathbb{R}^{d}$, we obtain that
            $$
           (\tilde{\varphi}_k)_t(t_{k},x^{k})={\frac{1}{2}}(\hat{t}-t_{k})^{-\frac{1}{2}}\rightarrow\infty,\ \mbox{as}\ k\rightarrow\infty,
           $$
            which  contradicts  that $(\tilde{\varphi}_k)_t(t_{k},x^{k})\rightarrow b_{1}\in \mathbb{R}$. 
               \par
 Now we consider the functionals,
              for $(t,\gamma_t), (s,\eta_s)\in [\hat{t},T]\times{{\Lambda}}$,
\begin{eqnarray}\label{4.1111}
                 \Gamma^1_{k}(\gamma_t):=w_{1}(\gamma_t)-\tilde{\varphi}_k(t,\gamma_t(t)),\ \ \ \Gamma^2_{k}(\eta_s):= w_{2}(\eta_s)
                -\tilde{\psi}_k(s,\eta_s(s)).
\end{eqnarray}
 It is clear that $\Gamma^1_{k}, \Gamma^2_{k}$ are upper semicontinuous functionals  bounded from above on ${\Lambda}^{\hat{t}}$.  
 Define a sequence of positive numbers $\{\delta_i\}_{i\geq0}$  by 
        $\delta_i=\frac{1}{2^i}$ for all $i\geq0$. Since 
         $\overline{\Upsilon}(\cdot,\cdot)$ is a gauge-type function on $(\Lambda^{\hat{t}}, d_{\infty})$, for every $k$ and $j>0$,
           applying  Lemma \ref{theoremleft} to $\Gamma^1_{k}$ and  $\Gamma^2_{k}$, respectively, we have that,
  for every  $(\check{t}_{0},\check{\gamma}^{0}_{\check{t}_{0}}, \check{s}_{0},\check{\eta}^{0}_{\check{s}_{0}})\in ([\hat{t},T]\times  \Lambda^{\hat{t}})^2$ satisfying
\begin{eqnarray}\label{20210509a}
\Gamma^1_{k}(\check{\gamma}^{0}_{\check{t}_{0}})\geq \sup_{(t,\gamma_t)\in [\hat{t},T]\times \Lambda^{\hat{t}}}\Gamma^1_{k}(\gamma_t)-\frac{1}{j},\ \
\Gamma^2_{k}(\check{\eta}^{0}_{\check{s}_{0}})\geq \sup_{(s,\eta_s)\in [\hat{t},T]\times \Lambda^{\hat{t}}}\Gamma^2_{k}(\eta_s)-\frac{1}{j},\ \
\end{eqnarray}
  there exist $(t_{k,j},{\gamma}^{k,j}_{t_{k,j}}), (s_{k,j},{\eta}^{k,j}_{s_{k,j}})\in [\hat{t},T]\times \Lambda^{\hat{t}}$ and two sequences $\{(\check{t}_{i},\check{\gamma}^{i}_{\check{t}_{i}})\}_{i\geq1}\subset
  [\check{t}_{0},T]\times \Lambda^{\hat{t}}$, $\{(\check{s}_{i},\check{\eta}^{i}_{\check{s}_{i}})\}_{i\geq1}\subset
  [\check{s}_{0},T]\times \Lambda^{\hat{t}}$ such that
  \begin{description}
        \item{(i)} $\overline{\Upsilon}(\check{\gamma}^{0}_{\check{t}_{0}},{\gamma}^{k,j}_{t_{k,j}})\vee\overline{\Upsilon}(\check{\eta}^0_{\check{s}_0},{\eta}^{k,j}_{s_{k,j}})\leq \frac{1}{j}$,
         $\overline{\Upsilon}(\check{\gamma}^{i}_{\check{t}_{i}},{\gamma}^{k,j}_{t_{k,j}})
         \vee\overline{\Upsilon}(\check{\eta}^{i}_{\check{s}_{i}},{\eta}^{k,j}_{s_{k,j}})\leq \frac{1}{2^ij}$
          and $\check{t}_{i}\uparrow t_{k,j}$, $\check{s}_{i}\uparrow s_{k,j}$ as $i\rightarrow\infty$,
        \item{(ii)}  $\Gamma^1_k({\gamma}^{k,j}_{t_{k,j}})
            -\sum_{i=0}^{\infty}\frac{1}{2^i}\overline{\Upsilon}(\check{\gamma}^{i}_{\check{t}_{i}},{\gamma}^{k,j}_{t_{k,j}})
        \geq \Gamma^1_{k}(\check{\gamma}^{0}_{\check{t}_{0}})$ ,
        $\Gamma^2_k({\eta}^{k,j}_{s_{k,j}})
            -\sum_{i=0}^{\infty}\frac{1}{2^i}\overline{\Upsilon}(\check{\eta}^{i}_{\check{s}_{i}},{\eta}^{k,j}_{s_{k,j}})\geq \Gamma^2_{k}(\check{\eta}^{0}_{\check{s}_{0}})$,
         and
        \item{(iii)}    for all $(t,\gamma_t)\in [t_{k,j},T]\times \Lambda^{t_{k,j}}\setminus \{(t_{k,j},{\gamma}^{k,j}_{t_{k,j}})\}$,
        \begin{eqnarray*}
        \Gamma^1_{k}(\gamma_t)
        -\sum_{i=0}^{\infty}
        \frac{1}{2^i}\overline{\Upsilon}(\check{\gamma}^{i}_{\check{t}_{i}},\gamma_t)
            <\Gamma^1_{k}({\gamma}^{k,j}_{t_{k,j}})
            -\sum_{i=0}^{\infty}\frac{1}{2^i}\overline{\Upsilon}(\check{\gamma}^{i}_{\check{t}_{i}},{\gamma}^{k,j}_{t_{k,j}}),
        \end{eqnarray*}
          and for all $(s,\eta_s)\in [s_{k,j},T]\times \Lambda^{s_{k,j}}\setminus \{(s_{k,j},{\eta}^{k,j}_{s_{k,j}})\}$,
        \begin{eqnarray*}
        \Gamma^2_{k}(\eta_s)
        -\sum_{i=0}^{\infty}
        \frac{1}{2^i}\overline{\Upsilon}(\check{\eta}^{i}_{\check{s}_{i}},\eta_s)
            <\Gamma^2_{k}({\eta}^{k,j}_{s_{k,j}})
            -\sum_{i=0}^{\infty}\frac{1}{2^i}\overline{\Upsilon}(\check{\eta}^{i}_{\check{s}_{i}},{\eta}^{k,j}_{s_{k,j}}).
        \end{eqnarray*}
        \end{description}
  By the following Lemma \ref{lemma4.40615}, we have
\begin{eqnarray}\label{4.22}
        (t_{k,j}, {\gamma}^{k,j}_{t_{k,j}}(t_{k,j}))\rightarrow (t_k,x^k),\ (s_{k,j}, {\eta}^{k,j}_{s_{k,j}}(s_{k,j}))\rightarrow (s_k,y^k) \ \mbox{as}\ j\rightarrow\infty,
       \end{eqnarray}
\begin{eqnarray}\label{05231}
        \tilde{w}_{1}^{\hat{t},*}(t_{k,j}, {\gamma}^{k,j}_{t_{k,j}}(t_{k,j}))\rightarrow  \tilde{w}_{1}^{\hat{t},*}(t_k,x^k),\ \tilde{w}_{2}^{\hat{t},*}(s_{k,j},  {\eta}^{k,j}_{s_{k,j}}(s_{k,j}))\rightarrow  \tilde{w}_{2}^{\hat{t},*}(s_k,y^k) \ \mbox{as}\ j\rightarrow\infty,
       \end{eqnarray}
        and 
\begin{eqnarray}\label{05232}
        w_1({\gamma}^{k,j}_{t_{k,j}})\rightarrow  \tilde{w}_{1}^{\hat{t},*}(t_k,x^k),\ w_2({\eta}^{k,j}_{s_{k,j}})\rightarrow  \tilde{w}_{2}^{\hat{t},*}(s_k,y^k) \ \mbox{as}\ j\rightarrow\infty.
       \end{eqnarray}
        Using these and (\ref{0607a}) and (\ref{0607b}) we can therefore select a subsequence $j_k$ such that
\begin{eqnarray*}
       &&\left(t_{k,j_k}, {\gamma}^{k,j_k}_{t_{k,j_k}}(t_{k,j_k}), w_1({\gamma}^{k,j_k}_{t_{k,j_k}}),((\tilde{\varphi}_k)_t,\nabla_x\tilde{\varphi}_k,\nabla^2_x\tilde{\varphi}_k)(t_{k,j_k}, {\gamma}^{k,j_k}_{t_{k,j_k}}(t_{k,j_k}))\right)\\
       &&\underrightarrow{k\rightarrow\infty}\left({\hat{t}},\hat{{\gamma}}_{{\hat{t}}}(\hat{t}), w_1(\hat{{\gamma}}_{{\hat{t}}}),(b_1, \nabla_{x_1}\varphi(\hat{{\gamma}}_{{\hat{t}}}({{\hat{t}}}),\hat{{\eta}}_{{\hat{t}}}({{\hat{t}}})), X)\right),
\end{eqnarray*}
\begin{eqnarray*}
       &&\left(s_{k,j_k}, {\eta}^{k,j_k}_{s_{k,j_k}}(s_{k,j_k}), w_2({\eta}^{k,j_k}_{s_{k,j_k}}),((\tilde{\psi}_k)_t,\nabla_x\tilde{\psi}_k,\nabla^2_x\tilde{\psi}_k)(s_{k,j_k}, {\eta}^{k,j_k}_{s_{k,j_k}}(s_{k,j_k}))\right)\\
       &&\underrightarrow{k\rightarrow\infty}\left({\hat{t}},\hat{{\eta}}_{{\hat{t}}}(\hat{t}), w_2(\hat{{\eta}}_{{\hat{t}}}),(b_2, \nabla_{x_2}\varphi(\hat{{\gamma}}_{{\hat{t}}}({{\hat{t}}}),\hat{{\eta}}_{{\hat{t}}}({{\hat{t}}})), Y)\right).
\end{eqnarray*}
 We notice that, by   (\ref{220817a}), (\ref{220817a1})   
   and the property (i) of $(t_{k,j}, {\gamma}^{k,j}_{t_{k,j}}$, $s_{k,j}, {\eta}^{k,j}_{s_{k,j}})$, there exists a generic constant $C>0$ such that
  \begin{eqnarray*}
  2\sum_{i=0}^{\infty}\frac{1}{2^i}[(s_{k,j_k}-\check{s}_{i})+(t_{k,j_k}-\check{t}_{i})]
  \leq Cj_k^{-\frac{1}{2}};
    \end{eqnarray*}
     \begin{eqnarray*}
  \bigg{|}\partial_x\left[\sum_{i=0}^{\infty}\frac{1}{2^i}
                      \Upsilon({\gamma}^{k,j_k}_{t_{k,j_k}}-\check{\gamma}^{i}_{\check{t}_{i},t_{k,j_k}})\right]\bigg{|}
                      +\bigg{|}\partial_x\left[\sum_{i=0}^{\infty}\frac{1}{2^i}
                      \Upsilon({\eta}^{k,j_k}_{s_{k,j_k}}-\check{\eta}^{i}_{\check{s}_{i},s_{k,j_k}})\right]\bigg{|}
\leq Cj_k^{-\frac{5}{6}};
                        \end{eqnarray*}
                        and
\begin{eqnarray*}
                         \bigg{|}\partial_{xx}\left[\sum_{i=0}^{\infty}\frac{1}{2^i}
                      \Upsilon({\gamma}^{k,j_k}_{t_{k,j_k}}-\check{\gamma}^{i}_{\check{t}_{i},t_{k,j_k}})\right]\bigg{|}
                      +\bigg{|}\partial_{xx}\left[\sum_{i=0}^{\infty}\frac{1}{2^i}
                      \Upsilon({\eta}^{k,j_k}_{s_{k,j_k}}-\check{\eta}^{i}_{\check{s}_{i},s_{k,j_k}})\right]\bigg{|}
\leq Cj_k^{-\frac{2}{3}}.
  \end{eqnarray*}
Therefore the lemma holds with $\varphi_k(\gamma_t):=\tilde{\varphi}_k(t,\gamma_t(t))+\sum_{i=0}^{\infty}\frac{1}{2^i}\overline{\Upsilon}(\check{\gamma}^{i}_{\check{t}_{i}},\gamma_t)
-\tilde{\varphi}_k(t_{k,j_k},{\gamma}^{k,j_k}_{t_{k,j_k}}(t_{k,j_k}))+w_1({\gamma}^{k,j_k}_{t_{k,j_k}})-\sum_{i=0}^{\infty}\frac{1}{2^i}\overline{\Upsilon}(\check{\gamma}_{\check{t}_{i}}^{i},{\gamma}^{k,j_k}_{t_{k,j_k}})
        $, $\psi_k(\eta_s):=\tilde{\psi}_k(s,\eta_s(s))+\sum_{i=0}^{\infty}\frac{1}{2^i}\overline{\Upsilon}(\check{\eta}^{i}_{\check{s}_{i}},\eta_s)
        -\tilde{\psi}_k(s_{k,j_k},{\eta}^{k,j_k}_{s_{k,j_k}}(s_{k,j_k}))+w_2({\eta}^{k,j_k}_{s_{k,j_k}})-\sum_{i=0}^{\infty}\frac{1}{2^i}\overline{\Upsilon}(\check{\eta}_{\check{s}_{i}}^{i},{\eta}^{k,j_k}_{s_{k,j_k}})$ and $t_{k}:=t_{k,j_k}, {\gamma}^{k}_{t_k}:={\gamma}^{k,j_k}_{t_{k,j_k}}, s_{k}:=s_{k,j_k}, {\eta}^{k}_{s_k}:={\eta}^{k,j_k}_{s_{k,j_k}}$.\ \ $\Box$
        \begin{remark}\label{remarkv0528}
As mentioned in Remark 6.1 in Chapter V of \cite{fle1},  Condition (\ref{05281}) is stated with reverse inequality in Theorem 8.3 of \cite{cran2}. However, we  immediately obtain  results (\ref{0608v})-(\ref{II0615}) from Theorem
8.3 of \cite{cran2} by considering the functions $u_1(t,x) :=\tilde{w}_{1}^{\hat{t},*}(T-t,x)$ and
$u_2(t,x) := \tilde{w}_{2}^{\hat{t},*}(T-t,x)$.
\end{remark}

        \par
 To complete the  proof of Theorem \ref{theorem0513}, it remains to state and prove the following  two lemmas.
\begin{lemma}\label{lemma4.30615}\ \ Let all the conditions in Theorem  \ref{theorem0513} hold. Recall that $(\hat{t},\hat{{\gamma}}_{{\hat{t}}},\hat{{\eta}}_{{\hat{t}}})$ is given  in Theorem  \ref{theorem0513},  $\tilde{w}_{1}^{\hat{t},*}$ and $ \tilde{w}_{2}^{\hat{t},*}$ are defined in Definition \ref{definition0607} with respect to $w_1$ and $w_2$ given in Theorem  \ref{theorem0513}, respectively.
               Then $\tilde{w}_{1}^{\hat{t},*}(t,x)+ \tilde{w}_{2}^{\hat{t},*}(t,y)-\varphi(x,y)$ has a  maximum over $  [0,T]\times \mathbb{R}^{d}\times \mathbb{R}^{d} $
 at $(\hat{t},\hat{{\gamma}}_{{\hat{t}}}({{\hat{t}}}),\hat{{\eta}}_{{\hat{t}}}({{\hat{t}}}))$.
 Moreover, we have
 \begin{eqnarray}\label{202002020}
  \tilde{w}_{1}^{\hat{t},*}(\hat{t},\hat{{\gamma}}_{{\hat{t}}}({{\hat{t}}}))={w}_{1}(\hat{{\gamma}}_{{\hat{t}}}), \ \ \tilde{w}_{2}^{\hat{t},*}(\hat{t},\hat{{\eta}}_{{\hat{t}}}({{\hat{t}}}))={w}_{2}(\hat{{\eta}}_{{\hat{t}}}).
  \end{eqnarray}
\end{lemma}
\par
  {\bf  Proof}. \ \
For every $\hat{t}\leq t\leq s\leq T$ and $x\in  \mathbb{R}^{d}$, 
from the definition of $\tilde{w}^{\hat{t}}_1$ 
 it follows that
\begin{eqnarray}\label{220831}
                             \tilde{w}^{\hat{t}}_1(t,x)-\tilde{w}^{\hat{t}}_1(s,x)
                             =\sup_{\gamma_t\in \Lambda^{\hat{t}},\gamma_t(t)=x}
                             w_{1}(\gamma_t)
                             -\sup_{\eta_s\in \Lambda^{\hat{t}},\eta_s(s)=x}
                             w_{1}(\eta_s).
\end{eqnarray}
 By (\ref{05131}), there exist constants {$M_1>0$} and $\varepsilon>0$ such that 
 $$
                          \frac{w_1(\gamma_t)}{||\gamma_t||_0}\leq -\varepsilon,\ \ \mbox{if}\ ||\gamma_t||_0 \geq M_1, \ \gamma_t\in \Lambda^{\hat{t}}.
 $$
 Thus,
\begin{eqnarray}\label{07061}
                          w_1(\gamma_t)\leq -\varepsilon||\gamma_t||_0,\ \ \mbox{if}\ ||\gamma_t||_0 \geq M_1, \ \gamma_t\in \Lambda^{\hat{t}}.
\end{eqnarray}
For every $t\in[\hat{t},T]$, define $\hat{\xi}_t\in \Lambda^{\hat{t}}$ by
$$
              \hat{\xi}_t(l)=x,\ \ l\in [0,t],
$$
 then, by (\ref{0608a}),
\begin{eqnarray*}
                             w_1(\hat{\xi}_{\hat{t}})-\sup_{l\in[\hat{t},T]}\rho_1(|l-\hat{t}|,|x|)\leq w_1(\hat{\xi}_{t})\leq \sup_{\xi_t\in \Lambda^{\hat{t}},\xi_t(t)=x}w_{1}(\xi_t).
\end{eqnarray*}
Notice that $w_1(\hat{\xi}_{\hat{t}})-\sup_{l\in[\hat{t},T]}\rho_1(|l-\hat{t}|,|x|)$ depends only on $x$, there exists a constant $C^1_{x}>0$ depending only on   $x$ such that, for all $(t,\gamma_t)\in [\hat{t},T]\times\Lambda^{\hat{t}}$ satisfying $||\gamma_t||_0\geq C^1_{x}$,
\begin{eqnarray}\label{07062}
                             -\varepsilon||\gamma_t||_0< w_1(\hat{\xi}_{\hat{t}})-\sup_{l\in[\hat{t},T]}\rho_1(|l-\hat{t}|,|x|)\leq \sup_{\xi_t\in \Lambda^{\hat{t}},\xi_t(t)=x}w_{1}(\xi_t).
\end{eqnarray}
 Taking $C_{x}=M_1\vee C^1_{x}$, by (\ref{07061}) and (\ref{07062}),
 \begin{eqnarray}\label{220831a}
                   w_1(\gamma_t)<\sup_{\xi_t\in \Lambda^{\hat{t}},\xi_t(t)=x}w_{1}(\xi_t),\ \ \mbox{if} \ ||\gamma_t||_0\geq C_{x}, \ \gamma_t\in \Lambda^{\hat{t}}.
\end{eqnarray}
Combining (\ref{220831}) and (\ref{220831a}),  from (\ref{0608a}) we have
 \begin{eqnarray}\label{202105083}
                             \tilde{w}^{\hat{t}}_1(t,x)-\tilde{w}^{\hat{t}}_1(s,x)
                            &=&\sup_{\gamma_t\in \Lambda^{\hat{t}},||\gamma_t||_0\leq C_x,\gamma_t(t)=x}
                            w_{1}(\gamma_t)
                             -\sup_{\eta_s\in \Lambda^{\hat{t}},\eta_s(s)=x}
                             w_{1}(\eta_s)
                            \nonumber\\
                            &\leq& \sup_{\gamma_t\in \Lambda^{\hat{t}},||\gamma_t||_0\leq C_x,\gamma_t(t)=x}
                            {[}w_{1}(\gamma_t)
                            -w_{1}(\gamma_{t,s})
                            {]}\nonumber\\
                            &\leq&  \sup_{\gamma_t\in \Lambda^{\hat{t}},||\gamma_t||_0\leq C_x,\gamma_t(t)=x}\rho_1(|s-t|,||\gamma_t||_0)
                             \leq \rho_1(|s-t|,C_x).
\end{eqnarray}
%
%
%
Clearly, if $0\leq t\leq s\leq \hat{t}$, we have
\begin{eqnarray}\label{202105084}
                             \tilde{w}^{\hat{t}}_1(t,x)-\tilde{w}^{\hat{t}}_1(s,x)=-(\hat{t}-t)^{\frac{1}{2}}+(\hat{t}-s)^{\frac{1}{2}}\leq 0,
\end{eqnarray}
and, if $0\leq t\leq  \hat{t}\leq s\leq T$ , we have
\begin{eqnarray}\label{202105085}
                             \tilde{w}^{\hat{t}}_1(t,x)-\tilde{w}^{\hat{t}}_1(s,x)\leq\tilde{w}^{\hat{t}}_1(\hat{t},x)-\tilde{w}^{\hat{t}}_1(s,x) \leq \rho_1(|s-\hat{t}|,C_x).
\end{eqnarray}
On the other hand, 
 for every $(t,x,y)\in [0,T]\times \mathbb{R}^{d}\times \mathbb{R}^{d}$,
  by the definitions of $\tilde{w}_{1}^{\hat{t},*}(t,x)$ and $ \tilde{w}_{2}^{\hat{t},*}(t,y)$, there exist  sequences  $(l_i,x_i), (\tau_i,y_i)\in [0,T]\times \mathbb{R}^d$ such that
   $(l_{i},x_i)\rightarrow (t,x)$ and $(\tau_{i},y_i)\rightarrow (t,y)$
                 as $i\rightarrow\infty$ and
   \begin{eqnarray}\label{202105081}
       \tilde{w}_{1}^{\hat{t},*}(t,x)=\lim_{i\rightarrow\infty}\tilde{w}^{\hat{t}}_{1}(l_i,x_i), \ \ \ \tilde{w}_{2}^{\hat{t},*}(t,y)=\lim_{i\rightarrow\infty}\tilde{w}^{\hat{t}}_2(\tau_i,y_i).
\end{eqnarray}
Without loss of generality, we may assume  $l_i\leq \tau_i$ for all $i>0$.
 By (\ref{202105083})-(\ref{202105085}), we have
 \begin{eqnarray}\label{202105082}
       \tilde{w}_{1}^{\hat{t},*}(t,x)=\lim_{i\rightarrow\infty}\tilde{w}^{\hat{t}}_{1}(l_i,x_i)\leq \liminf_{i\rightarrow\infty}[\tilde{w}^{\hat{t}}_{1}(\tau_i,x_i)+\rho_1(|\tau_i-l_i|,C_{x_i})].
\end{eqnarray}
Here $C_{x_i}>0$ is the constant  that makes the following formula true:
$$
                                  w_1(\gamma_t)<\sup_{\xi_t\in \Lambda^{\hat{t}},\xi_t(t)=x_i}w_{1}(\xi_t),\ \ \mbox{if} \ ||\gamma_t||_0\geq C_{x_i}, \ \gamma_t\in \Lambda^{\hat{t}}.
                                 %
$$
We claim that we can assume that there exists a constant $M_2>0$ such that $C_{x_i}\leq M_2$ for all $i\geq 1$. Indeed, if not, for every $n$, there exists  $i_n$ such that
\begin{eqnarray}
\tilde{w}^{\hat{t}}_{1}(l_{i_n},x_{i_n})
=\begin{cases}\sup_{\gamma_{l_{i_n}}\in \Lambda^{\hat{t}},||\gamma_{l_{i_n}}||_0> n,\gamma_{l_{i_n}}(l_{i_n})=x_{i_n}}
                             {[}w_{1}(\gamma_{l_{i_n}}){]}, \ \ \ \ \ \  \ l_{i_n}\geq \hat{t};\\
                             \sup_{\gamma_{\hat{t}}\in \Lambda^{\hat{t}},||\gamma_{\hat{t}}||_0> n,\gamma_{\hat{t}}({\hat{t}})=x_{i_n}}
                             {[}w_{1}(\gamma_{\hat{t}}){]}-({\hat{t}}-l_{i_n})^{\frac{1}{2}}, \ \ l_{i_n}< \hat{t}.
\end{cases}
\end{eqnarray}
Letting $n\rightarrow\infty$, by (\ref{05131}), we get that
$$
\tilde{w}^{\hat{t}}_{1}(l_{i_n},x_{i_n})\rightarrow-\infty \ \mbox{as}\ n\rightarrow\infty,
$$
 which contradicts the convergence that $ \tilde{w}_{1}^{\hat{t},*}(t,x)=\lim_{i\rightarrow\infty}\tilde{w}^{\hat{t}}_{1}(l_i,x_i)$. Then, by  (\ref{202105082}),
 \begin{eqnarray}\label{20210704a}
                           \tilde{w}_{1}^{\hat{t},*}(t,x)
                           \leq \liminf_{i\rightarrow\infty}[\tilde{w}^{\hat{t}}_{1}(\tau_i,x_i)+\rho_1(|\tau_i-l_i|,M_2)] = \liminf_{i\rightarrow\infty}\tilde{w}^{\hat{t}}_{1}(\tau_i,x_i).
\end{eqnarray}
Therefore, by (\ref{202105081}), (\ref{20210704a}) and the definitions of $\tilde{w}^{\hat{t}}_{1}$ and $\tilde{w}^{\hat{t}}_{2}$,
\begin{eqnarray}\label{20210508b}
                            &&\tilde{w}_{1}^{\hat{t},*}(t,x)+\tilde{w}_{2}^{\hat{t},*}(t,y)-\varphi(x,y)\nonumber\\
                          &\leq&\liminf_{i\rightarrow\infty}[\tilde{w}^{\hat{t}}_{1}(\tau_i,x_i)+\tilde{w}^{\hat{t}}_{2}(\tau_i,y_i)-\varphi(x_i,y_i)]\nonumber\\
                           &\leq&\sup_{(l,x_0,y_0)\in [0,T]\times \mathbb{R}^{d}\times \mathbb{R}^{d}}[\tilde{w}^{\hat{t}}_{1}(l,x_0)+\tilde{w}^{\hat{t}}_{2}(l,y_0)-\varphi(x_0,y_0)]\nonumber\\
                            &=&\sup_{(l,x_0,y_0)\in [\hat{t},T]\times \mathbb{R}^{d}\times \mathbb{R}^{d}}[\tilde{w}^{\hat{t}}_{1}(l,x_0)+\tilde{w}^{\hat{t}}_{2}(l,y_0)-\varphi(x_0,y_0)].
\end{eqnarray}
We also have, for $(l,x_0,y_0)\in [\hat{t},T]\times \mathbb{R}^d\times \mathbb{R}^d$,
\begin{eqnarray}\label{0608aa}
                         &&\tilde{w}^{\hat{t}}_{1}(l,x_0)+  \tilde{w}^{\hat{t}}_{2}(l,y_0)-\varphi(x_0,y_0)\nonumber   \\
                          &=&\sup_{
                          \gamma_l,\eta_l\in\Lambda^{\hat{t}},\gamma_l(l)=x_0,
                          \eta_l(l)=y_0}
                          \left[w_{1}(\gamma_l)
                          +w_{2}(\eta_l)
                          -\varphi(\gamma_l(l),\eta_l(l))\right]\nonumber \\
                          &\leq&w_{1}(\hat{{\gamma}}_{{\hat{t}}})+w_{2}(\hat{{\eta}}_{{\hat{t}}})-\varphi(\hat{{\gamma}}_{{\hat{t}}}({\hat{t}}),\hat{{\eta}}_{{\hat{t}}}({\hat{t}})),
\end{eqnarray}
                       where the  inequality becomes equality if  $l={\hat{t}}$ 
                       and $x_0=\hat{{\gamma}}_{{\hat{t}}}({\hat{t}}),y_0=\hat{{\eta}}_{{\hat{t}}}({\hat{t}})$.
                       Combining  (\ref{20210508b}) and (\ref{0608aa}), we obtain that
\begin{eqnarray}\label{0608a1}
                            &&\tilde{w}_{1}^{\hat{t},*}(t,x)+\tilde{w}_{2}^{\hat{t},*}(t,y)-\varphi(x,y)\leq w_{1}(\hat{{\gamma}}_{{\hat{t}}})+w_{2}(\hat{{\eta}}_{{\hat{t}}})-\varphi(\hat{{\gamma}}_{{\hat{t}}}({\hat{t}}),\hat{{\eta}}_{{\hat{t}}}({\hat{t}})).
\end{eqnarray}
By the definitions of $\tilde{w}_1^{\hat{t},*}$ and $\tilde{w}_2^{\hat{t},*}$, we have $\tilde{w}_1^{\hat{t},*}(t,x)\geq \tilde{w}^{\hat{t}}_1(t,x), \tilde{w}_{2}^{\hat{t},*}(t,y) \geq \tilde{w}_{2}(t,y)$.  Then by also (\ref{0608aa}) and (\ref{0608a1}), for every $(t,x,y)\in [0,T]\times \mathbb{R}^{d}\times \mathbb{R}^{d}$,
\begin{eqnarray}\label{0608abc}
                            &&\tilde{w}_{1}^{\hat{t},*}(t,x)+\tilde{w}_{2}^{\hat{t},*}(t,y)-\varphi(x,y)\leq w_{1}(\hat{{\gamma}}_{{\hat{t}}})+w_{2}(\hat{{\eta}}_{{\hat{t}}})-\varphi(\hat{{\gamma}}_{{\hat{t}}}({\hat{t}}),\hat{{\eta}}_{{\hat{t}}}({\hat{t}}))\nonumber\\
                            &=&\tilde{w}^{\hat{t}}_{1}(\hat{t},\hat{{\gamma}}_{{\hat{t}}}({\hat{t}}))+  \tilde{w}^{\hat{t}}_{2}(\hat{t},\hat{{\eta}}_{{\hat{t}}}({\hat{t}}))-\varphi(\hat{{\gamma}}_{{\hat{t}}}({\hat{t}}),\hat{{\eta}}_{{\hat{t}}}({\hat{t}}))\nonumber\\
                            &\leq&\tilde{w}_{1}^{\hat{t},*}(\hat{t},\hat{{\gamma}}_{{\hat{t}}}({\hat{t}}))+  \tilde{w}_{2}^{\hat{t},*}(\hat{t},\hat{{\eta}}_{{\hat{t}}}({\hat{t}}))-\varphi(\hat{{\gamma}}_{{\hat{t}}}({\hat{t}}),\hat{{\eta}}_{{\hat{t}}}({\hat{t}})).
\end{eqnarray}
Thus 
we obtain that (\ref{202002020}) holds true, and  $ \tilde{w}_{1}^{\hat{t},*}(t,x)+\tilde{w}_{2}^{\hat{t},*}(t,y)-\varphi(x, y)$ has a maximum at $({\hat{t}},\hat{{\gamma}}_{\hat{t}}({\hat{t}}),
\hat{{\eta}}_{\hat{t}}({\hat{t}}))$
 on $[0,T]\times \mathbb{R}^{d}\times \mathbb{R}^{d}$.
The proof is now complete. \ \ $\Box$
\begin{lemma}\label{lemma4.40615}\ \
Let all the conditions in Theorem  \ref{theorem0513} hold. Recall that $\Gamma^1_{k}$ and $\Gamma^2_{k}$ are defined in  (\ref{4.1111}).
Then the maximum points ${\gamma}^{k,j}_{t_{k,j}}$ 
of $\Gamma^1_{k}(\gamma_t)
        -\sum_{i=0}^{\infty}
        \frac{1}{2^i}\overline{\Upsilon}(\check{\gamma}^{i}_{\check{t}_{i}},\gamma_t)$ and the maximum points
         $ {\eta}^{k,j}_{s_{k,j}}$ 
of $\Gamma^2_{k}(\eta_s)
        -\sum_{i=0}^{\infty}
        \overline{\Upsilon}(\check{\eta}^{i}_{\check{s}_{i}},\eta_s)$
        satisfy  conditions (\ref{4.22}), (\ref{05231}) and (\ref{05232}).
\end{lemma}
\par
  {\bf  Proof}. \ \
   Recall that $\tilde{w}_{1}^{\hat{t},*}\geq \tilde{w}^{\hat{t}}_{1}, \tilde{w}_{2}^{\hat{t},*} \geq \tilde{w}^{\hat{t}}_{2}$,
by  
the definitions of $\tilde{w}^{\hat{t}}_1$ and $\tilde{w}^{\hat{t}}_2$, we get that
 \begin{eqnarray*}
                           && \tilde{w}_{1}^{\hat{t},*}(t_{k,j},{\gamma}^{k,j}_{t_{k,j}}(t_{k,j}))
                           -\tilde{\varphi}_k(t_{k,j},{\gamma}^{k,j}_{t_{k,j}}(t_{k,j}))
                           \geq w_1({\gamma}^{k,j}_{t_{k,j}})-\tilde{\varphi}_k(t_{k,j},{\gamma}^{k,j}_{t_{k,j}}(t_{k,j}))
                           =\Gamma^1_{k}({\gamma}^{k,j}_{t_{k,j}}),\\
                            &&\tilde{w}_{2}^{\hat{t},*}(s_{k,j},{\eta}^{k,j}_{s_{k,j}}(s_{k,j}))
                           -\tilde{\psi}_k(s_{k,j},{\eta}^{k,j}_{s_{k,j}}(s_{k,j}))
                           \geq w_2({\eta}^{k,j}_{s_{k,j}})-\tilde{\psi}_k(s_{k,j},{\eta}^{k,j}_{s_{k,j}}(s_{k,j}))
                           =\Gamma^2_{k}({\eta}^{k,j}_{s_{k,j}}).
  \end{eqnarray*}
We notice that, from (\ref{20210509a}) and the property (ii) of $(t_{k,j},{\gamma}^{k,j}_{t_{k,j}}, s_{k,j},{\eta}^{k,j}_{s_{k,j}})$,
\begin{eqnarray*}
                 \Gamma^1_{k}({\gamma}^{k,j}_{t_{k,j}})\geq\Gamma^1_{k}(\check{\gamma}^{0}_{\check{t}_{0}})
                 \geq \sup_{(t,\gamma_t)\in [\hat{t},T]\times \Lambda^{\hat{t}}}\Gamma^1_{k}(\gamma_t)-\frac{1}{j}, \ \ \
                  \Gamma^2_{k}({\eta}^{k,j}_{s_{k,j}})\geq\Gamma^2_{k}(\check{\eta}^{0}_{\check{s}_{0}})
                 \geq \sup_{(s,\eta_s)\in [\hat{t},T]\times \Lambda^{\hat{t}}}\Gamma^2_{k}(\eta_s)-\frac{1}{j},
\end{eqnarray*}
and by  
the  definitions of $\tilde{w}_{1}^{\hat{t},*}$ and $\tilde{w}_{2}^{\hat{t},*}$,
$$
\sup_{(t,\gamma_t)\in [\hat{t},T]\times \Lambda^{\hat{t}}}\Gamma^1_{k}(\gamma_t)\geq\tilde{w}_{1}^{\hat{t},*}(t_k,{x}^k)-\tilde{\varphi}_k(t_k,{x}^k), \ \
\sup_{(s,\eta_s)\in [\hat{t},T]\times \Lambda^{\hat{t}}}\Gamma^2_{k}(\eta_s)\geq\tilde{w}_{2}^{\hat{t},*}(s_k,{y}^k)-\tilde{\psi}_k(s_k,{y}^k).
$$
Therefore,
\begin{eqnarray}\label{0525b}
\tilde{w}_{1}^{\hat{t},*}(t_{k,j},{\gamma}^{k,j}_{t_{k,j}}(t_{k,j}))-\tilde{\varphi}_k(t_{k,j},{\gamma}^{k,j}_{t_{k,j}}(t_{k,j}))
                  \geq\Gamma^1_{k}({\gamma}^{k,j}_{t_{k,j}})\geq
                  \tilde{w}_{1}^{\hat{t},*}(t_k,{x}^k)-\tilde{\varphi}_k(t_k,{x}^k)-\frac{1}{j},
\end{eqnarray}
\begin{eqnarray}\label{0525b0920}
\tilde{w}_{2}^{\hat{t},*}(s_{k,j},{\eta}^{k,j}_{s_{k,j}}(s_{k,j}))-\tilde{\psi}_k(s_{k,j},{\eta}^{k,j}_{s_{k,j}}(s_{k,j}))
                  \geq\Gamma^2_{k}({\eta}^{k,j}_{s_{k,j}})\geq
                  \tilde{w}_{2}^{\hat{t},*}(s_k,{y}^k)-\tilde{\psi}_k(s_k,{y}^k)-\frac{1}{j}.
\end{eqnarray}
By (\ref{05131}) and  that $\tilde{\varphi}_k$ and $\tilde{\psi}_k$ are bounded from below, 
 there exists a constant  $M_3>0$  that is sufficiently  large   that
 $$
           \Gamma^1_k(\gamma_t)<\sup_{(l,\xi_l)\in [\hat{t},T]\times \Lambda^{\hat{t}}}\Gamma^1_{k}(\xi_l)-1, \ \ t\in [\hat{t},T],\ ||\gamma_t||_0\geq M_3,
           $$ and
 $$
 \Gamma^2_{k}(\eta_s)<\sup_{(r,\varsigma_r)\in [\hat{t},T]\times \Lambda^{\hat{t}}}\Gamma^2_{k}(\varsigma_r)-1,\ s\in [\hat{t},T], \ ||\eta_s||_0\geq M_3.
 $$
  Thus, we have $||{\gamma}^{k,j}_{t_{k,j}}||_0\vee
           ||{\eta}^{k,j}_{s_{k,j}}||_0<M_3$. In particular, $|{\gamma}^{k,j}_{t_{k,j}}(t_{k,j})|\vee|{\eta}^{k,j}_{s_{k,j}}(s_{k,j})|<M_3$. We note that $M_3$ is independent of $j$.
 Then letting $j\rightarrow\infty$ in (\ref{0525b}) and (\ref{0525b0920}), 
   we obtain (\ref{4.22}).  Indeed, if not, we may assume that  there exist $(\grave{t},\grave{x},\grave{s},\grave{y})\in [0, T]\times \mathbb{R}^{d}\times [0, T]\times \mathbb{R}^{d}$ and  a subsequence of $(t_{k,j},{\gamma}^{k,j}_{t_{k,j}}(t_{k,j}),s_{k,j},{\eta}^{k,j}_{s_{k,j}}(s_{k,j}))$ still denoted by itself  such that
       $$
       (t_{k,j},{\gamma}^{k,j}_{t_{k,j}}(t_{k,j}),s_{k,j},{\eta}^{k,j}_{s_{k,j}}(s_{k,j}))\rightarrow  (\grave{t},\grave{x},\grave{s},\grave{y})\neq (t_k,x^k,s_k,y^k).
       $$
Letting $j\rightarrow\infty$ in (\ref{0525b}) and (\ref{0525b0920}), by the upper semicontinuity of $\tilde{w}_{1}^{\hat{t},*}+\tilde{w}_{2}^{\hat{t},*}-\tilde{\varphi}_k-\tilde{\psi}_k$, we have
$$
           \tilde{w}_{1}^{\hat{t},*}(\grave{t},\grave{x})+\tilde{w}_{2}^{\hat{t},*}(\grave{s},\grave{y})-\tilde{\varphi}_k(\grave{t},\grave{x})-\tilde{\psi}_k(\grave{s},\grave{y})
                 \geq\tilde{w}_{1}^{\hat{t},*}(t_k,{x}^k)+\tilde{w}_{2}^{\hat{t},*}(s_k,{y}^k) -\tilde{\varphi}_k(t_k,{x}^k)-\tilde{\psi}_k(s_k,{y}^k),
$$
       which  contradicts that 
       $(t_k,x^{k},s_k,y^{k})$  is the  strict  maximum point of $\tilde{w}_{1}^{\hat{t},*}(t,x)+\tilde{w}_{2}^{\hat{t},*}(s,y)-\tilde{\varphi}_k(t,x)-\tilde{\psi}_k(s,y)$ on $[0,T]\times\mathbb{R}^d\times [0,T]\times\mathbb{R}^d$.\par
       By (\ref{4.22}), the  upper semicontinuity of $\tilde{w}_{1}^{\hat{t},*}$ and $\tilde{w}_{2}^{\hat{t},*}$ and  the continuity of $\tilde{\varphi}_k$ and $\tilde{\psi}_k$, letting $j\rightarrow\infty$ in (\ref{0525b})  and (\ref{0525b0920}), we obtain
       \begin{eqnarray*}
                  &&\tilde{w}_{1}^{\hat{t},*}(t_k,x_k)
                  \geq \limsup_{j\rightarrow\infty}\tilde{w}_{1}^{\hat{t},*}(t_{k,j},{\gamma}_{t_{k,j}}^{k,j}(t_{k,j}))
                   \geq \liminf_{j\rightarrow\infty}\tilde{w}_{1}^{\hat{t},*}(t_{k,j},{\gamma}_{t_{k,j}}^{k,j}(t_{k,j}))
                  \geq\tilde{w}_{1}^{\hat{t},*}(t_k,x_k),\\
                  &&\tilde{w}_{2}^{\hat{t},*}(s_k,y_k)
                  \geq \limsup_{j\rightarrow\infty}\tilde{w}_{2}^{\hat{t},*}(s_{k,j},{\eta}_{t_{k,j}}^{k,j}(t_{k,j}))
                   \geq \liminf_{j\rightarrow\infty}\tilde{w}_{2}^{\hat{t},*}(s_{k,j},{\eta}_{t_{k,j}}^{k,j}(t_{k,j}))
                  \geq\tilde{w}_{2}^{\hat{t},*}(s_k,y_k).
\end{eqnarray*}
  Thus, we get (\ref{05231}) holds true. 
  Letting $j\rightarrow\infty$ in (\ref{0525b}) and (\ref{0525b0920}), by (\ref{05231}) and the definitions of $\Gamma^1_{k}$ and $\Gamma^2_{k}$,
  $$
                             \tilde{w}_{1}^{\hat{t},*}(t_k,x_k)=\lim_{j\rightarrow\infty}w_1({\gamma}_{t_{k,j}}^{k,j}), \ \ \
                             \tilde{w}_{2}^{\hat{t},*}(s_k,y_k)=\lim_{j\rightarrow\infty}w_2({\eta}_{t_{k,j}}^{k,j}).
  $$
  Thus, we obtain (\ref{05232}).
  The proof is now complete. \ \ $\Box$

\section{Viscosity solution to  PHJB equation: Uniqueness theorem.}

\par
%
             This section is devoted to a  proof of uniqueness of  viscosity
                   solutions to (\ref{hjb1}). This result, together with
                  the results from  Section 4, will be used to characterize
                   the value functional defined by (\ref{value1}).
                   \par
By    \cite[Proposition 11.2.13]{zhang},  without loss of generality we assume that there exists a constant $K>0$ such that,
for all $(t,\gamma_t, p,\iota)\in [0,T]\times{\Lambda}\times \mathbb{R}^{d}\times \mathcal{S}(\mathbb{R}^{d})$ and $r_1,r_2\in \mathbb{R}$ such that $r_1<r_2$,
\begin{eqnarray}\label{5.1}
{\mathbf{H}}(\gamma_t,r_1,p,\iota)-{\mathbf{H}}(\gamma_t,r_2,p,\iota)\geq K(r_2-r_1).
\end{eqnarray}
We  now state the main result of this section.
\begin{theorem}\label{theoremhjbm}
Suppose Hypothesis \ref{hypstate}   holds.
                         Let $W_1\in C^0({\Lambda})$ $(\mbox{resp}., W_2\in C^0({\Lambda}))$ be  a viscosity subsolution (resp., supersolution) to equation (\ref{hjb1}) and  let  there exist  constant $L>0$ and a local modulus of continuity  $\rho_2$
                        such that, for any  
                        $(t,\gamma_t),(s,\eta_s)\in[0,T]\times{\Lambda}$,
\begin{eqnarray}\label{w}
                                   |W_1(\gamma_t)|\vee |W_2(\gamma_t)|\leq L (1+||\gamma_t||_0);
                                   \end{eqnarray}
\begin{eqnarray}\label{w1}
                               |W_1(\gamma_t)-W_1(\eta_s)|\vee|W_2(\gamma_t)-W_2(\eta_s)|\leq
                       \rho_2(|s-t|,||\gamma_t||_0\vee||\eta_s||_0)+L||\gamma_{t}-\eta_s||_0.
\end{eqnarray}
                   Then  $W_1\leq W_2$. 
%
%
%
%
\end{theorem}
\par
                      Theorems    \ref{theoremvexist} and \ref{theoremhjbm} lead to the result (given below) that the viscosity solution to   PHJB equation given in (\ref{hjb1})
                      corresponds to the value functional  $V$ of our optimal control problem given in (\ref{state1}), (\ref{fbsde1}) and (\ref{value1}).
\begin{theorem}\label{theorem52}\ \
                 Let Hypothesis \ref{hypstate}  hold. Then the value
                          functional $V$ defined by (\ref{value1}) is the unique viscosity
                          solution to (\ref{hjb1}) in the class of functionals satisfying (\ref{w}) and (\ref{w1}).
\end{theorem}
\par
  {\bf  Proof}. \ \   Theorem \ref{theoremvexist} shows that $V$ is a viscosity solution to equation (\ref{hjb1}).  Thus, our conclusion follows from 
   Theorems  \ref{theorem3.9} and
    \ref{theoremhjbm}.  \ \ $\Box$

\par
  Next, we prove Theorem \ref{theoremhjbm}.   
 We  note that for $\varrho>0$, the functional
                    defined by $\tilde{W}:=W_1-\frac{\varrho}{t+1}$ is a viscosity subsolution
                   for
 \begin{eqnarray}\label{0619c}
\begin{cases}
{\partial_t} \tilde{W}(\gamma_t)+{\mathbf{H}}(\gamma_t, \tilde{W}(\gamma_t), \partial_x \tilde{W}(\gamma_t),\partial_{xx} \tilde{W}(\gamma_t))
          = \frac{\varrho}{(t+1)^2}, \ \  (t,\gamma_t)\in [0,T)\times {\Lambda}, \\
\tilde{W}(\gamma_T)=\phi(\gamma_T), \ \  \gamma_T\in \Lambda_T.
\end{cases}
\end{eqnarray}
We mention  that  $\tilde{W}$ is also a viscosity subsolution of (\ref{0619c}) if the second argument of $\mathbf{H}$ is ${W_1}(\gamma_t)$ instead of $\tilde{W}(\gamma_t)$.
                As $W_1\leq W_2$ follows from $\tilde{W}\leq W_2$ in
                the limit $\varrho\downarrow0$, it suffices to prove
                $W_1\leq W_2$ under the additional assumption given below:
$$
{\partial_t} {W_1}(\gamma_t)+{\mathbf{H}}(\gamma_t, {W_1}(\gamma_t), \partial_x {W_1}(\gamma_t),\partial_{xx}{W_1}(\gamma_t))
          \geq c,\ \ c:=\frac{\varrho}{(T+1)^2}, \ \  (t,\gamma_t)\in [0,T)\times {\Lambda}.
$$
\par
   {\bf  Proof of Theorem \ref{theoremhjbm}}. \ \   The proof of this theorem  is rather long. Thus, we split it into several
        steps.
\par
            $Step\  1.$ Definitions of auxiliary functionals.
            \par
 We only need to prove that $W_1(\gamma_t)\leq W_2(\gamma_t)$ for all $(t,\gamma_t)\in
[T-\bar{a},T)\times
       {\Lambda}$.
        Here,
        $$\bar{a}=\frac{1}{2(342L+36)L}\wedge{T}.$$
         Then, we can  repeat the same procedure for the case
        $[T-i\bar{a},T-(i-1)\bar{a})$.  Thus, we assume the converse result that $(\tilde{t},\tilde{\gamma}_{\tilde{t}})\in (T-\bar{a},T)\times
      {\Lambda}$ exists  such that
        $\tilde{m}:=W_1(\tilde{\gamma}_{\tilde{t}})-W_2(\tilde{\gamma}_{\tilde{t}})>0$. 
\par
         Consider  a small number $\varepsilon >0$  such that
 $$
 W_1(\tilde{\gamma}_{\tilde{t}})-W_2(\tilde{\gamma}_{\tilde{t}})-2\varepsilon \frac{\nu T-\tilde{t}}{\nu
 T}\Upsilon(\tilde{\gamma}_{\tilde{t}})
 >\frac{\tilde{m}}{2},
 $$
      and
\begin{eqnarray}\label{5.3}
                          \frac{\varepsilon}{\nu T}\leq\frac{c}{4},
\end{eqnarray}
             where
$$
            \nu=1+\frac{1}{2T(342L+36)L}.
$$
 Next, recalling that $\Lambda^t\otimes\Lambda^t:=\{(\gamma_s,\eta_s)|\gamma_s,\eta_s\in \Lambda^t\}$ for all $ t\in [0,T]$, we define for any  $\beta\in (0,\infty)$ and   $(\gamma_t,\eta_t)\in {\Lambda}^{T-\bar{a}}\otimes{\Lambda}^{T-\bar{a}}$,
\begin{eqnarray*}
                 \Psi(\gamma_t,\eta_t)=W_1(\gamma_t)-W_2(\eta_t)-{\beta}\Upsilon(\gamma_{t},\eta_{t})-\beta^{\frac{1}{3}}|\gamma_{t}(t)-\eta_{t}(t)|^2
                 -\varepsilon\frac{\nu T-t}{\nu
                 T}(\Upsilon(\gamma_t)+\Upsilon(\eta_t)).
\end{eqnarray*}
 By (\ref{s0}) and (\ref{w}), it is clear that $\Psi$ is bounded from above on ${\Lambda}^{T-\bar{a}}\otimes{\Lambda}^{T-\bar{a}}$. Moreover,  by Lemma \ref{theoremS}, $\Psi$ is a continuous  functional.
        Define a sequence of positive numbers $\{\delta_i\}_{i\geq0}$  by 
        $\delta_i=\frac{1}{2^i}$ for all $i\geq0$.
           Since 
            $\overline{\Upsilon}^{3,3,2}(\cdot,\cdot)$ is a gauge-type function on $(\Lambda^{\tilde{t}}\otimes \Lambda^{\tilde{t}},d_{1,\infty})$, from Lemma \ref{theoremleft1} it follows that,
  for every  $(\gamma^0_{t_0},\eta^0_{t_0})\in  \Lambda^{\tilde{t}}\otimes \Lambda^{\tilde{t}}$ satisfying
$$
\Psi(\gamma^0_{t_0},\eta^0_{t_0})\geq \sup_{(s,(\gamma_s,\eta_s))\in [\tilde{t},T]\times (\Lambda^{\tilde{t}}\otimes \Lambda^{\tilde{t}})}\Psi(\gamma_s,\eta_s)-\frac{1}{\beta},\
\    \mbox{and} \ \ \Psi(\gamma^0_{t_0},\eta^0_{t_0})\geq \Psi(\tilde{\gamma}_{\tilde{t}},\tilde{\gamma}_{\tilde{t}}) >\frac{\tilde{m}}{2},
 $$
  there exist $(\hat{t},(\hat{\gamma}_{\hat{t}},\hat{\eta}_{\hat{t}}))\in [\tilde{t},T]\times (\Lambda^{\tilde{t}}\otimes \Lambda^{\tilde{t}})$ and a sequence $\{(t_i,(\gamma^i_{t_i},\eta^i_{t_i}))\}_{i\geq1}\subset
  [t_0,T]\times (\Lambda^{\tilde{t}}\otimes \Lambda^{\tilde{t}})$ such that
  \begin{description}
        \item{(i)} $\Upsilon(\gamma^0_{t_0},\hat{\gamma}_{\hat{t}})+\Upsilon(\eta^0_{t_0},\hat{\eta}_{\hat{t}})+|\hat{t}-t_0|^2\leq \frac{1}{\beta}$,
         $\Upsilon(\gamma^i_{t_i},\hat{\gamma}_{\hat{t}})+\Upsilon(\eta^i_{t_i},\hat{\eta}_{\hat{t}})+|\hat{t}-t_i|^2
         \leq \frac{1}{\beta2^i}$ and $t_i\uparrow \hat{t}$ as $i\rightarrow\infty$,
        \item{(ii)}  $\Psi_1(\hat{\gamma}_{\hat{t}},\hat{\eta}_{\hat{t}})
\geq \Psi(\gamma^0_{t_0},\eta^0_{t_0})$, and
        \item{(iii)}    for all $(s,(\gamma_s,\eta_s))\in [\hat{t},T]\times (\Lambda^{\hat{t}}\otimes \Lambda^{\hat{t}})\setminus \{(\hat{t},(\hat{\gamma}_{\hat{t}},\hat{\eta}_{\hat{t}}))\}$,
        \begin{eqnarray}\label{iii4}
        \Psi_1(\gamma_s,\eta_s)
           <\Psi_1(\hat{\gamma}_{\hat{t}},\hat{\eta}_{\hat{t}}),
        \end{eqnarray}
        \end{description}
        where we define
$$
     \Psi_1(\gamma_t,\eta_t):=  \Psi(\gamma_t,\eta_t)
        -\sum_{i=0}^{\infty}
        \frac{1}{2^i}[{\Upsilon}(\gamma^i_{t_i},\gamma_t)+{\Upsilon}(\eta^i_{t_i},\eta_t)+|t-t_i|^2], \ \  \  (\gamma_t,\eta_t)\in  \Lambda^{\tilde{t}}\otimes \Lambda^{\tilde{t}}.
$$
             We should note that the point
             $({\hat{t}},(\hat{{\gamma}}_{{\hat{t}}},\hat{{\eta}}_{{\hat{t}}}))$ depends on $\beta$ and
              $\varepsilon$.
\par
$Step\ 2.$
There exists ${{M}_0}>0$ independent of $\beta$
    such that
                   \begin{eqnarray}\label{5.10jiajiaaaa}||\hat{\gamma}_{\hat{t}}||_0\vee||\hat{\eta}_{\hat{t}}||_0<M_0,
                   \end{eqnarray} and
   the following result  holds true:
 \begin{eqnarray}\label{5.10}
                      \beta ||\hat{{\gamma}}_{{\hat{t}}}-\hat{{\eta}}_{{\hat{t}}}||_{0}^6
                         +\beta|\hat{{\gamma}}_{{\hat{t}}}(\hat{t})-\hat{{\eta}}_{{\hat{t}}}(\hat{t})|^4
                         \rightarrow0 \ \mbox{as} \ \beta\rightarrow\infty.
 \end{eqnarray}
  Let us show the above. First,   noting $\nu$ is independent of  $\beta$, by the definition of  ${\Psi}$,
 there exists a constant  ${M}_0>0$ independent of $\beta$ that is sufficiently  large   that
           $
           \Psi(\gamma_t, \eta_t)<0
           $ for all $t\in [T-\bar{a},T]$ and $||\gamma_t||_0\vee||\eta_t||_0\geq M_0$. Thus, we have $||\hat{\gamma}_{\hat{t}}||_0\vee||\hat{\eta}_{\hat{t}}||_0\vee
           ||{\gamma}^{0}_{t_{0}}||_0\vee||{\eta}^{0}_{t_{0}}||_0<M_0$.  We  note that $M_0$ depends on
              $\varepsilon$.
           \par
   Second, by (\ref{iii4}), we have
 \begin{eqnarray}\label{5.56789}
                        2\Psi_1(\hat{\gamma}_{\hat{t}},\hat{\eta}_{\hat{t}})
            \geq  \Psi_1(\hat{\gamma}_{\hat{t}},\hat{\gamma}_{\hat{t}})
            +\Psi_1(\hat{\eta}_{\hat{t}},\hat{\eta}_{\hat{t}}).
 \end{eqnarray}
 This implies that
 \begin{eqnarray}\label{5.6}
                         &&2{\beta}\Upsilon(\hat{\gamma}_{\hat{t}},\hat{\eta}_{\hat{t}})+
                         2{\beta}^{\frac{1}{3}}|\hat{\gamma}_{\hat{t}}(\hat{t})-\hat{\eta}_{\hat{t}}(\hat{t})|^2\nonumber\\
                         &
                         \leq&|W_1(\hat{\gamma}_{\hat{t}})-W_1(\hat{\eta}_{\hat{t}})|
                                   +|W_2(\hat{\gamma}_{\hat{t}})-W_2(\hat{\eta}_{\hat{t}})|+
                                   \sum_{i=0}^{\infty}\frac{1}{2^i}[\Upsilon(\eta^i_{t_i},\hat{\gamma}_{\hat{t}})+\Upsilon(\gamma^i_{t_i},\hat{\eta}_{\hat{t}})].
 \end{eqnarray}
  On the other hand,  notice that
 $$
                \Upsilon(\gamma_t,\eta_s)=\Upsilon(\gamma_t-\eta_{s,t}), \  \ \gamma_t,\eta_s\in \Lambda, \ 0\leq s\leq t\leq T,
 $$
  by Lemma \ref{theoremS000} and   the property (i) of $(\hat{t},(\hat{\gamma}_{\hat{t}},\hat{\eta}_{\hat{t}}))$,
 \begin{eqnarray}\label{4.7jiajia130}
 &&\sum_{i=0}^{\infty}\frac{1}{2^i}[\Upsilon(\eta^i_{t_i},\hat{\gamma}_{\hat{t}})+\Upsilon(\gamma^i_{t_i},\hat{\eta}_{\hat{t}})]\nonumber\\
 &\leq&2^5\sum_{i=0}^{\infty}\frac{1}{2^i}[\Upsilon(\eta^i_{t_i},\hat{\eta}_{\hat{t}})
 +\Upsilon(\gamma^i_{t_i},\hat{\gamma}_{\hat{t}})+2\Upsilon(\hat{\gamma}_{\hat{t}},\hat{\eta}_{\hat{t}})]
 \leq\frac{2^6}{\beta}+{2^7}\Upsilon(\hat{\gamma}_{\hat{t}},\hat{\eta}_{\hat{t}}).
 \end{eqnarray}
  Combining (\ref{5.6}) and (\ref{4.7jiajia130}),  from  (\ref{w}) and (\ref{5.10jiajiaaaa})  we have
 \begin{eqnarray}\label{5.jia6}
                         &&(2{\beta}-2^7)\Upsilon(\hat{\gamma}_{\hat{t}},\hat{\eta}_{\hat{t}})+
                          2\beta^{\frac{1}{3}}|\hat{\gamma}_{\hat{t}}(\hat{t})-\hat{\eta}_{\hat{t}}(\hat{t})|^2
                         \nonumber\\
                         &\leq& |W_1(\hat{\gamma}_{\hat{t}})-W_1(\hat{\eta}_{\hat{t}})|
                                   +|W_2(\hat{\gamma}_{\hat{t}})-W_2(\hat{\eta}_{\hat{t}})|+\frac{2^6}{\beta}\nonumber\\
                                   &\leq& 2L(2+||\hat{\gamma}_{\hat{t}}||_0+||\hat{\eta}_{\hat{t}}||_0)+\frac{2^6}{\beta}
                                   \leq 4L(1+M_0)+\frac{2^6}{\beta}.
 \end{eqnarray}
   Letting $\beta\rightarrow\infty$, we have
                \begin{eqnarray*}
                \Upsilon(\hat{\gamma}_{\hat{t}},\hat{\eta}_{\hat{t}})
                \leq \frac{1}{2{\beta}-2^7}\left[4L(1+M_0)+\frac{2^6}{\beta}\right]\rightarrow0\ 
                          \mbox{as} \ \beta\rightarrow\infty.
                         \end{eqnarray*}
In view of  (\ref{s0}), we have
\begin{eqnarray}\label{5.66666123}
||\hat{\gamma}_{\hat{t}}-\hat{\eta}_{\hat{t}}||_0 \rightarrow0\ 
                          \mbox{as} \ \beta\rightarrow\infty.
 \end{eqnarray}
  From  (\ref{s0}),  (\ref{w1}), (\ref{5.6}), (\ref{4.7jiajia130}) and (\ref{5.66666123}), we conclude
                           that
 \begin{eqnarray}\label{5.10112345}
                        &&{\beta}||\hat{\gamma}_{\hat{t}}-\hat{\eta}_{\hat{t}}||^6_{0}+\beta^{\frac{1}{3}}|\hat{\gamma}_{\hat{t}}(\hat{t})-\hat{\eta}_{\hat{t}}(\hat{t})|^2
                        \leq{\beta}\Upsilon(\hat{\gamma}_{\hat{t}},\hat{\eta}_{\hat{t}})
                        +\beta^{\frac{1}{3}}|\hat{\gamma}_{\hat{t}}(\hat{t})-\hat{\eta}_{\hat{t}}(\hat{t})|^2\nonumber\\
                        &\leq&\frac{1}{2}[|W_1(\hat{\gamma}_{\hat{t}})-W_1(\hat{\eta}_{\hat{t}})|
                                   +|W_2(\hat{\gamma}_{\hat{t}})-W_2(\hat{\eta}_{\hat{t}})|]+\frac{2^5}{\beta}+{2^6}\Upsilon(\hat{\gamma}_{\hat{t}},\hat{\eta}_{\hat{t}})\nonumber\\
                        &\leq& L||\hat{\gamma}_{\hat{t}}-\hat{\eta}_{\hat{t}}||_0
                                   +\frac{2^5}{\beta}+{2^{8}}||\hat{\gamma}_{\hat{t}}-\hat{\eta}_{\hat{t}}||_{0}^6
                                   \rightarrow0 \ 
                                   \mbox{as} \ \beta\rightarrow\infty.
 \end{eqnarray}
       Multiply the leftmost and rightmost sides of inequality (\ref{5.10112345}) by $\beta^{\frac{1}{6}}$, we obtain that
   \begin{eqnarray}\label{5.10123445567890}
                      \beta^{\frac{1}{2}}|\hat{{\gamma}}_{{\hat{t}}}(\hat{t})-\hat{{\eta}}_{{\hat{t}}}(\hat{t})|^2
                     \leq
                               L\beta^{\frac{1}{6}}||\hat{\gamma}_{\hat{t}}-\hat{\eta}_{\hat{t}}||_0
                                   +\frac{2^5}{\beta^{\frac{5}{6}}}+{2^{8}}\beta^{\frac{1}{6}}||\hat{\gamma}_{\hat{t}}-\hat{\eta}_{\hat{t}}||_{0}^6.
 \end{eqnarray}
By also (\ref{5.10112345}), the right side of above inequality converges to 0 as $\beta\rightarrow\infty$. Then we have that
 \begin{eqnarray*}
                      \beta^{\frac{1}{2}}|\hat{{\gamma}}_{{\hat{t}}}(\hat{t})-\hat{{\eta}}_{{\hat{t}}}(\hat{t})|^2
                     \rightarrow0 \ 
                                   \mbox{as} \ \beta\rightarrow\infty.
 \end{eqnarray*}
Combining with (\ref{5.10112345}), we have (\ref{5.10}).
 \par
   $Step\ 3.$ There exists 
   $N>0$
    such that
                   $\hat{t}\in [\tilde{t},T)$
                for all $\beta\geq N$.
 \par
By (\ref{5.66666123}), we can let $N>0$ be a large number such that
$$
                         L||\hat{\gamma}_{\hat{t}}-\hat{\eta}_{\hat{t}}||_0
                         \leq
                         \frac{\tilde{m}}{4},
$$
               for all $\beta\geq N$.
            Then we have $\hat{t}\in [\tilde{t},T)$ for all $\beta\geq N$. Indeed, if $\hat{t}=T$,  we  deduce the following contradiction:
 \begin{eqnarray*}
                         \frac{\tilde{m}}{2}\leq\Psi(\hat{\gamma}_{\hat{t}},\hat{\eta}_{\hat{t}})\leq \phi(\hat{\gamma}_{\hat{t}})-\phi(\hat{\eta}_{\hat{t}})\leq
                        L||\hat{\gamma}_{\hat{t}}-\hat{\eta}_{\hat{t}}||_0
                         \leq
                         \frac{\tilde{m}}{4}.
 \end{eqnarray*}
%
%
%
%
\par
 $Step\ 4.$   Crandall-Ishii lemma.
\par
          From above all,  for the fixed   $N>0$ in step 3, 
          we  find
$(\hat{t},\hat{\gamma}_{\hat{t}}), (\hat{t},\hat{\eta}_{\hat{t}})\in [\tilde{t}, T]\times
                 \Lambda^{\tilde{t}}$   satisfying $\hat{t}\in [\tilde{t},T)$  for all $\beta\geq N$
           such that
\begin{eqnarray}\label{psi4}
            \Psi_1(\hat{\gamma}_{\hat{t}},\hat{\eta}_{\hat{t}})\geq \Psi(\tilde{\gamma}_{\tilde{t}},\tilde{\gamma}_{\tilde{t}}) \ \ \mbox{and} \ \      \Psi_1(\hat{\gamma}_{\hat{t}},\hat{\eta}_{\hat{t}})\geq
                   {\Psi}_1(\gamma_t,\eta_t),
                   \  (\gamma_t,\eta_t)\in  \Lambda^{\hat{t}}\otimes \Lambda^{\hat{t}}.
\end{eqnarray}
                  We define,  
  for $(t,\gamma_t,\eta_t)\in [0,T]\times {\Lambda}\times {\Lambda}$,
\begin{eqnarray}\label{06091}
                             w_{1}(\gamma_t)=W_1(\gamma_t)-2^5\beta \Upsilon(\gamma_{t},\hat{\xi}_{\hat{t}})-\varepsilon\frac{\nu T-t}{\nu
                 T}\Upsilon(\gamma_t)
                 -\varepsilon \overline{\Upsilon}(\gamma_t,\hat{\gamma}_{\hat{t}})
               -\sum_{i=0}^{\infty}
        \frac{1}{2^i}\overline{\Upsilon}(\gamma^i_{t_i},\gamma_t),
        \end{eqnarray}
        \begin{eqnarray}\label{06092}
                             w_{2}(\eta_t)=-W_2(\eta_t)-2^5\beta \Upsilon(\eta_{t},\hat{\xi}_{\hat{t}})-\varepsilon\frac{\nu T-t}{\nu
                 T}\Upsilon(\eta_t)
                 -\varepsilon \overline{\Upsilon}(\eta_t,\hat{\eta}_{\hat{t}})-\sum_{i=0}^{\infty}
        \frac{1}{2^i}{\Upsilon}(\eta^i_{t_i},\eta_t),
\end{eqnarray}
   where $\hat{\xi}_{\hat{t}}=\frac{\hat{\gamma}_{\hat{t}}+\hat{\eta}_{\hat{t}}}{2}$.   We  note that $w_1,w_2$ depend on $\hat{\xi}_{{\hat{t}}}$, and thus on $\beta$ and
              $\varepsilon$. We also note that the last term in (\ref{06092}) is $\sum_{i=0}^{\infty}
        \frac{1}{2^i}{\Upsilon}(\eta^i_{t_i},\eta_t)$  rather than $\sum_{i=0}^{\infty}
        \frac{1}{2^i}\overline{\Upsilon}(\eta^i_{t_i},\eta_t)$. This is because  we divide the term $\sum_{i=0}^{\infty}
        \frac{1}{2^i}[{\Upsilon}(\gamma^i_{t_i},\gamma_t)+{\Upsilon}(\eta^i_{t_i},\eta_t)+|t-t_i|^2]$ in $\Psi_1$  into two terms $\sum_{i=0}^{\infty}
        \frac{1}{2^i}\overline{\Upsilon}(\gamma^i_{t_i},\gamma_t)$ and $\sum_{i=0}^{\infty}
        \frac{1}{2^i}{\Upsilon}(\eta^i_{t_i},\eta_t)$. Define $\varphi\in C^2(\mathbb{R}^d\times \mathbb{R}^d)$ by
              \begin{eqnarray}\label{07063}
                             \varphi(x,y)=\beta^{\frac{1}{3}}|x-y|^2,\ \ (x,y)\in \mathbb{R}^d\times \mathbb{R}^d.
             \end{eqnarray}
By the following Lemma \ref{0611a}, $w_1$ and $w_2$ satisfy the conditions of  Theorem \ref{theorem0513} with $\varphi$ defined by (\ref{07063}).
 Then by Theorem \ref{theorem0513}, there exist
sequences  $(l_{k},\check{\gamma}^{k}_{l_k}), (s_{k},\check{\eta}^{k}_{s_k})\in [\hat{t},T]\times \Lambda^{\hat{t}}$ and
 sequences of functionals $\varphi_k\in C_p^{1,2}(\Lambda^{l_k}),\psi_k\in C_p^{1,2}(\Lambda^{s_k})$ bounded from below such that 
\begin{eqnarray}\label{0609a}
 w_{1}(\gamma_t)-\varphi_k(\gamma_t)
\end{eqnarray}
has a strict  maximum $0$ at  $\check{\gamma}^{k}_{l_k}$ over $\Lambda^{l_k}$,
\begin{eqnarray}\label{0609b}
w_{2}(\eta_t)-\psi_k(\eta_{t})
\end{eqnarray}
has a strict  maximum $0$ at $\check{\eta}^{k}_{s_k}$over $ \Lambda^{s_k}$, and
\begin{eqnarray}\label{0608v1}
       &&\left(l_{k}, \check{\gamma}^{k}_{l_{k}}(l_{k}), w_1(\check{\gamma}^{k}_{l_{k}}),\partial_t\varphi_k( \check{\gamma}^{k}_{l_{k}}),\partial_x\varphi_k( \check{\gamma}^{k}_{l_{k}}),\partial_{xx}\varphi_k( \check{\gamma}^{k}_{l_{k}})\right)\nonumber\\
       &&\underrightarrow{k\rightarrow\infty}\left({\hat{t}},\hat{{\gamma}}_{{\hat{t}}}(\hat{t}), w_1(\hat{{\gamma}}_{{\hat{t}}}), b_1, 2\beta^{\frac{1}{3}} (\hat{{\gamma}}_{{\hat{t}}}({{\hat{t}}})-\hat{{\eta}}_{{\hat{t}}}({{\hat{t}}})), X\right),
\end{eqnarray}
\begin{eqnarray}\label{0608vw1}
       &&\left(s_{k}, \check{\eta}^{k}_{s_{k}}(s_{k}), w_2(\check{\eta}^{k}_{s_{k}}),\partial_t\psi_k( \check{\eta}^{k}_{s_{k}}), \partial_x\psi_k( \check{\eta}^{k}_{s_{k}}), \partial_{xx}\psi_k( \check{\eta}^{k}_{s_{k}})\right)\nonumber\\
       &&\underrightarrow{k\rightarrow\infty}\left({\hat{t}},\hat{{\eta}}_{{\hat{t}}}(\hat{t}), w_2(\hat{{\eta}}_{{\hat{t}}}), b_2, 2\beta^{\frac{1}{3}}(\hat{{\eta}}_{{\hat{t}}}({{\hat{t}}})-\hat{{\gamma}}_{{\hat{t}}}({{\hat{t}}})), Y\right),
\end{eqnarray}
 where $b_{1}+b_{2}=0$ and $X,Y\in \mathcal{S}(\mathbb{R}^{d})$ satisfy the following inequality:
\begin{eqnarray}\label{II}
                              {-6\beta^{\frac{1}{3}}}\left(\begin{array}{cc}
                                    I&0\\
                                    0&I
                                    \end{array}\right)\leq \left(\begin{array}{cc}
                                    X&0\\
                                    0&Y
                                    \end{array}\right)\leq  6\beta^{\frac{1}{3}} \left(\begin{array}{cc}
                                    I&-I\\
                                    -I&I
                                    \end{array}\right).
\end{eqnarray}
            We note that  
              sequence  $(\check{\gamma}^{k}_{l_k},\check{\eta}^{k}_{s_k},l_{k},s_{k},\varphi_k,\psi_k)$ and $b_{1},b_{2},X,Y$  depend on  $\beta$ and
              $\varepsilon$. We also note that  (\ref{II}) follows from (\ref{II0615}) choosing $\kappa=\frac{1}{2}\beta^{-\frac{1}{3}}$. In fact, by (\ref{07063}),
              $$
                           A= \nabla_x^2\varphi(\hat{\gamma}(0),\hat{\eta}(0))=2\beta^{\frac{1}{3}}\left(\begin{array}{cc}
                                    I&-I\\
                                    -I&I
                                    \end{array}\right),
              $$
              and thus, if $\kappa=\frac{1}{2}\beta^{-\frac{1}{3}}$,
               $$A+\kappa A^2=(1+4\kappa \beta^{\frac{1}{3}})A=3A,$$ and
              $$-\left(\frac{1}{\kappa}+|A|\right)=-\left(2\beta^{\frac{1}{3}}+4\beta^{\frac{1}{3}}\right)=-6\beta^{\frac{1}{3}}.$$
               Then from (\ref{II0615}) it follows that  (\ref{II}) holds true.
  By the following Lemma \ref{lemma4.4}, we have
\begin{eqnarray}\label{4.23}
\lim_{k\rightarrow\infty}[d_\infty(\check{\gamma}^{k}_{l_k},\hat{\gamma}_{\hat{t}})
+d_\infty(\check{\eta}^{k}_{s_k},\hat{\eta}_{\hat{t}})]=0.
\end{eqnarray}
 For every $(t,\gamma_t),(s,\eta_s)\in [T-\bar{a},T]\times{{\Lambda}^{T-\bar{a}}}$, let
\begin{eqnarray*}
\chi^{k}(\gamma_t)
        :=
                \varepsilon\frac{\nu T-t}{\nu
                 T}\Upsilon(\gamma_t)
                +\varepsilon \overline{\Upsilon}(\gamma_t,\hat{\gamma}_{\hat{t}})
                +\sum_{i=0}^{\infty}
        \frac{1}{2^i}\overline{\Upsilon}(\gamma^i_{t_i},\gamma_t)+2^5\beta\Upsilon(\gamma_{t},\hat{{\xi}}_{{\hat{t}}})
                +\varphi_k(\gamma_t),
  \end{eqnarray*}
  \begin{eqnarray*}
\hbar^{k}(\eta_s)
        :=-\varepsilon\frac{\nu T-s}{\nu
                 T}\Upsilon(\eta_s)
                 -\varepsilon \overline{\Upsilon}(\eta_s,\hat{\eta}_{\hat{t}})
                 -\sum_{i=0}^{\infty}
        \frac{1}{2^i}{\Upsilon}(\eta^i_{t_i},\eta_s)-2^5\beta\Upsilon(\eta_{s},\hat{{\xi}}_{{\hat{t}}})
        -\psi_k(\eta_s).
  \end{eqnarray*}
  It is clear that  $\chi^{k}(\cdot)\in C^{1,2}_p(\Lambda^{l_k}),\hbar^{k}(\cdot)\in C^{1,2}_p(\Lambda^{s_k})$. 
  Moreover, by  (\ref{0609a}), (\ref{0609b}) and definitions of $w_1$ and $w_2$,
  $$
                         (W_1-\chi^{k})(\check{\gamma}^{k}_{l_k})=\sup_{(t,\gamma_t)\in [l_k,T]\times\Lambda^{l_k}}
                         (W_1-\chi^{k})(\gamma_t),
$$
$$
                         (W_2-\hbar^{k})(\check{\eta}^{k}_{s_k})=\inf_{(s,\eta_s)\in [s_k,T]\times\Lambda^{s_k}}
                         (W_2-\hbar^{k})(\eta_s).
$$
From $l_{k}\rightarrow {\hat{t}}$, $s_{k}\rightarrow {\hat{t}}$ as $k\rightarrow\infty$ and ${\hat{t}}<T$
 for $\beta>N$, it follows that, for every fixed $\beta>N$,   constant $ K_\beta>0$ exists such that
$$
            l_{k}\vee s_{k}<T, 
            \ \ \mbox{for all}    \ \ k\geq K_\beta.
$$
Now, for every $\beta>N$ and $k>K_\beta$, from the definition of viscosity solutions it follows that
  \begin{eqnarray}\label{vis1}
                      \partial_t\chi^{k}(\check{\gamma}^{k}_{l_{k}})
  +{\mathbf{H}}{(}\check{\gamma}^{k}_{l_{k}},  W_1(\check{\gamma}^{k}_{l_{k}}), \partial_x\chi^{k}(\check{\gamma}^{k}_{l_{k}}),
                                       \partial_{xx}\chi^{k}(\check{\gamma}^{k}_{l_{k}})
                                    {)}\geq c,
 \end{eqnarray}
 and
  \begin{eqnarray}\label{vis2}
                     \partial_t\hbar^{k}(\check{\eta}^{k}_{s_{k}})+{\mathbf{H}}{(}\check{\eta}^{k}_{s_{k}}, W_2(\check{\eta}^{k}_{s_{k}}),
                     \partial_x\hbar^{k}(\check{\eta}^{k}_{s_{k}}),\partial_{xx}\hbar^{k}(\check{\eta}^{k}_{s_{k}}){)}\leq0,
  \end{eqnarray}
 where, for every $(t,\gamma_t)\in [l_k,T]\times{{\Lambda}^{l_k}}$ and  $ (s,\eta_s)\in [s_k,T]\times{{\Lambda}^{s_k}}$, from 
    Lemma \ref{theoremS},
  \begin{eqnarray*}
\partial_t\chi^{k}(\gamma_t)=
                       \partial_t\varphi_k(\gamma_{{t}})-\frac{\varepsilon}{\nu T}\Upsilon(\gamma_{{t}})
                     +2\varepsilon({t}-{\hat{t}})+2\sum_{i=0}^{\infty}\frac{1}{2^i}(t-t_{i}),
  \end{eqnarray*}
  \begin{eqnarray*}
\partial_x\chi^{k}(\gamma_t)&=&\partial_{x}\varphi_k(\gamma_{{t}})
                     +\varepsilon\frac{\nu T-{t}}{\nu T}\partial_x\Upsilon(\gamma_{{t}}) +\varepsilon\partial_x\Upsilon(\gamma_{{t}}-\hat{\gamma}_{{\hat{t},t}})
+2^5\beta\partial_x\Upsilon(\gamma_{{t}}-\hat{\xi}_{{\hat{t},t}})
\\
 &&+\partial_x\left[\sum_{i=0}^{\infty}\frac{1}{2^i}
                      \Upsilon(\gamma_{t}-\gamma^{i}_{t_{i},t})
                     \right],
  \end{eqnarray*}
   \begin{eqnarray*}
\partial_{xx}\chi^{k}(\gamma_t)&=&\partial_{xx}(\varphi_k)(\gamma_{{t}})+\varepsilon\frac{\nu T-{t}}{\nu T}\partial_{xx}\Upsilon(\gamma_{{t}}) +\varepsilon\partial_{xx}\Upsilon(\gamma_{{t}}-\hat{\gamma}_{{\hat{t},t}})
+2^5\beta\partial_{xx}\Upsilon(\gamma_{{t}}-\hat{\xi}_{{\hat{t},t}})
\\
                      &&
                      +\partial_{xx}\left[\sum_{i=0}^{\infty}\frac{1}{2^i}
                      \Upsilon(\gamma_{t}-\gamma^{i}_{t_{i},t})
                      \right],
  \end{eqnarray*}
   \begin{eqnarray*}
\partial_t\hbar^{k}(\eta_s)=
                      -\partial_t\psi_k(\eta_{{s}})+\frac{\varepsilon}{\nu T}\Upsilon(\eta_{{s}})-2\varepsilon({s}-{\hat{t}}),
  \end{eqnarray*}
  \begin{eqnarray*}
\partial_x\hbar^{k}(\eta_s)
                          &=& -\partial_{x}\psi_k(\eta_{{s}})
                      -\varepsilon\frac{\nu T-{s}}{\nu T}
                     \partial_x\Upsilon(\eta_{{s}})-\varepsilon\partial_x\Upsilon(\eta_{{s}}-\hat{\eta}_{{\hat{t}},s})
                     -2^5\beta\partial_x\Upsilon(\eta_{{s}}-\hat{\xi}_{{\hat{t}},s})\\
                     &&
                      -\partial_x\left[\sum_{i=0}^{\infty}\frac{1}{2^i}
                      \Upsilon(\eta_{s}-{\eta}^{i}_{{t}_{i},s})
                      \right],
  \end{eqnarray*}
   \begin{eqnarray*}
\partial_{xx}\hbar^{k}(\eta_s)
                                &=&-\partial_{xx}\psi_k( \eta_{{s}})-\varepsilon\frac{\nu T-{s}}{\nu T}
                     \partial_{xx}\Upsilon(\eta_{{s}})-\varepsilon\partial_{xx}\Upsilon(\eta_{{s}}-\hat{\eta}_{{\hat{t}},s})-2^5\beta\partial_{xx}\Upsilon(\eta_{{s}}-\hat{\xi}_{{\hat{t}},s})\\
                     &&
                      -\partial_{xx}\left[\sum_{i=0}^{\infty}\frac{1}{2^i}
                      \Upsilon(\eta_{s}-{\eta}^{i}_{{t}_{i},s})
                      \right].
  \end{eqnarray*}
 \par
 $Step\ 5.$    Calculation and completion of the proof.
 \par
  We notice that,  by  (\ref{220817a}) and  (\ref{220817a1}), 
   there exists a generic constant $C>0$ such that
    \begin{eqnarray*}
    |\partial_x\Upsilon(\check{\gamma}^{k}_{l_k}-\hat{\gamma}_{{\hat{t}},l_k})|
    +|\partial_x\Upsilon(\check{\eta}^{k}_{s_k}-\hat{\eta}_{{\hat{t}},s_k})|\leq C|\hat{\gamma}_{{\hat{t}}}(\hat{t})-\check{\gamma}^{k}_{l_k}(l_k)|^5
                      +C|\hat{\eta}_{{\hat{t}}}(\hat{t})-\check{\eta}^{k}_{s_k}(s_k)|^5;
    \end{eqnarray*}
    \begin{eqnarray*}
    |\partial_{xx}\Upsilon(\check{\gamma}^{k}_{l_k}-\hat{\gamma}_{{\hat{t}},l_k})|
    +|\partial_{xx}\Upsilon(\check{\eta}^{k}_{s_k}-\hat{\eta}_{{\hat{t}},s_k})|\leq C|\hat{\gamma}_{{\hat{t}}}(\hat{t})-\check{\gamma}^{k}_{l_k}(l_k)|^4
                      +C|\hat{\eta}_{{\hat{t}}}(\hat{t})-\check{\eta}^{k}_{s_k}(s_k)|^4.
    \end{eqnarray*}
    Letting  $k\rightarrow\infty$ in (\ref{vis1}) and (\ref{vis2}), and using (\ref{0608v1}), (\ref{0608vw1}) and (\ref{4.23}),    we obtain
    \begin{eqnarray}\label{03103}
                   b_1-\frac{\varepsilon}{\nu T}\Upsilon(\hat{{\gamma}}_{{\hat{t}}})+2\sum_{i=0}^{\infty}\frac{1}{2^i}(\hat{t}-t_i)
                     +{\mathbf{H}}(\hat{{\gamma}}_{{\hat{t}}},W_1(\hat{{\gamma}}_{{\hat{t}}}),
                    \partial_x\chi(\hat{\gamma}_{\hat{t}}),\partial_{xx}\chi(\hat{\gamma}_{\hat{t}}))
                     \geq c;
\end{eqnarray}
and
 \begin{eqnarray}\label{03104}
                     -b_2+ \frac{\varepsilon}{\nu T}\Upsilon(\hat{{\eta}}_{{\hat{t}}})+{\mathbf{H}}(\hat{{\eta}}_{{\hat{t}}},W_2(\hat{{\eta}}_{{\hat{t}}}),\partial_x\hbar(\hat{\eta}_{{\hat{t}}}),\partial_{xx}\hbar(\hat{\eta}_{\hat{t}}))
                     \leq0,
\end{eqnarray}
  where
  \begin{eqnarray*}
\partial_x\chi(\hat{\gamma}_{\hat{t}})
                                     :=
2\beta^{\frac{1}{3}}(\hat{{\gamma}}_{{\hat{t}}}({{\hat{t}}})-\hat{{\eta}}_{{\hat{t}}}({{\hat{t}}}))+2^5\beta\partial_x\Upsilon(\hat{\gamma}_{{\hat{t}}}-\hat{\xi}_{{\hat{t}}})
                                      +\varepsilon\frac{\nu T-{\hat{t}}}{\nu T}\partial_x\Upsilon(\hat{{\gamma}}_{{\hat{t}}})
                                    +\partial_x\left[\sum_{i=0}^{\infty}\frac{1}{2^i}\Upsilon(\hat{\gamma}_{\hat{t}}-\gamma^i_{t_i,\hat{t}})\right],
  \end{eqnarray*}
   \begin{eqnarray*}
\partial_{xx}\chi(\hat{\gamma}_{\hat{t}})
:= X+2^5\beta\partial_{xx}\Upsilon(\hat{\gamma}_{{\hat{t}}}-\hat{\xi}_{{\hat{t}}})
                                     +\varepsilon\frac{\nu T-{\hat{t}}}{\nu T}\partial_{xx}\Upsilon(\hat{{\gamma}}_{{\hat{t}}})
                                     +\partial_{xx}\left[\sum_{i=0}^{\infty}\frac{1}{2^i}\Upsilon(\hat{\gamma}_{\hat{t}}-\gamma^i_{t_i,\hat{t}})\right],
  \end{eqnarray*}
  \begin{eqnarray*}
\partial_x\hbar(\hat{\eta}_{{\hat{t}}})
                                     :=
2\beta^{\frac{1}{3}}(\hat{{\gamma}}_{{\hat{t}}}({{\hat{t}}})-\hat{{\eta}}_{{\hat{t}}}({{\hat{t}}}))-2^5\beta\partial_x\Upsilon(\hat{\eta}_{{\hat{t}}}-\hat{\xi}_{{\hat{t}}})-\varepsilon\frac{\nu T-{\hat{t}}}{\nu T}\partial_x\Upsilon(\hat{{\eta}}_{{\hat{t}}})-\partial_x\left[\sum_{i=0}^{\infty}\frac{1}{2^i}
                                     \Upsilon(\hat{\eta}_{\hat{t}}-\eta^i_{t_i,\hat{t}})\right],
  \end{eqnarray*}
  and
   \begin{eqnarray*}
\partial_{xx}\hbar(\hat{\eta}_{\hat{t}})
:= -Y-2^5\beta\partial_{xx}\Upsilon(\hat{\eta}_{{\hat{t}}}-\hat{\xi}_{{\hat{t}}})-\varepsilon\frac{\nu T-{\hat{t}}}{\nu T}\partial_{xx}\Upsilon(\hat{{\eta}}_{{\hat{t}}})
                      -\partial_{xx}\left[\sum_{i=0}^{\infty}\frac{1}{2^i}
                      \Upsilon(\hat{\eta}_{\hat{t}}-\eta^i_{t_i,\hat{t}})\right].
  \end{eqnarray*}
  Notice that $b_1+b_2=0$ and $\hat{\xi}_{{\hat{t}}}=\frac{\hat{{\gamma}}_{\hat{t}}+\hat{{\eta}}_{\hat{t}}}{2}$,
combining  (\ref{03103}) and (\ref{03104}),  we have
 \begin{eqnarray}\label{vis112}
                    && c+ \frac{\varepsilon}{\nu T}(\Upsilon(\hat{{\gamma}}_{{\hat{t}}})+\Upsilon(\hat{{\eta}}_{{\hat{t}}})
                     )-2\sum_{i=0}^{\infty}\frac{1}{2^i}(\hat{t}-t_i)\nonumber\\
                     &\leq&{\mathbf{H}}(\hat{{\gamma}}_{{\hat{t}}},W_1(\hat{{\gamma}}_{{\hat{t}}}),\partial_x\chi(\hat{\gamma}_{\hat{t}}),\partial_{xx}\chi(\hat{\gamma}_{\hat{t}}))
                     -{\mathbf{H}}(\hat{{\eta}}_{{\hat{t}}},W_2(\hat{{\eta}}_{{\hat{t}}}),\partial_x\hbar(\hat{\eta}_{{\hat{t}}}),\partial_{xx}\hbar(\hat{\eta}_{\hat{t}})).
\end{eqnarray}
 On the other hand, by (\ref{5.1}) and via a simple calculation we obtain
 \begin{eqnarray}\label{v4}
                &&{\mathbf{H}}(\hat{{\gamma}}_{{\hat{t}}},W_1(\hat{{\gamma}}_{{\hat{t}}}),\partial_x\chi(\hat{\gamma}_{\hat{t}}),\partial_{xx}\chi(\hat{\gamma}_{\hat{t}}))
                     -{\mathbf{H}}(\hat{{\eta}}_{{\hat{t}}},W_2(\hat{{\eta}}_{{\hat{t}}}),\partial_x\hbar(\hat{\eta}_{{\hat{t}}}),\partial_{xx}\hbar(\hat{\eta}_{\hat{t}}))\nonumber\\
                     &\leq&{\mathbf{H}}(\hat{{\gamma}}_{{\hat{t}}},W_2(\hat{{\eta}}_{{\hat{t}}}),\partial_x\chi(\hat{\gamma}_{\hat{t}}),\partial_{xx}\chi(\hat{\gamma}_{\hat{t}}))
                     -{\mathbf{H}}(\hat{{\eta}}_{{\hat{t}}},W_2(\hat{{\eta}}_{{\hat{t}}}),\partial_x\hbar(\hat{\eta}_{{\hat{t}}}),\partial_{xx}\hbar(\hat{\eta}_{\hat{t}}))\nonumber\\
                &\leq&\sup_{u\in U}(J_{1}+J_{2}+J_{3}),
\end{eqnarray}
  where from Hypothesis \ref{hypstate}, (\ref{220817a}), (\ref{220817a1}) 
    and (\ref{II}),  we have
\begin{eqnarray}\label{j1}
                               J_{1}&=&\langle {b}(\hat{{\gamma}}_{{\hat{t}}},u),\partial_x\chi(\hat{\gamma}_{\hat{t}})\rangle  -\langle {b}(\hat{{\eta}}_{{\hat{t}}},u),\partial_x\hbar(\hat{\eta}_{{\hat{t}}})\rangle\nonumber\\
                                  &\leq&2\beta^{\frac{1}{3}}{L}|\hat{{\gamma}}_{{\hat{t}}}({\hat{t}})-\hat{{\eta}}_{{\hat{t}}}({\hat{t}})|\times
                                 ||\hat{{\gamma}}_{{\hat{t}}}-\hat{{\eta}}_{{\hat{t}}}||_0
                                 +18\beta|\hat{{\gamma}}_{{\hat{t}}}({\hat{t}})-\hat{{\eta}}_{{\hat{t}}}({\hat{t}})|^5L(2+||\hat{{\gamma}}_{{\hat{t}}}||_0
                                 +||\hat{{\eta}}_{{\hat{t}}}||_0)
                                \nonumber\\
                                            &&
                                            +18L\sum_{i=0}^{\infty}\frac{1}{2^i}\left[|\gamma^i_{t_i}(t_i)-\hat{\gamma}_{\hat{t}}(\hat{t})|^5
                                            +|\eta^i_{t_i}(t_i)-\hat{\eta}_{\hat{t}}(\hat{t})|^5\right]
                                            (1+||\hat{{\gamma}}_{{\hat{t}}}||_0+||\hat{{\eta}}_{{\hat{t}}}||_0)\nonumber\\
                                           && +36\varepsilon \frac{\nu T-{\hat{t}}}{\nu T} L(1+||\hat{{\gamma}}_{{\hat{t}}}||^6_0+||\hat{{\eta}}_{{\hat{t}}}||^6_0);
\end{eqnarray}
\begin{eqnarray}\label{j2}
                               J_{2}&=&\frac{1}{2}\mbox{tr}{[}\partial_{xx}\chi(\hat{\gamma}_{\hat{t}}){\sigma}(\hat{{\gamma}}_{{\hat{t}}},u)
                                        {\sigma}^\top(\hat{{\gamma}}_{{\hat{t}}},u){]}-\frac{1}{2}\mbox{tr}{[}\partial_{xx}\hbar(\hat{\eta}_{\hat{t}})\sigma(\hat{{\eta}}_{{\hat{t}}},u)\sigma^\top(\hat{{\eta}}_{{\hat{t}}},u){]}\nonumber\\
                               &\leq&3
                               \beta^{\frac{1}{3}}|{\sigma}(\hat{{\gamma}}_{{\hat{t}}},u)-\sigma(\hat{{\eta}}_{{\hat{t}}},u)|_2^2
                                      +306\beta|\hat{{\gamma}}_{{\hat{t}}}({\hat{t}})-\hat{{\eta}}_{{\hat{t}}}({\hat{t}})|^4(|{\sigma}(\hat{{\gamma}}_{{\hat{t}}},u)|_2^2+
                                      |\sigma(\hat{{\eta}}_{{\hat{t}}},u)|_2^2)\nonumber\\
                                      &&+153\varepsilon\frac{\nu T-{\hat{t}}}{\nu    T}(|\hat{{\gamma}}_{{\hat{t}}}({\hat{t}})|^4|{\sigma}(\hat{{\gamma}}_{{\hat{t}}},u)|_2^2+|\hat{{\eta}}_{{\hat{t}}}({\hat{t}})|^4|{\sigma}(\hat{{\eta}}_{{\hat{t}}},u)|_2^2)
                                      \nonumber\\
                                      &&
                                      +153\sum_{i=0}^{\infty}\frac{1}{2^i}|\gamma^i_{t_i}(t_i)-\hat{\gamma}_{\hat{t}}(\hat{t})|^4
                                        |\sigma(\hat{{\gamma}}_{{\hat{t}}},u)|_2^2
                                        +153\sum_{i=0}^{\infty}\frac{1}{2^i}|\eta^i_{t_i}(t_i)-\hat{\eta}_{\hat{t}}(\hat{t})|^4
                                        |\sigma(\hat{{\eta}}_{{\hat{t}}},u)|_2^2
                                        \nonumber\\
                               &\leq&
                             3
                               \beta^{\frac{1}{3}}{L^2}||\hat{{\gamma}}_{{\hat{t}}}-\hat{{\eta}}_{{\hat{t}}}||_0^2+ 306\beta L^2|\hat{{\gamma}}_{{\hat{t}}}({\hat{t}})-\hat{{\eta}}_{{\hat{t}}}({\hat{t}})|^4(2+||\hat{{\gamma}}_{{\hat{t}}}||_0^2
                                             +||\hat{{\eta}}_{{\hat{t}}}||_0^2
                                             )\nonumber\\
                                             &&+306\varepsilon \frac{\nu T-{\hat{t}}}{\nu T}L^2(1+||\hat{{\gamma}}_{{\hat{t}}}||^6_0
                                 +||\hat{{\eta}}_{{\hat{t}}}||^6_0) \nonumber\\
                                             &&+153\left(
                                     \sum_{i=0}^{\infty}\frac{1}{2^i}\left[\left|\gamma^i_{t_i}(t_i)-\hat{\gamma}_{\hat{t}}(\hat{t})\right|^4
                                     +\left|\eta^i_{t_i}(t_i)-\hat{\eta}_{\hat{t}}(\hat{t})\right|^4\right]\right)L^2
                                     (1+||\hat{{\gamma}}_{{\hat{t}}}||_0^2
                                             +||\hat{{\eta}}_{{\hat{t}}}||_0^2
                                             );
\end{eqnarray}
and
\begin{eqnarray}\label{j3}
                                 J_{3}&=&q{(}\hat{{\gamma}}_{{\hat{t}}}, W_2(\hat{{\eta}}_{{\hat{t}}}), \sigma^\top(\hat{{\gamma}}_{{\hat{t}}},u)\partial_x\chi(\hat{\gamma}_{\hat{t}}),u{)}-
                                 q{(}\hat{{\eta}}_{{\hat{t}}}, W_2(\hat{{\eta}}_{{\hat{t}}}), \sigma(\hat{{\eta}}_{{\hat{t}}},u)^\top\partial_x\hbar(\hat{\eta}_{{\hat{t}}}),u{)}\nonumber\\
                                     &\leq&
                                 L||\hat{{\gamma}}_{{\hat{t}}}-\hat{{\eta}}_{{\hat{t}}}||_0
                               +2\beta^{\frac{1}{3}} L^2|\hat{{\gamma}}_{{\hat{t}}}({\hat{t}})
                                 -\hat{{\eta}}_{{\hat{t}}}({\hat{t}})|\times||\hat{{\gamma}}_{{\hat{t}}}-\hat{{\eta}}_{{\hat{t}}}||_0
                                 +18\beta L^2|\hat{{\gamma}}_{{\hat{t}}}({\hat{t}})-\hat{{\eta}}_{{\hat{t}}}({\hat{t}})|^5
                                 (2+||\hat{{\gamma}}_{{\hat{t}}}||_0
                                             +||\hat{{\eta}}_{{\hat{t}}}||_0
                                             )\nonumber\\
                               &&+18L^2\sum_{i=0}^{\infty}\frac{1}{2^i}\left[\left|\gamma^i_{t_i}(t_i)-\hat{\gamma}_{\hat{t}}(\hat{t})\right|^5
                                            +\left|\eta^i_{t_i}(t_i)-\hat{\eta}_{\hat{t}}(\hat{t})\right|^5\right]
                                            (1+||\hat{{\gamma}}_{{\hat{t}}}||_0+||\hat{{\eta}}_{{\hat{t}}}||_0)\nonumber\\
                                            &&+36\varepsilon \frac{\nu T-{\hat{t}}}{\nu T} L^2(1+||\hat{{\gamma}}_{{\hat{t}}}||_0^6
                                             +||\hat{{\eta}}_{{\hat{t}}}||_0^6
                                             ).
\end{eqnarray}
   We notice that, by  the property (i) of $(\hat{t},(\hat{\gamma}_{\hat{t}},\hat{\eta}_{\hat{t}}))$,
  \begin{eqnarray*}
2\sum_{i=0}^{\infty}\frac{1}{2^i}(\hat{t}-t_i)
  \leq2\sum_{i=0}^{\infty}\frac{1}{2^i}\bigg{(}\frac{1}{2^i\beta}\bigg{)}^{\frac{1}{2}}\leq 4{\bigg{(}\frac{1}{{\beta}}\bigg{)}}^{\frac{1}{2}},
    \end{eqnarray*}
     \begin{eqnarray*}
  \sum_{i=0}^{\infty}\frac{1}{2^i}\left[\left|\gamma^i_{t_i}(t_i)-\hat{\gamma}_{\hat{t}}(\hat{t})\right|^5
                                            +\left|\eta^i_{t_i}(t_i)-\hat{\eta}_{\hat{t}}(\hat{t})\right|^5\right]
                      \leq
                      2\sum_{i=0}^{\infty}\frac{1}{2^i}\bigg{(}\frac{1}{2^i\beta}\bigg{)}^{\frac{5}{6}}\leq 4{\bigg{(}\frac{1}{{\beta}}\bigg{)}}^{\frac{5}{6}},
                        \end{eqnarray*}
                        and
                         \begin{eqnarray*}
  \sum_{i=0}^{\infty}\frac{1}{2^i}\left[\left|\gamma^i_{t_i}(t_i)-\hat{\gamma}_{\hat{t}}(\hat{t})\right|^4
                                     +\left|\eta^i_{t_i}(t_i)-\hat{\eta}_{\hat{t}}(\hat{t})\right|^4\right]
                      \leq
                      2\sum_{i=0}^{\infty}\frac{1}{2^i}\bigg{(}\frac{1}{2^i\beta}\bigg{)}^{\frac{2}{3}}\leq 4{\bigg{(}\frac{1}{{\beta}}\bigg{)}}^{\frac{2}{3}}.
  \end{eqnarray*}
 Combining (\ref{vis112})-(\ref{j3}), then by (\ref{5.10jiajiaaaa}) and (\ref{5.10}) we can let $\beta>0$ be large enough such that
\begin{eqnarray}\label{vis122}
                     c
                                             \leq
                      -\frac{\varepsilon}{\nu T}(\Upsilon(\hat{{\gamma}}_{{\hat{t}}})
                     +\Upsilon(\hat{{\eta}}_{{\hat{t}}}))+ \varepsilon \frac{\nu T-{\hat{t}}}{\nu T} (342L+36)L(1+||\hat{{\gamma}}_{{\hat{t}}}||_0^6
                                             +||\hat{{\eta}}_{{\hat{t}}}||_0^6
                                             )+\frac{c}{4}.
\end{eqnarray}
Recalling $
         \nu=1+\frac{1}{2T(342L+36)L}
$ and $\bar{a}=\frac{1}{2(342L+36)L}\wedge{T}$, by (\ref{s0}) and (\ref{5.3}), the following contradiction is induced:
\begin{eqnarray*}\label{vis122}
                     c\leq
                             \frac{\varepsilon}{\nu
                              T}+\frac{c}{4}\leq \frac{c}{2}.
\end{eqnarray*}
 The proof is now complete.
 \ \ $\Box$
  \par
 To complete the previous proof, it remains to state and prove the following lemmas. In the following Lemmas of this section, let $\tilde{w}_{1}^{\hat{t}}, \tilde{w}_{1}^{\hat{t},*}$ and $\tilde{w}_{2}^{\hat{t}}, \tilde{w}_{2}^{\hat{t},*}$ be the definitions in Definition \ref{definition0607} with respect to $w_1$ defined by (\ref{06091}) and  $w_2$ defined by  (\ref{06092}), respectively.
 \begin{lemma}\label{lemma4.3}\ \
  There exists a local modulus of continuity  $\rho_1$ 
  such that  the functionals $w_1$ and $w_2$ defined by (\ref{06091}) and (\ref{06092}) satisfy condition (\ref{0608a}).
\end{lemma}
\par
  {\bf  Proof}. \ \
From (\ref{w1}) and the definition of $w_{1}$, we have that,
for every $\hat{t}\leq t\leq s\leq T$ and $\gamma_t\in \Lambda^{\hat{t}}$,
\begin{eqnarray*}
&&w_1(\gamma_t)-w_1(\gamma_{t,s})\\&=&W_1(\gamma_t)-2^5\beta \Upsilon(\gamma_{t},\hat{\xi}_{\hat{t}})-\varepsilon\frac{\nu T-t}{\nu
                 T}\Upsilon(\gamma_t)
                 -\varepsilon \overline{\Upsilon}(\gamma_t,\hat{\gamma}_{\hat{t}})
               -\sum_{i=0}^{\infty}
        \frac{1}{2^i}\overline{\Upsilon}(\gamma^i_{t_i},\gamma_t)\\
        &&-W_1(\gamma_{t,s})+2^5\beta \Upsilon(\gamma_{t,s},\hat{\xi}_{\hat{t}})+\varepsilon\frac{\nu T-s}{\nu
                 T}\Upsilon(\gamma_{t,s})
                 +\varepsilon \overline{\Upsilon}(\gamma_{t,s},\hat{\gamma}_{\hat{t}})
               +\sum_{i=0}^{\infty}
        \frac{1}{2^i}\overline{\Upsilon}(\gamma^i_{t_i},\gamma_{t,s})\\
        &=&W_1(\gamma_t)-W_1(\gamma_{t,s})+\varepsilon\frac{ t-s}{\nu
                 T}\Upsilon(\gamma_t)+\varepsilon((s-\hat{t})^2-(t-\hat{t})^2)+\sum_{i=0}^{\infty}
        \frac{1}{2^i}((s-t_i)^2-(t-t_i)^2)\\
                 &\leq& \rho_2(|s-t|,||\gamma_t||_0)+2T(2+\varepsilon)|s-t|.
\end{eqnarray*}
Taking $\rho_1(l,x)= \rho_2(l,x)+2T(2+\varepsilon)l$,  
  $(l,x)\in[0,\infty)\times [0,\infty)$, it is clear that $\rho_1$ is a local modulus of continuity 
   and $w_1$  satisfies condition (\ref{0608a}) with it.  In a similar way, we show that  $w_2$  satisfies condition (\ref{0608a}) with this $\rho_1$.
The proof is now complete. \ \ $\Box$
 \begin{lemma}\label{lemma4.344}\ \ $\tilde{w}_{1}^{\hat{t},*}\in \Phi(\hat{t},\hat{{\gamma}}_{{\hat{t}}}({\hat{t}}))$ and $\tilde{w}_{2}^{\hat{t},*}\in \Phi(\hat{t},\hat{{\eta}}_{{\hat{t}}}({\hat{t}}))$. 
\end{lemma}
\par
  {\bf  Proof}. \ \  We only prove $\tilde{w}_{1}^{\hat{t},*}\in \Phi(\hat{t},\hat{{\gamma}}_{{\hat{t}}}({\hat{t}}))$. 
                           $\tilde{w}_{2}^{\hat{t},*}\in \Phi(\hat{t},\hat{{\eta}}_{{\hat{t}}}({\hat{t}}))$ can be obtained by a symmetric  way.
 Set $r=\frac{1}{2}(|T-{\hat{t}}|\wedge \hat{t})
 $, for  given $L>0$, let $\varphi\in C^{1,2}([0,T]\times \mathbb{R}^{d})$ be a function such that
            $\tilde{w}_{1}^{\hat{t},*}(t,x)-\varphi(t,x)$  has a  maximum at $(\bar{{t}},\bar{x})\in (0, T)\times \mathbb{R}^{d}$, moreover, the following inequalities hold true:
\begin{eqnarray}\label{0815zhou1}
                    |\bar{{t}}-{\hat{t}}|+|\bar{x}-\hat{{\gamma}}_{{\hat{t}}}({\hat{t}})|<r=\frac{1}{2}(|T-{\hat{t}}|\wedge \hat{t}),
                    \end{eqnarray}
             \begin{eqnarray}\label{0815zhou2}
                    |\tilde{w}_{1}^{\hat{t},*}(\bar{{t}},\bar{{x}})|+|\nabla_x\varphi(\bar{{t}},\bar{{x}})|
                    +|\nabla^2_x\varphi(\bar{{t}},\bar{{x}})|\leq L.
\end{eqnarray}
By  Lemma 5.4 of Chapter 4 in \cite{yong}, we can modify $\varphi$  such  that $\varphi\in C^{1,2}([0,T]\times \mathbb{R}^d)$ bounded from below, ${\varphi}$, ${\varphi_t}$, $\nabla_{x}{\varphi}$ and $\nabla^2_{x}{\varphi}$  grow  in a polynomial way,
   $\tilde{w}_{1}^{\hat{t},*}(t,x)-\varphi(t,x)$  has a strict   maximum 0 at $(\bar{{t}},\bar{x})\in (0, T)\times \mathbb{R}^{d}$ on $[0, T]\times \mathbb{R}^{d}$ and the above two inequalities hold true.
If $\bar{{t}}<\hat{t}$, we have $\varphi_{t}(\bar{{t}},\bar{x})=\frac{1}{2}(\hat{t}-\bar{{t}})^{-\frac{1}{2}} > 0$.
%
%
%
%
%
%
%
%
If $\bar{{t}}\geq \hat{t}$,
 recall that $w_1$ is defined in (\ref{06091}), we consider the functional
$$
                \Gamma(\gamma_t)= w_{1}(\gamma_t)
                 -\varphi(t,\gamma_t(t)),\ (t,\gamma_t)\in [\hat{t},T]\times {{\Lambda}}.
$$
Note that $w_1$ is a continuous  functional bounded from above and $\varphi$ is a continuous  function bounded from below, 
 $\Gamma$ is a continuous  functional bounded from above on $\Lambda^{\hat{t}}$.
 Define a sequence of positive numbers $\{\delta_i\}_{i\geq0}$  by 
        $\delta_i=\frac{1}{2^i}$ for all $i\geq0$.  For every  $0<\delta<1$,
 by Lemma \ref{theoremleft} we have that,
 for every  $(\breve{t}_0,\breve{\gamma}^0_{\breve{t}_0})\in [\bar{t},T]\times \Lambda^{\bar{t}}$ satisfying
 \begin{eqnarray}\label{0615a}
\Gamma(\breve{\gamma}^0_{\breve{t}_0})\geq \sup_{(s,\gamma_s)\in [\bar{t},T]\times \Lambda^{\bar{t}}}\Gamma(\gamma_s)-\delta,
 \end{eqnarray}
  there exist $(\breve{t},\breve{\gamma}_{\breve{t}})\in [\bar{t},T]\times \Lambda^{\bar{t}}$ and a sequence $\{(\breve{t}_i,\breve{\gamma}^i_{\breve{t}_i})\}_{i\geq1}\subset
  [\breve{t}_0,T]\times \Lambda^{\bar{t}}$ such that
  \begin{description}
        \item{(i)} $\overline{\Upsilon}(\breve{\gamma}^0_{\breve{t}_0},\breve{\gamma}_{\breve{t}})\leq \delta$,
         $\overline{\Upsilon}(\breve{\gamma}^i_{\breve{t}_i},\breve{\gamma}_{\breve{t}})\leq \frac{\delta}{2^i}$ and $t_i\uparrow \breve{t}$ as $i\rightarrow\infty$,
        \item{(ii)}  $\Gamma(\breve{\gamma}_{\breve{t}})
            -\sum_{i=0}^{\infty}\frac{1}{2^i}\overline{\Upsilon}(\breve{\gamma}^i_{\breve{t}_i},\breve{\gamma}_{\breve{t}})
        \geq \Gamma(\breve{\gamma}^0_{\breve{t}_0})$, and
        \item{(iii)}    for all $(s,\gamma_s)\in [\breve{t},T]\times \Lambda^{\breve{t}}\setminus \{(\breve{t},\breve{\gamma}_{\breve{t}})\}$,
        \begin{eqnarray*}
        \Gamma(\gamma_s)
        -\sum_{i=0}^{\infty}
        \frac{1}{2^i}\overline{\Upsilon}(\breve{\gamma}^i_{\breve{t}_i},\gamma_s)
            <\Gamma(\breve{\gamma}_{\breve{t}})
            -\sum_{i=0}^{\infty}\frac{1}{2^i}\overline{\Upsilon}(\breve{\gamma}^i_{\breve{t}_i},\breve{\gamma}_{\breve{t}}).
        \end{eqnarray*}
        \end{description}
        We should note that the point
             $(\breve{t},\breve{\gamma}_{\breve{t}})$ depends on $\delta$.
        By the definitions of $\tilde{w}^{\hat{t}}_1$ and $\tilde{w}_{1}^{\hat{t},*}$, we have
 \begin{eqnarray}\label{220901a}
        &&\tilde{w}_{1}^{\hat{t},*}(\bar{{t}},\bar{x})-\varphi(\bar{{t}},\bar{x})\nonumber\\
        &=&\limsup_{s\geq \hat{t},(s,y)\rightarrow(\bar{{t}},\bar{x})}(\tilde{w}^{\hat{t}}_1(s,y)-\varphi(s,y))
        =\limsup_{s\geq \hat{t}, (s,y)\rightarrow(\bar{{t}},\bar{x})}\left(\sup_{\xi_s\in \Lambda^{\hat{t}},\xi_s(s)=y}
                             w_1(\xi_s)-\varphi(s,y)\right). \ \ \ \
 \end{eqnarray}
  Note that by Lemma \ref{lemma4.3}, $w_1$   satisfies condition (\ref{0608a}). Then, for every $(s,\xi_s)\in [\hat{t},\bar{t}]\times \Lambda$,
 \begin{eqnarray}\label{220901b}
                   w_1(\xi_s)\leq w_1(\xi_{s,\bar{t}})+\rho_1(|\bar{t}-s|,||\xi_s||_0).
 \end{eqnarray}
  By the definition of $w_1$, there exists a constant {$M_4>0$} such that
  \begin{eqnarray}\label{220901c}
  \sup_{\xi_s\in \Lambda^{\hat{t}},\xi_s(s)=y}
                             w_1(\xi_s)=\sup_{\xi_s\in \Lambda^{\hat{t}},\xi_s(s)=y, ||\xi_s||_0\leq {M_4}}
                             w_1(\xi_s).
 \end{eqnarray}
  Thus, by (\ref{220901b}) and (\ref{220901c}),
  \begin{eqnarray}\label{220901d}
  &&\limsup_{s\geq \bar{t}, (s,y)\rightarrow(\bar{{t}},\bar{x})}\sup_{\xi_s\in \Lambda^{\hat{t}},\xi_s(s)=y}
                             [w_1(\xi_s)-\varphi(s,y)]\nonumber\\
                             &\leq& \limsup_{s\geq \hat{t}, (s,y)\rightarrow(\bar{{t}},\bar{x})}\left(\sup_{\xi_s\in \Lambda^{\hat{t}},\xi_s(s)=y 
                             }
                             w_1(\xi_s)-\varphi(s,y)\right)\nonumber\\
                             &\leq& \limsup_{s\geq \hat{t}, (s,y)\rightarrow(\bar{{t}},\bar{x})}\sup_{\xi_s\in \Lambda^{\hat{t}},\xi_s(s)=y, ||\xi_s||_0\leq {M_4}}
                             {[}w_1(\xi_{s,s\vee\bar{t}})+\rho_1(|s\vee\bar{t}-s|,{M_4})-\varphi(s,y){]}\nonumber\\
                             &\leq&\limsup_{s\geq \bar{t}, (s,y)\rightarrow(\bar{{t}},\bar{x})}\sup_{\xi_s\in \Lambda^{\hat{t}},\xi_s(s)=y}
                             [w_1(\xi_s)-\varphi(s,y)].
  \end{eqnarray}
  Therefore, by (\ref{220901a}) and (\ref{220901d}),
  \begin{eqnarray*}
        \tilde{w}_{1}^{\hat{t},*}(\bar{{t}},\bar{x})-\varphi(\bar{{t}},\bar{x})
                             =\limsup_{s\geq \bar{t}, (s,y)\rightarrow(\bar{{t}},\bar{x})}\sup_{\xi_s\in \Lambda^{\hat{t}},\xi_s(s)=y}
                             [w_1(\xi_s)-\varphi(s,y)]\leq\sup_{(s,\gamma_s)\in [\bar{t},T]\times \Lambda^{\bar{t}}}\Gamma(\gamma_s).
 \end{eqnarray*}
 Combining with (\ref{0615a}),
 \begin{eqnarray}\label{220901e}
 \Gamma(\breve{\gamma}^0_{\breve{t}_0})\geq \sup_{(s,\gamma_s)\in [\bar{t},T]\times \Lambda^{\bar{t}}}\Gamma(\gamma_s)-\delta
\geq\tilde{w}_{1}^{\hat{t},*}(\bar{{t}},\bar{x})-\varphi(\bar{{t}},\bar{x})-\delta.
 \end{eqnarray}
  Recall that $\tilde{w}_{1}^{\hat{t},*}\geq\tilde{w}^{\hat{t}}_1$.  Then, by 
   the definition of $\tilde{w}^{\hat{t}}_1$, the property (ii) of $(\breve{t},\breve{\gamma}_{\breve{t}})$ and (\ref{220901e}),
\begin{eqnarray}\label{20210509}
                 \tilde{w}_{1}^{\hat{t},*}(\breve{t},\breve{\gamma}_{\breve{t}}(\breve{t}))
                   -\varphi(\breve{t},\breve{\gamma}_{\breve{t}}(\breve{t}))
                   &\geq&\tilde{w}^{\hat{t}}_1(\breve{t},\breve{\gamma}_{\breve{t}}(\breve{t}))
                   -\varphi(\breve{t},\breve{\gamma}_{\breve{t}}(\breve{t}))
                 \geq w_1(\breve{\gamma}_{\breve{t}})
                -\varphi(\breve{t},\breve{\gamma}_{\breve{t}}(\breve{t}))\nonumber\\
                  &\geq& \Gamma(\breve{\gamma}^0_{\breve{t}_0})
                 \geq\tilde{w}_{1}^{\hat{t},*}(\bar{{t}},\bar{x})-\varphi(\bar{{t}},\bar{x})-\delta=-\delta.
\end{eqnarray}
Noting that $\nu$ is independent of  $\delta$ and $\varphi$ is a continuous  function bounded from below, 
 by the definitions of  $\Gamma$ and $w_1$, (\ref{w}) and (\ref{s0})
 there exists a constant  $M_5>0$   depending only on $\varphi$ 
  that is sufficiently  large   that
$
           \Gamma(\gamma_t)<\sup_{(s,\gamma_s)\in [\bar{t},T]\times \Lambda^{\bar{t}}}\Gamma(\gamma_s)-1
           $ for all $t\in [\bar{t},T]$ and $||\gamma_t||_0\geq M_5$. Thus, we have $||\breve{\gamma}_{\breve{t}}||_0\vee
           ||{\breve{\gamma}}^{0}_{\breve{t}_{0}}||_0<M_5$. In particular, $|\breve{\gamma}_{\breve{t}}(\breve{t})|<M_5$.  
   %
%
Letting $\delta\rightarrow0$, by  the similar proof procedure of (\ref{4.22}) and (\ref{05231}), we obtain
\begin{eqnarray}\label{delta0}
       && \breve{t}\rightarrow \bar{t},\ \breve{\gamma}_{\breve{t}}(\breve{t})\rightarrow \bar{x},\  \tilde{w}_{1}^{\hat{t},*}(\breve{t},\breve{\gamma}_{\breve{t}}(\breve{t}))\rightarrow\tilde{w}_{1}^{\hat{t},*}(\bar{{t}},\bar{x}) \ \mbox{as}\ \delta\rightarrow0.
       \end{eqnarray}
       Noting that $\tilde{w}_{1}^{\hat{t},*}(t,x)-\varphi(t,x)$  has a strict   maximum 0 at $(\bar{{t}},\bar{x})\in (0, T)\times \mathbb{R}^{d}$ on $[0, T]\times \mathbb{R}^{d}$, by (\ref{0815zhou2})
        and (\ref{delta0}) there exists a constant $0<\Delta<1$ such that
       for all $0<\delta<\Delta$,
       $$
       \varphi(\breve{t},\breve{\gamma}_{\breve{t}}(\breve{t}))\geq \tilde{w}_{1}^{\hat{t},*}(\breve{t},\breve{\gamma}_{\breve{t}}(\breve{t}))\geq \tilde{w}_{1}^{\hat{t},*}(\bar{{t}},\bar{x})-1\geq -(L+1).
       $$
       Then, by the definitions of  $\Gamma$ and $w_1$,  (\ref{s0}), (\ref{w}) and (\ref{20210509}),
 there exists a constant  $M_6>0$   depending only on $L$ 
 such that, for all $0<\delta<\Delta$ and $||{\gamma}_{\breve{t}}||_0\geq M_6$ satisfying ${\gamma}_{\breve{t}}(\breve{t})=\breve{\gamma}_{\breve{t}}(\breve{t})$,
$$
           \Gamma(\gamma_{\breve{t}})<-2\leq \Gamma(\breve{\gamma}^0_{\breve{t}_0})-1 \leq\sup_{(s,\gamma_s)\in [\bar{t},T]\times \Lambda^{\bar{t}}}\Gamma(\gamma_s)-1.
$$
Thus, we have $||\breve{\gamma}_{\breve{t}}||_0<M_6$ for all $0<\delta<\Delta$.
    %
%
  From (\ref{0815zhou1}), it follows that $\bar{t}<{\hat{t}}+\frac{|T-{\hat{t}}|}{2}\leq T$. 
 Then, by (\ref{delta0}),  we have  $\breve{t}<T$ 
  provided that  $\delta>0$ is small enough.
      Thus,  the definition of the viscosity subsolution can be used to obtain the following result:
 \begin{eqnarray}\label{5.15}
                      \partial_t\Im(\breve{\gamma}_{{\breve{t}}})
                     +{\mathbf{H}}{(}\breve{\gamma}_{{\breve{t}}},  W_1(\breve{\gamma}_{{\breve{t}}}), \partial_x\Im(\breve{\gamma}_{{\breve{t}}}),\partial_{xx}\Im(\breve{\gamma}_{{\breve{t}}}){)}\geq c.
 \end{eqnarray}
  where, for every $(t,\gamma_t)\in [\bar{t},T]\times{{\Lambda}}$,
  \begin{eqnarray*}
\Im(\gamma_t)
        :=\varepsilon\frac{\nu T-t}{\nu
                 T}\Upsilon(\gamma_t)
                +\varepsilon \overline{\Upsilon}(\gamma_t,\hat{\gamma}_{\hat{t}})
                +\sum_{i=0}^{\infty}
        \frac{1}{2^i}\overline{\Upsilon}(\gamma^i_{t_i},\gamma_t)+2^5\beta\Upsilon(\gamma_{t},\hat{{\xi}}_{{\hat{t}}})
                +\sum_{i=0}^{\infty}
        \frac{1}{2^i}\overline{\Upsilon}(\breve{\gamma}^{i}_{\breve{t}_{i}},\gamma_t)+\varphi(t,\gamma_t(t)),
  \end{eqnarray*}
 \begin{eqnarray*}
\partial_t\Im(\gamma_t):=-\frac{\varepsilon}{\nu T}\Upsilon(\gamma_{{t}})
                     +2\varepsilon({t}-{\hat{t}})+2\sum_{i=0}^{\infty}\frac{1}{2^i}[(t-t_{i})+(t-\breve{t}_i)]+\varphi_t({t},\gamma_{{t}}({t})),
  \end{eqnarray*}
  \begin{eqnarray*}
\partial_x\Im(\gamma_t)&:=&
                     \varepsilon\frac{\nu T-{t}}{\nu T}\partial_x\Upsilon(\gamma_{{t}}) +\varepsilon\partial_x\Upsilon(\gamma_{{t}}-\hat{\gamma}_{{\hat{t},t}})
+2^5\beta\partial_x\Upsilon(\gamma_{{t}}-\hat{\xi}_{{\hat{t},t}})
\\
 &&+\partial_x\left[\sum_{i=0}^{\infty}\frac{1}{2^i}
                      \Upsilon(\gamma_{t}-\gamma^{i}_{t_{i},t})
                      +\sum_{i=0}^{\infty}\frac{1}{2^i}\Upsilon(\gamma_{t}-\breve{\gamma}^i_{\breve{t}_i,t})\right]+\nabla_{x}\varphi({t}, \gamma_{{t}}(t)),
  \end{eqnarray*}
   \begin{eqnarray*}
\partial_{xx}\Im(\gamma_t)&:=&\varepsilon\frac{\nu T-{t}}{\nu T}\partial_{xx}\Upsilon(\gamma_{{t}}) +\varepsilon\partial_{xx}\Upsilon(\gamma_{{t}}-\hat{\gamma}_{{\hat{t},t}})+2^5\beta\partial_{xx}\Upsilon(\gamma_{{t}}-\hat{\xi}_{{\hat{t},t}})
\\
                      &&
                      +\partial_{xx}\left[\sum_{i=0}^{\infty}\frac{1}{2^i}
                      \Upsilon(\gamma_{t}-\gamma^{i}_{t_{i},t})
                      +\sum_{i=0}^{\infty}\frac{1}{2^i}\Upsilon(\gamma_{t}-\breve{\gamma}^i_{\breve{t}_i,t})\right]+\nabla^2_{x}\varphi({t}, \gamma_{{t}}(t)).
  \end{eqnarray*}
  We notice that  $||\breve{\gamma}_{{\breve{t}}}||_0< M_6$ for all $0<\delta<\Delta$ and $M_6$ only depends on $L$.
Then letting $\delta\rightarrow0$ in (\ref{5.15}),   by the definition of ${\mathbf{H}}$, (\ref{0815zhou2}) and (\ref{delta0}), it follows that there exists a constant $C_0\geq0$ depending only on $L$ such that $ \varphi_{t}(\bar{t},\bar{x})
\geq -C_0$.
The proof is now complete. \ \ $\Box$
 \begin{lemma}\label{0611a}\ \
 The functionals $w_1$ and $w_2$ defined by (\ref{06091}) and (\ref{06092}) satisfy the conditions of  Theorem \ref{theorem0513} with  $\varphi$, where $\varphi$ is the function defined by (\ref{07063}).
\end{lemma}
\par
  {\bf  Proof}. \ \
From (\ref{s0}) and (\ref{w}),  $w_1$ and $w_2$ are continuous functionals bounded from above and satisfy (\ref{05131}).
By  Lemmas \ref{lemma4.3} and \ref{lemma4.344}, $w_1$ and $w_2$ satisfy condition  (\ref{0608a}), and  $\tilde{w}_{1}^{\hat{t},*}\in \Phi(\hat{t},\hat{{\gamma}}_{{\hat{t}}}({\hat{t}}))$ and $\tilde{w}_{2}^{\hat{t},*}\in \Phi(\hat{t},\hat{{\eta}}_{{\hat{t}}}({\hat{t}}))$.
Moreover, let $\varphi$ be the function defined by (\ref{07063}). By Lemma \ref{theoremS000} and (\ref{iii4}) we obtain that, for all 
   $(t,(\gamma_t,\eta_t))\in [\hat{t},T]\times (\Lambda^{\hat{t}}\otimes  \Lambda^{\hat{t}})$,
\begin{eqnarray}\label{wv}
                         &&w_{1}(\gamma_t)+w_{2}(\eta_t)-\varphi(\gamma_t(t),\eta_t(t))=w_{1}(\gamma_t)+w_{2}(\eta_t)-\beta^{\frac{1}{3}}|\gamma_t(t)-\eta_t(t)|^2\nonumber   \\
                          &\leq&
                   {\Psi}_1(\gamma_t,\eta_t)\leq\Psi_1(\hat{\gamma}_{\hat{t}},\hat{\eta}_{\hat{t}})=w_{1}(\hat{{\gamma}}_{{\hat{t}}})+w_{2}(\hat{{\eta}}_{{\hat{t}}})
                          -\beta^{\frac{1}{3}}|\hat{{\gamma}}_{{\hat{t}}}({\hat{t}})-\hat{{\eta}}_{{\hat{t}}}({\hat{t}})|^2\nonumber\\
                          &=&w_{1}(\hat{{\gamma}}_{{\hat{t}}})+w_{2}(\hat{{\eta}}_{{\hat{t}}})
                          -\varphi(\hat{{\gamma}}_{{\hat{t}}}({\hat{t}}),\hat{{\eta}}_{{\hat{t}}}({\hat{t}})),
\end{eqnarray}
                       where the last inequality becomes equality if and only if $t={\hat{t}}$, $\gamma_t=\hat{{\gamma}}_{{\hat{t}}}, \eta_t=\hat{{\eta}}_{{\hat{t}}}$.
                         Then we obtain that $
                w_1(\gamma_t)+w_2(\eta_t)-\varphi(\gamma_t(t),\eta_t(t))
$
has a 
 maximum over $\Lambda^{\hat{t}}\otimes \Lambda^{\hat{t}}$ at a point $(\hat{\gamma}_{\hat{t}},\hat{\eta}_{\hat{t}})$ with $\hat{t}\in (0,T)$. Thus $w_1$ and $w_2$ satisfy the conditions of  Theorem \ref{theorem0513}  with  $\varphi$ defined by (\ref{07063}).  \ \ $\Box$
\begin{lemma}\label{lemma4.4}\ \ The maximum points $(\check{\gamma}^{k}_{l_{k}}, \check{\eta}^{k}_{s_{k}})$
 satisfy  condition (\ref{4.23}).
\end{lemma}
\par
  {\bf  Proof}. \ \
 Without loss of generality, we may assume $s_{k}\leq l_{k}$, by  (\ref{up}), (\ref{w1}), (\ref{iii4}) and the definitions of $w_1$ and $w_2$, we have that
\begin{eqnarray}\label{06101}
                      &&w_1(\check{\gamma}^{k}_{l_{k}})+w_2(\check{\eta}^{k}_{s_{k}})
               -\beta^{\frac{1}{3}} |\check{\gamma}^{k}_{l_{k}}(l_{k})-\check{\eta}^{k}_{s_{k}}(s_{k})|^2\nonumber\\
                      &\leq&\Psi_1(\check{\gamma}^{k}_{l_{k}},\check{\eta}^k_{s_k,l_k})  -W_2(\check{\eta}^{k}_{s_{k}})
                      +W_2(\check{\eta}^k_{s_k,l_k})
                       -\varepsilon[\overline{\Upsilon}(\check{\gamma}^{k}_{l_{k}},\hat{\gamma}_{\hat{t}})
+\overline{\Upsilon}(\check{\eta}^k_{s_k},\hat{\eta}_{\hat{t}})
]\nonumber\\
                       &\leq&\Psi_1(\hat{\gamma}_{\hat{t}},\hat{\eta}_{\hat{t}})
                    +\rho_2(|l_k- s_k|,||\check{\eta}^k_{s_k}||_0) 
                      -\varepsilon[\overline{\Upsilon}(\check{\gamma}^{k}_{l_{k}},\hat{\gamma}_{\hat{t}})
+\overline{\Upsilon}(\check{\eta}^k_{s_k},\hat{\eta}_{\hat{t}})
].
\end{eqnarray}
By (\ref{0608vw1}), $w_2(\check{\eta}^k_{s_k})\rightarrow w_2(\hat{\eta}_{\hat{t}})$ as $k\rightarrow\infty$. Then by  that $w_2$ satisfies condition (\ref{05131}),  
 there exists a constant  $M_7>0$  that is sufficiently  large   that
$$
 ||\check{\eta}^k_{s_k}||_0\leq M_7,\ \mbox{for all} \ k>0.
$$
Letting $k\rightarrow\infty$ in (\ref{06101}), by (\ref{0608v1}) and (\ref{0608vw1}) we have that
\begin{eqnarray*}
                       \Psi_1(\hat{\gamma}_{\hat{t}},\hat{\eta}_{\hat{t}})=
                       w_1(\hat{\gamma}_{\hat{t}})+w_2(\hat{\eta}_{\hat{t}})
               -\beta^{\frac{1}{3}} |\hat{\gamma}_{\hat{t}}(\hat{t})-\hat{\eta}_{\hat{t}}(\hat{t})|^2
                      \leq \Psi_1(\hat{\gamma}_{\hat{t}},\hat{\eta}_{\hat{t}})
                            -\varepsilon\limsup_{k\rightarrow\infty}[\overline{\Upsilon}(\check{\gamma}^{k}_{l_{k}},\hat{\gamma}_{\hat{t}})
+\overline{\Upsilon}(\check{\eta}^k_{s_k},\hat{\eta}_{\hat{t}})
].
\end{eqnarray*}
Thus,
$$
\lim_{k\rightarrow\infty}[\overline{\Upsilon}(\check{\gamma}^{k}_{l_{k}},\hat{\gamma}_{\hat{t}})
+\overline{\Upsilon}(\check{\eta}^k_{s_k},\hat{\eta}_{\hat{t}})
]=0.
$$
Then by (\ref{s0}) we get
(\ref{4.23}) holds true.
  The proof is now complete. \ \ $\Box$
  \begin{remark}\label{remarks}
   The continuity condition (\ref{w1}) is crucial in our current approach. In fact, condition (\ref{w1}) is used to prove (\ref{5.10}) and  Lemmas \ref{lemma4.3} and \ref{lemma4.4}.  In particular, since $\partial_{xx}S^{ a_{\hat{t}}}(\cdot)$ is not equal to $\mathbf{0}$ (see Lemma \ref{theoremS}),
the convergence property (\ref{5.10}) is the key to prove  comparison theorem.  It will be very interesting to structure a new smooth function to see if we can apply it to prove comparison theorem without condition (\ref{w1}).
\end{remark}

\section{Application to BSHJB equations.}
\par
In this section, we
show that our PHJB equations includes backward stochastic HJB (BSHJB) equations as a special case (see also   \cite[Example 4.5]{tang1}). In the following, we let $n=d$.
\par
 We consider the  controlled state
             equation:
\begin{eqnarray}\label{state6}
          \bar{X}^{t,x,u}(s)=x+\int_{t}^{s}\bar{b}(W_l,\bar{X}^{t,x,u}(l),u(l))dl
          +\int_{t}^{s}\bar{\sigma}(W_l,\bar{X}^{t,x,u}(l),u(l))dW(l),\
             s\in [t,T],
\end{eqnarray}
and the associated BSDE:
\begin{eqnarray}\label{fbsde6}
\bar{Y}^{{t,x,u}}(s)&=&\bar{\phi}(W_T,\bar{X}^{t,x,u}(T))+\int^{T}_{s}\bar{q}(W_l,\bar{X}^{t,x,u}(l),
          \bar{Y}^{t,x,u}(l),\bar{Z}^{t,x,u}(l),u(l))dl\nonumber\\
                 &&-\int^{T}_{s}\bar{Z}^{t,x,u}(l)dW(l),\ \ s\in [t,T],
\end{eqnarray}
with $\bar{b}:{\Lambda}\times \mathbb{R}^m\times U\rightarrow \mathbb{R}^m$, $\bar{\sigma}:{\Lambda}\times \mathbb{R}^m\times U\rightarrow \mathbb{R}^{m\times d}$,
$\bar{q}:{\Lambda}\times \mathbb{R}^m\times \mathbb{R}\times \mathbb{R}^{d}\times U\rightarrow \mathbb{R}$ and $\bar{\phi}:{\Lambda}_T\times \mathbb{R}^m\rightarrow \mathbb{R}$.
The value functional of the optimal control is defined by
\begin{eqnarray}\label{value6}
\bar{V}(t,x):=\mathop{\esssup}\limits_{u(\cdot)\in{\cal{U}}[t,T]} \bar{Y}^{t,x,u}(t),\ \  (t,x)\in [0,T]\times \mathbb{R}^m.
\end{eqnarray}
This problem is path dependent on $\omega_t$ and state dependent on $\bar{X}(t)$. Now we transform this problem into the path-dependent case.
\par
In this section, for each  $t\in[0,T]$,
          define
         ${\Lambda}^{d+m}_t$ as  the set of  continuous  $\mathbb{R}^{d+m}$-valued
         functions on $[0,t]$.
        We denote ${\Lambda}^{d+m}=\bigcup_{t\in[0,T]}{\Lambda}^{d+m}_{t}$.
For any $(\omega_t,\xi_t),(\omega_T,\xi_T)\in {\Lambda}^{d+m}$, $(y,z)\in \mathbb{R}\times \mathbb{R}^{d}$ and $ u\in U$, we define
$b:{\Lambda}^{d+m}\times U\rightarrow \mathbb{R}^{d+m}$, $\sigma:{\Lambda}^{d+m}\times U\rightarrow \mathbb{R}^{(d+m)\times d}$,
$q:{\Lambda}^{d+m}\times \mathbb{R}\times \mathbb{R}^{d}\times U\rightarrow \mathbb{R}$ and $\phi:{\Lambda}_T^{d+m}\rightarrow \mathbb{R}$ as
\begin{eqnarray*}
                   b((\omega_t,\xi_t),u)&:=&\left(\begin{array}{cc}
                                    \mathbf{0}\\
                                    \bar{b}(\omega_t,\xi_t(t),u)
                                    \end{array}\right),\ \ \ \
                                     \sigma((\omega_t,\xi_t),u):=\left(\begin{array}{cc}
                                    {I}\\
                                    \bar{\sigma}(\omega_t,\xi_t(t),u)
                                    \end{array}\right),\\
                                     q((\omega_t,\xi_t),y,z,u)&:=&
                                    \bar{q}(\omega_t,\xi_t(t),y,z,u),\ \ \ \
                                    \ \ \ \ \ \ \   \phi(\omega_T,\xi_T)  :=   \bar{\phi}(\omega_T,\xi_T(T)).
\end{eqnarray*}
We assume ${b}, {\sigma}, {q}, {\phi}$ satisfy Hypothesis \ref{hypstate}, then following (\ref{state1}), (\ref{fbsde1}) and (\ref{value1}), for any $(\omega_t,\xi_t)\in {\Lambda}^{d+m}$ and $u(\cdot)\in {\cal{U}}[t,T]$
 we can define $X^{(\omega_t,\xi_t),u}$, $Y^{(\omega_t,\xi_t),u}$ and $V(\omega_t,\xi_t):=\mathop{\esssup}\limits_{u(\cdot)\in{\cal{U}}[t,T]}Y^{(\omega_t,\xi_t),u}(t)$.
 Noting $V(\omega_t,\xi_t)$ only depends on the state $x=\xi_t(t)$ of the path $\xi_t$ at time $t$, we can rewrite
 $X^{(\omega_t,\xi_t),u}$, $Y^{(\omega_t,\xi_t),u}$ and $V(\omega_t,\xi_t)$ into $X^{\omega_t,x,u}$, $Y^{\omega_t,x,u}$ and $V(\omega_t,x)$, respectively.
 Then, in view of Theorem \ref{theorem52},  $V(\omega_t, x)$ is a unique viscosity  solution to the
 PHJB  equation:
 \begin{eqnarray}\label{hjb6666}
\begin{cases}
\partial_tV(\omega_t,x)+\sup_{u\in{
                                         {U}}}[\langle\nabla_x V(\omega_t,x),\bar{b}(\omega_t,x,u)\rangle+
                                         \frac{1}{2}\mbox{tr}(\nabla^2_{x}V(\omega_t,x)\bar{\sigma}(\omega_t,x,u)\bar{\sigma}^\top(\omega_t,x,u))\\
                                         \ \ \ \ \ \ \ \ \ \ \ \ \ \ \ \ \ +
                                         \frac{1}{2}\mbox{tr}\partial_{\gamma\gamma}V(\omega_t,x)+
                                        \mbox{tr}(\bar{\sigma}^\top(\omega_t,x,u)\partial_{x\gamma}V(\omega_t,x))+
                                         \bar{q}(\omega_t,x,V(\omega_t,x),\partial_\gamma V(\omega_t,x)\\
                                          \ \ \ \ \ \ \ \ \ \ \ \ \ \ \ \ \  +\bar{\sigma}^\top(\omega_t,x,u)\nabla_xV(\omega_t,x),u)
                                         ]= 0,\ \ \  (t,x,\omega)\in
                               [0,T)\times \mathbb{R}^m\times \Omega,\\
  V(\omega_T,x)=\bar{\phi}(\omega_T,x),\ \ \  (x,\omega)\in
                               \mathbb{R}^m\times \Omega.
                               \end{cases}
\end{eqnarray}
Here, $\partial_\gamma$ and $\partial_{\gamma\gamma}$ are the spatial derivatives in $\gamma_t\in \Lambda$, and $\nabla_x$ and $\nabla^2_{x}$ are the classical partial derivatives in the state variable $x$.
\par
If $V(\omega_t,x)$ is smooth enough, applying functional It\^o formula to $V(W_t,x)$, we obtain
\begin{eqnarray*}
                             dV(W_t,x)=[\partial_tV(W_t,x)+\frac{1}{2}\mbox{tr}\partial_{\gamma\gamma}V(W_t,x)]dt+\partial_\gamma V(W_t,x)dW(t),  \ \  P\mbox{-}a.s.,
                             \ \  (t,x)\in
                               [0,T]\times \mathbb{R}^m.
\end{eqnarray*}
Define the pair of ${\cal{F}}_t$-adapted processes
\begin{eqnarray}\label{value62}
(\bar{V}(t,x), p(t,x)):=(V(W_t,x),\partial_\gamma V(W_t,x)), \ \ (t,x)\in
                               [0,T]\times \mathbb{R}^m,
\end{eqnarray}
and combine with (\ref{hjb6666}), we have
\begin{eqnarray}\label{hjb666611}
\begin{cases}
d\bar{V}(t,x)=-\sup_{u\in{
                                         {U}}}[\langle\nabla_x \bar{V}(t,x),\bar{b}(W_t,x,u)\rangle+
                                         \frac{1}{2}\mbox{tr}(\nabla^2_{x}\bar{V}(t,x)\bar{\sigma}(W_t,x,u)\bar{\sigma}^\top(W_t,x,u))\\
                                         ~~~~~~~~~~~~~~~~~~+
                                        \mbox{tr}(\bar{\sigma}^\top(W_t,x,u)\nabla_{x}p(t,x))+
                                         \bar{q}(W_t,x,\bar{V}(t,x),p(t,x)\cr
                                          ~~~~~~~~~~~~~~~~~~+\bar{\sigma}^\top(W_t,x,u)\nabla_x\bar{V}(t,x),u)
                                         ]+p(t,x)dW(t),\  (t,x)\in
                               [0,T]\times \mathbb{R}^m,  \ \ P\mbox{-}a.s.,\\
 \bar{V}(T,x)=\bar{\phi}(W_T,x), \ \ x\in \mathbb{R}^m, \  \ P\mbox{-}a.s.
 \end{cases}
\end{eqnarray}
Thus we obtain that
\begin{theorem}\label{theorem61}
If the value functional  $V\in C_{p}^{1,2}({\Lambda}^{d+m})$, then the pair of $\{{\cal{F}}_t\}_{t\geq0}$-adapted processes $(\bar{V}(t,x), p(t,x))$ defined by (\ref{value62}) is a
                         classical  solution to (\ref{hjb666611}).
 \end{theorem}
 Notice that $(\bar{V}(t,x), p(t,x))$ is only dependent on $V$ which is a unique viscosity solution to PHJB equation (\ref{hjb6666}), then we can give the definition of viscosity solutions to BSHJB equation (\ref{hjb666611}).
 \begin{definition}\label{definition6234} \ \
 If $V\in C^0({\Lambda}^{d+m})$ is a viscosity solution to PHJB (\ref{hjb6666}), then we call $\{{\cal{F}}_t\}_{t\geq0}$-adapted process $\bar{V}(t,x):=V(W_t,x)$ defined by (\ref{value62}) is a
 viscosity solution to BSHJB equation (\ref{hjb666611}).
 \end{definition}
 By Theorem \ref{theorem52}, we obtain that
 \begin{theorem}\label{theorem62}\ \
                 Let ${b}, {\sigma}, {q}, {\phi}$ satisfy Hypothesis \ref{hypstate}. Then the ${\cal{F}}_t$-adapted process $\bar{V}(t,x):=V(W_t,x)$ defined by (\ref{value62}) is a unique
 viscosity solution to  BSHJB equation (\ref{hjb666611}).
\end{theorem}
 \begin{remark}\label{remark6}
  If the  coefficients in (\ref{hjb666611}) are independent of $x$ and $u$, the  BSHJB equation (\ref{hjb666611}) reduces to a
 BSDE:
 \begin{eqnarray}\label{bsde6}
\begin{cases}
 d\bar{V}(t)=-
                                         \bar{q}(W_t,\bar{V}(t),p(t))dt+p(t)dW(t),\ \ \  t\in
                               [0,T), \\
 \bar{V}(T)=\bar{\phi}(W_T),\end{cases}\ \ \ P\mbox{-}a.s.
\end{eqnarray}
 We refer to the seminal paper by Pardoux and Peng \cite{par0} for the wellposedness of such BSDEs. Moreover, for any $(t,\gamma_t)\in [0,T]\times \Lambda$,
 by \cite{par0} the following BSDE on $[t,T]$ has a unique solution:
 \begin{eqnarray}\label{bsde677}
\begin{cases}
d\bar{V}^{\gamma_t}(s)=-
                                         \bar{q}(W^{\gamma_t}_s,\bar{V}^{\gamma_t}(s),p^{\gamma_t}(s))ds+p^{\gamma_t}(s)dW(s),\   s\in
                               [0,T),\\
 \bar{V}^{\gamma_t}(T)=\bar{\phi}(W^{\gamma_t}_T),
\end{cases} \ \  P\mbox{-}a.s,
\end{eqnarray}
 where
 \begin{eqnarray*}
W^{\gamma_t}_s(l)=\begin{cases}
\gamma_t(l),\ \  \ \ \ \ \ \ \ \ \ \ \ \ \ \ \ \  \ \ \ l\in [0,t],\\
  \gamma_t(t)+W(l)-W(t),\ \ l\in (t,T].\end{cases}
\end{eqnarray*}
 Define $V(\gamma_t):=\bar{V}^{\gamma_t}(t)$,
 in the similar (even easier) process of the proof of Theorem \ref{theoremvexist}, we show that $V(W_t)$ is a viscosity solution to BSDE (\ref{bsde6}) in our definition.
 On the other hand, it is easy to  show that, for any $t\in [0,T]$
    $$
                           {V}(W_t)=\bar{V}(t), \ \ P\mbox{-}a.a.
    $$
    Thus,  the viscosity solution to BSDE (\ref{bsde6}) coincidences with its classical solution.
    Therefore, our definition of viscosity solution to BSHJB equation (\ref{hjb666611}) is a natural  extension of classical solution to BSDE (\ref{bsde6}).

 \end{remark}

\appendix
     \section{ Borwein-Preiss variational principle}

\setcounter{equation}{0}
\renewcommand{\theequation}{A.\arabic{equation}}

\par
   {\bf  Proof (of Lemma \ref{theoremleft})}. \ \ Define sequences $\{(t_i,\gamma^i_{t_i})\}_{i\geq1}$ and $\{B_i\}_{i\geq1}$ inductively starting with
\begin{eqnarray}\label{left1}
 B_0:=\{(s,\gamma_s)\in [t_0,T]\times \Lambda^{t_0}| \ f(\gamma_s)-\delta_0\rho(\gamma_s,\gamma^0_{t_0})\geq f(\gamma^0_{t_0})\}.
\end{eqnarray}
 Since $(t_0,\gamma^0_{t_0})\in B_0$, $B_0$ is nonempty. Moreover it is closed because both $f$ and $-\rho(\cdot,\gamma^0_{t_0})$ are upper semicontinuous functionals. We also have that, for
  all $(s,\gamma_s)\in B_0$,
\begin{eqnarray}\label{left2}
\delta_0\rho(\gamma_s,\gamma^0_{t_0})\leq f(\gamma_s)-f(\gamma^0_{t_0})\leq \sup_{(s,\gamma_s)\in [t,T]\times \Lambda^t}f(\gamma_s)-f(\gamma^0_{t_0})\leq \varepsilon.
\end{eqnarray}
 Take $(t_1,\gamma^1_{t_1})\in B_0$ such that
 \begin{eqnarray}\label{left21111}
f(\gamma^1_{t_1})-\delta_0\rho(\gamma^1_{t_1},\gamma^0_{t_0})\geq \sup_{(s,\gamma_s)\in B_0}[f(\gamma_s)-\delta_0\rho(\gamma_s,\gamma^0_{t_0})]-\frac{\delta_1\varepsilon}{2\delta_0},
\end{eqnarray}
 and define similarly
 \begin{eqnarray}\label{left3}
 B_1:=\bigg{\{}(s,\gamma_s)\in B_0\cap [t_1,T]\times \Lambda^{t_1}\ \bigg{|} \ f(\gamma_s)-\sum_{k=0}^{1}\delta_k\rho(\gamma_s,\gamma^k_{t_k})\geq f(\gamma^1_{t_1})-\delta_0\rho(\gamma^1_{t_1},\gamma^0_{t_0})\bigg{\}}.
\end{eqnarray}
 In general, suppose that we have defined $(t_j,\gamma^j_{t_j})$, $B_j$ for $j=1,2,\ldots, i-1$ satisfying
  \begin{eqnarray}\label{left4}
f(\gamma^j_{t_j})-\sum_{k=0}^{j-1}\delta_k\rho(\gamma^j_{t_j},\gamma^k_{t_k})\geq \sup_{(s,\gamma_s)\in B_{j-1}}\bigg{[}f(\gamma_s)-\sum_{k=0}^{j-1}\delta_k\rho(\gamma_{s},\gamma^k_{t_k})\bigg{]}-\frac{\delta_j\varepsilon}{2^j\delta_0},
\end{eqnarray}
 and
 \begin{eqnarray}\label{left5}
 B_j:=\bigg{\{}(s,\gamma_s)\in B_{j-1}\cap [t_j,T]\times \Lambda^{t_j}\ \bigg{|} \ f(\gamma_s)-\sum_{k=0}^{j}\delta_k\rho(\gamma_s,\gamma^k_{t_k})\geq f(\gamma^j_{t_j})-\sum_{k=0}^{j-1}\delta_k\rho(\gamma^j_{t_j},\gamma^k_{t_k})\bigg{\}}.
\end{eqnarray}
 We choose $(t_i,\gamma^i_{t_i})\in B_{i-1}$ such that
  \begin{eqnarray}\label{left6}
f(\gamma^i_{t_i})-\sum_{k=0}^{i-1}\delta_k\rho(\gamma^i_{t_i},\gamma^k_{t_k})\geq \sup_{(s,\gamma_s)\in B_{i-1}}\bigg{[}f(\gamma_s)-\sum_{k=0}^{i-1}\delta_k\rho(\gamma_{s},\gamma^k_{t_k})\bigg{]}-\frac{\delta_i\varepsilon}{2^i\delta_0},
\end{eqnarray}
 and we define
  \begin{eqnarray}\label{left7}
 B_i:=\bigg{\{}(s,\gamma_s)\in B_{i-1}\cap [t_i,T]\times \Lambda^{t_i}\ \bigg{|} \ f(\gamma_s)-\sum_{k=0}^{i}\delta_k\rho(\gamma_s,\gamma^k_{t_k})\geq f(\gamma^i_{t_i})-\sum_{k=0}^{i-1}\delta_k\rho(\gamma^i_{t_i},\gamma^k_{t_k})\bigg{\}}.
\end{eqnarray}
 We can see that  for every $i=1,2,\ldots$,  $B_i$ is a closed and nonempty set. It follows from (\ref{left6}) and (\ref{left7}) that, for all $(s,\gamma_s)\in B_i$,
 \begin{eqnarray*}
 \delta_i\rho(\gamma_s,\gamma^i_{t_i})&\leq& \bigg{[}f(\gamma_s)-\sum_{k=0}^{i-1}\delta_k\rho(\gamma_s,\gamma^k_{t_k})\bigg{]}-\bigg{[}f(\gamma^i_{t_i})-\sum_{k=0}^{i-1}\delta_k\rho(\gamma^i_{t_i},\gamma^k_{t_k})\bigg{]}\\
 &\leq&\sup_{(s,\gamma_s)\in B_{i-1}}\bigg{[}f(\gamma_s)-\sum_{k=0}^{i-1}\delta_k\rho(\gamma_s,\gamma^k_{t_k})\bigg{]}-\bigg{[}f(\gamma^i_{t_i})-\sum_{k=0}^{i-1}\delta_k\rho(\gamma^i_{t_i},\gamma^k_{t_k})\bigg{]}\leq \frac{\delta_i\varepsilon}{2^i\delta_0},
\end{eqnarray*}
which implies that
 \begin{eqnarray}\label{left8}
\rho(\gamma_s,\gamma^i_{t_i})\leq \frac{\varepsilon}{2^i\delta_0},\ \ \mbox{for all}\  (s,\gamma_s)\in B_i.
\end{eqnarray}
Since $\rho$ is a gauge-type function, inequality (\ref{left8}) implies that $\sup_{(s,\gamma_s)\in B_i}d_\infty(\gamma_s, \gamma^i_{t_i})\rightarrow0$ as $i\rightarrow \infty$, and therefore,
$\sup_{(s,\gamma_s),(l,\eta_l)\in B_i}d_\infty(\gamma_s, \eta_l)\rightarrow0$ as $i\rightarrow \infty$. Since $[t,T]\times\Lambda^t$ is complete, by Cantor's intersection theorem there exists a unique
 $(\hat{t},\hat{\gamma}_{\hat{t}})\in \bigcap_{i=0}^{\infty}B_i$. Obviously, we have $d_\infty(\gamma^i_{t_i}, \hat{\gamma}_{\hat{t}})\rightarrow0$ and $t_i\uparrow \hat{t}$  as $i\rightarrow \infty$. Then $(\hat{t},\hat{\gamma}_{\hat{t}})$
  satisfies (i) by (\ref{left2}) and  (\ref{left8}). For any $(s,\gamma_s)\in [\hat{t},T]\times \Lambda^{\hat{t}}$ and $(s,\gamma_s)\neq (\hat{t},\hat{\gamma}_{\hat{t}})$, we have $(s,\gamma_s)\notin\bigcap_{i=0}^{\infty}B_i$, and therefore, for some $j$,
\begin{eqnarray}\label{left9}
f(\gamma_s)-\sum_{k=0}^{\infty}\delta_k\rho(\gamma_s,\gamma^k_{t_k})\leq f(\gamma_s)-\sum_{k=0}^{j}\delta_k\rho(\gamma_s,\gamma^k_{t_k})<
 f(\gamma^j_{t_j})-\sum_{k=0}^{j-1}\delta_k\rho(\gamma^j_{t_j},\gamma^k_{t_k}).
\end{eqnarray}
On the other hand, it follows from (\ref{left1}), (\ref{left7}) and $(\hat{t},\hat{\gamma}_{\hat{t}})\in \bigcap_{i=0}^{\infty}B_i$ that, for any $q\geq j$,
\begin{eqnarray}\label{left10}
f(\gamma^0_{t_0})&\leq& f(\gamma^j_{t_j})-\sum_{k=0}^{j-1}\delta_k\rho(\gamma^j_{t_j},\gamma^k_{t_k})\leq f(\gamma^q_{t_q})-\sum_{k=0}^{q-1}\delta_k\rho(\gamma^q_{t_q},\gamma^k_{t_k})\nonumber\\
&\leq&f(\hat{\gamma}_{\hat{t}})-\sum_{k=0}^{q}\delta_k\rho(\hat{\gamma}_{\hat{t}},\gamma^k_{t_k}).
\end{eqnarray}
Letting $q\rightarrow\infty$ in (\ref{left10}), we obtain
\begin{eqnarray}\label{left11}
f(\gamma^0_{t_0})\leq  f(\gamma^j_{t_j})-\sum_{k=0}^{j-1}\delta_k\rho(\gamma^j_{t_j},\gamma^k_{t_k})\leq f(\hat{\gamma}_{\hat{t}})-\sum_{k=0}^{\infty}\delta_k\rho(\hat{\gamma}_{\hat{t}},\gamma^k_{t_k}),
\end{eqnarray}
which verifies (ii). Combining  (\ref{left9}) and (\ref{left11}) yields (iii). \ \ $\Box$

   \section{ Existence and consistency for viscosity solutions.}

\setcounter{equation}{0}
\renewcommand{\theequation}{B.\arabic{equation}}

 {\bf  Proof (of Theorem \ref{theoremvexist})}. \ \
 We let  $\varphi\in \mathcal{A}^+(\hat{\gamma}_{\hat{t}},V)$
                  with
                   $(\hat{t},\hat{\gamma}_{\hat{t}})\in [0,T)\times \Lambda$. 
For  $0< \delta\leq T-\hat{t}$, we have $\hat{t}< \hat{t}+\delta \leq T$, then by the DPP (Theorem \ref{theoremddp}), we obtain the following result:
 \begin{eqnarray}\label{4.9}
                           && 0=V(\hat{\gamma}_{\hat{t}})-{{\varphi}} (\hat{\gamma}_{\hat{t}})
                           =\mathop{\esssup}\limits_{u(\cdot)\in {\cal{U}}[\hat{t},\hat{t}+\delta]} G^{\hat{\gamma}_{\hat{t}},u}_{{\hat{t}},t+\delta}[V(X^{\hat{\gamma}_{\hat{t}},u}_{{\hat{t}}+\delta})]
                           -{{\varphi}} (\hat{\gamma}_{\hat{t}}).
\end{eqnarray}
Then, for any $\varepsilon>0$ and $0<\delta\leq T-\hat{t}$,  we can  find a control  ${u}^{\varepsilon}(\cdot)\equiv u^{{\varepsilon},\delta}(\cdot)\in {\cal{U}}[\hat{t},\hat{t}+\delta]$ such
   that the following result holds:
\begin{eqnarray}\label{4.10}
    -{\varepsilon}\delta
    \leq G^{\hat{\gamma}_{\hat{t}},{u}^{\varepsilon}}_{{\hat{t}},\hat{t}+\delta}[V(X^{\hat{\gamma}_{\hat{t}},{u}^{\varepsilon}}_{{\hat{t}}+\delta})]-{{\varphi}} (\hat{\gamma}_{\hat{t}}).
\end{eqnarray}
 We note that
                     $G^{\hat{\gamma}_{\hat{t}},{u}^{\varepsilon}}_{s,\hat{t}+\delta}[V(X^{\hat{\gamma}_{\hat{t}},{u}^{\varepsilon}}_{{\hat{t}}+\delta})]$
                     is defined in terms of the solution of the
                     BSDE:
 \begin{eqnarray}\label{bsde4.10}
\begin{cases}
dY^{\hat{\gamma}_{\hat{t}},{u}^{\varepsilon}}(s) =
               -q(X^{\hat{\gamma}_{\hat{t}},{u}^{\varepsilon}}_s,Y^{\hat{\gamma}_{\hat{t}},{u}^{\varepsilon}}(s),Z^{\hat{\gamma}_{\hat{t}},
               {u}^{\varepsilon}}(s),{u}^{\varepsilon}(s))ds+Z^{\hat{\gamma}_{\hat{t}},{u}^{\varepsilon}}(s)dW(s),\ \  s\in[\hat{t},\hat{t}+\delta], \\
 ~Y^{\hat{\gamma}_{\hat{t}},{u}^{\varepsilon}}(\hat{t}+\delta)=V(X^{\hat{\gamma}_{\hat{t}},
 {u}^{\varepsilon}}_{{\hat{t}}+\delta}),\end{cases}
\end{eqnarray}
                     by the following formula:
$$
                        G^{\hat{\gamma}_{\hat{t}},{u}^{\varepsilon}}_{s,\hat{t}+\delta}[V(X^{\hat{\gamma}_{\hat{t}},{u}^{\varepsilon}}_{{\hat{t}}+\delta})]
                        =Y^{\hat{\gamma}_{\hat{t}},{u}^{\varepsilon}}(s),\
                        \ \ s\in[\hat{t},\hat{t}+\delta].
$$
                         Applying functional It\^{o} formula (\ref{statesop}) to
                         ${\varphi}(X^{\hat{\gamma}_{\hat{t}},{u}^{\varepsilon}}_s)$,   we get that
\begin{eqnarray}\label{bsde4.21}
                             {\varphi}(X^{\hat{\gamma}_{\hat{t}},{u}^{\varepsilon}}_s)
                             &=& {\varphi}(\hat{\gamma}_{\hat{t}})+\int^{s}_{{\hat{t}}} (\tilde{\cal{L}}{\varphi})(X^{\hat{\gamma}_{\hat{t}},{u}^{\varepsilon}}_l,
                                      u^{{\varepsilon}}(l))dl -\int^{s}_{{\hat{t}}}q(X^{\hat{\gamma}_{\hat{t}},{u}^{\varepsilon}}_l,{\varphi}(X^{\hat{\gamma}_{\hat{t}},
                                      {u}^{\varepsilon}}_l),\sigma^\top(X^{\hat{\gamma}_{\hat{t}},{u}^{\varepsilon}}_l,u^{{\varepsilon}}(l))\nonumber\\
                             &&
                             \times \partial_x{\varphi}(X^{\hat{\gamma}_{\hat{t}},{u}^{\varepsilon}}_l), u^{\varepsilon}(l))dl
                             +\int^{s}_{{\hat{t}}}[
                             \sigma^\top(X^{\hat{\gamma}_{\hat{t}},{u}^{\varepsilon}}_l,u^{\varepsilon}(l))\partial_x{\varphi}(X^{\hat{\gamma}_{\hat{t}},{u}^{\varepsilon}}_l)]dW(l),\ \
\end{eqnarray}
where
\begin{eqnarray*}
                       (\tilde{\cal{L}}{\varphi})(\gamma_t,u)
                       &=&\partial_t{\varphi}(\gamma_t)+
                                         \langle\partial_x {\varphi}(\gamma_t),b(\gamma_t,u)\rangle+\frac{1}{2}\mbox{tr}[\partial_{xx}{\varphi}(\gamma_t)\sigma(\gamma_t,u)\sigma^{\top}(\gamma_t,u)]\\
                       &&+
                                         q(\gamma_t,{\varphi}(\gamma_t),\sigma^{\top}(\gamma_t,u){\partial_x{\varphi}(\gamma_t)},u),
                                         ~~~~ \ (t, \gamma_t,u)\in [0,T]\times {\Lambda}\times U.
\end{eqnarray*}
Set
\begin{eqnarray*}
                          &&Y^{2,{\hat{\gamma}_{\hat{t}},{u}^{\varepsilon}}}(s):=
                          {\varphi}(X_s^{\hat{\gamma}_{\hat{t}},{u}^{\varepsilon}})-Y^{\hat{\gamma}_{\hat{t}},u^{\varepsilon}}(s),  \ \ s\in[\hat{t},\hat{t}+\delta],\\
                            &&Z^{2,{\hat{\gamma}_{\hat{t}},{u}^{\varepsilon}}}(s)  :=
                             \sigma^\top(X^{\hat{\gamma}_{\hat{t}},{u}^{\varepsilon}}_s,u^{\varepsilon}(s))\partial_x{\varphi}(X^{\hat{\gamma}_{\hat{t}},{u}^{\varepsilon}}_s)-Z^{\hat{\gamma}_{\hat{t}},u^{\varepsilon}}(s), \ \ s\in[\hat{t},\hat{t}+\delta].
\end{eqnarray*}
Comparing (\ref{bsde4.10}) and (\ref{bsde4.21}), we have, $P$-a.s.,
\begin{eqnarray*}
                      &&dY^{2,{\hat{\gamma}_{\hat{t}},{u}^{\varepsilon}}}(s)\\
                      &=&[(\tilde{\cal{L}}{\varphi})(X^{\hat{\gamma}_{\hat{t}},{u}^{\varepsilon}}_s,u^{\varepsilon}(s))
                       -q(X^{\hat{\gamma}_{\hat{t}},{u}^{\varepsilon}}_s,{\varphi}(X^{\hat{\gamma}_{\hat{t}},{u}^{\varepsilon}}_s),
                             \sigma^\top(X^{\hat{\gamma}_{\hat{t}},{u}^{\varepsilon}}_s,{u}^{\varepsilon}(s))\partial_x{\varphi}(X^{\hat{\gamma}_{\hat{t}},{u}^{\varepsilon}}_s), {u}^{\varepsilon}(s))\\
                             &&+q(X^{\hat{\gamma}_{\hat{t}},{u}^{\varepsilon}}_s,Y^{\hat{\gamma}_{\hat{t}},{u}^{\varepsilon}}(s),
                                  Z^{\hat{\gamma}_{\hat{t}},{u}^{\varepsilon}}(s),{u}^{\varepsilon}(s))]ds
                                +Z^{2,{\hat{\gamma}_{\hat{t}},{u}^{\varepsilon}}}(s)dW(s)\\
                                &=&[(\tilde{\cal{L}}{\varphi})(X^{\hat{\gamma}_{\hat{t}},{u}^{\varepsilon}}_s,u^{\varepsilon}(s))
                                -A(s)Y^{2,{\hat{\gamma}_{\hat{t}},{u}^{\varepsilon}}}(s)-(\bar{A}(s),Z^{2,{\hat{\gamma}_{\hat{t}},{u}^{\varepsilon}}}(s))_{\mathbb{R}^n}]ds
                                +Z^{2,{\hat{\gamma}_{\hat{t}},{u}^{\varepsilon}}}(s)dW(s),
\end{eqnarray*}
where $|A|\vee|\bar{A}|\leq L$. 
Therefore, we obtain (see Proposition 2.2 in \cite{el1})
\begin{eqnarray}\label{4.14}
                      Y^{2,{\hat{\gamma}_{\hat{t}},{u}^{\varepsilon}}}(\hat{t})
                      =\mathbb{E}\bigg{[}Y^{2,{\hat{\gamma}_{\hat{t}},{u}^{\varepsilon}}}(\hat{t}+\delta)\Gamma^{\hat{t}}(\hat{t}+\delta)
                      -
                      \int^{\hat{t}+\delta}_{{\hat{t}}}\Gamma^{\hat{t}}(l)(\tilde{\cal{L}}{\varphi})
                       (X^{\hat{\gamma}_{\hat{t}},{u}^{\varepsilon}}_l,u^{\varepsilon}(l))dl\bigg{|}
                      {\cal{F}}_{\hat{t}}\bigg{]},\ \ \
\end{eqnarray}
where $\Gamma^{\hat{t}}(\cdot)$ solves the linear SDE
$$
               d\Gamma^{\hat{t}}(s)=\Gamma^{\hat{t}}(s)(A(s)ds+\bar{A}(s)dW(s)),\ s\in [{\hat{t}},{\hat{t}}+\delta];\ \ \ \Gamma^{\hat{t}}({\hat{t}})=1.
$$
Obviously, $\Gamma^{\hat{t}}\geq 0$. Combining (\ref{4.10}) and (\ref{4.14}), we have
\begin{eqnarray}\label{4.15}
-\varepsilon
    &\leq& \frac{1}{\delta}\mathbb{E}\bigg{[}-Y^{2,{\hat{\gamma}_{\hat{t}},{u}^{\varepsilon}}}(\hat{t}+\delta)\Gamma^{\hat{t}}(\hat{t}+\delta)+
    \int^{\hat{t}+\delta}_{{\hat{t}}}\Gamma^{\hat{t}}(l)(\tilde{\cal{L}}{\varphi})(X^{\hat{\gamma}_{\hat{t}},{u}^{\varepsilon}}_l,u^{\varepsilon}(l))dl\bigg{]}\nonumber\\
    &=&-\frac{1}{\delta}\mathbb{E}\bigg{[}Y^{2,{\hat{\gamma}_{\hat{t}},{u}^{\varepsilon}}}(\hat{t}+\delta)\Gamma^{\hat{t}}(\hat{t}+\delta)\bigg{]}+
      \frac{1}{\delta}\mathbb{E}\bigg{[}\int^{\hat{t}+\delta}_{{\hat{t}}}(\tilde{\cal{L}}{\varphi})(\hat{\gamma}_{{\hat{t}}},u^{\varepsilon}(l))dl\bigg{]}\nonumber\\
                      &&+\frac{1}{\delta}\mathbb{E}\bigg{[}\int^{\hat{t}+\delta}_{{\hat{t}}}[(\tilde{\cal{L}}{\varphi})(X^{\hat{\gamma}_{\hat{t}},{u}^{\varepsilon}}_l,u^{\varepsilon}(l))-
                      (\tilde{\cal{L}}{\varphi})(\hat{\gamma}_{{\hat{t}}},u^{\varepsilon}(l))]dl\bigg{]}\nonumber\\
                      &&+\frac{1}{\delta}\mathbb{E}\bigg{[}\int^{\hat{t}+\delta}_{{\hat{t}}}(\Gamma^{\hat{t}}(l)-1)(\tilde{\cal{L}}{\varphi})(X^{\hat{\gamma}_{\hat{t}},{u}^{\varepsilon}}_l,
                                u^{\varepsilon}(l))dl\bigg{]}\nonumber\\
    &:=&I+II+III+IV.
\end{eqnarray}
Since the coefficients in $\tilde{\cal{L}}$  satisfy  linear growth  condition,
 combining the regularity of $\varphi\in C^{1,2}_{p}({\Lambda}^{\hat{t}})$, there exist a integer
$\bar{p}\geq1$ and a constant $C>0$ independent of $u\in U$ such that, for all $(t,\gamma_t,u)\in [0,T]\times\Lambda\times U$,
\begin{eqnarray}\label{4.4444}|{\varphi}(\gamma_{t})|\vee  |
                      (\tilde{\cal{L}}{\varphi})(\gamma_{t},u)|
                      \leq  C(1+||\gamma_{t}||_0)^{\bar{p}}.
\end{eqnarray}
In view of Lemma \ref{lemmaexist}, we also have
\begin{eqnarray*}
                          \sup_{u(\cdot)\in{\cal{U}}[\hat{t},\hat{t}+\delta]}\mathbb{E}[\sup_{{\hat{t}}\leq s\leq \hat{t}+\delta}|\Gamma^{\hat{t}}(s)-1|^2]\leq C\delta.
\end{eqnarray*}
Thus, by $\varphi\in \mathcal{A}^+(\hat{\gamma}_{\hat{t}},V)$,
\begin{eqnarray}\label{4.1611}
I&=& -\frac{1}{\delta}\mathbb{E}\bigg{[}\bigg{(}{\varphi}(X_{\hat{t}+\delta}^{\hat{\gamma}_{\hat{t}},{u}^{\varepsilon}})-Y^{\hat{\gamma}_{\hat{t}},u^{\varepsilon}}(\hat{t}+\delta)\bigg{)}\Gamma^{\hat{t}}(\hat{t}+\delta)\bigg{]}\nonumber\\
&=&\frac{1}{\delta}\mathbb{E}\bigg{[}\bigg{(}V(X^{\hat{\gamma}_{\hat{t}},
 {u}^{\varepsilon}}_{{\hat{t}}+\delta})-{\varphi}(X_{\hat{t}+\delta}^{\hat{\gamma}_{\hat{t}},{u}^{\varepsilon}})\bigg{)}\Gamma^{\hat{t}}(\hat{t}+\delta)\bigg{]}\leq 0;
\end{eqnarray}
\begin{eqnarray}\label{4.16}
             II\leq\frac{1}{\delta}\bigg{[}\int^{\hat{t}+\delta}_{{\hat{t}}}\sup_{u\in U}(\tilde{\cal{L}}{\varphi})({\hat{\gamma}_{\hat{t}},{u}})dl\bigg{]}
                      = {\mathcal{L}}{\varphi}(\hat{\gamma}_{\hat{t}}).
\end{eqnarray}
Now we estimate higher order terms  $III$ and $IV$.
By (\ref{4.4444}) and the dominated convergence theorem,
$$\lim_{\delta\rightarrow0}
 \mathbb{E}\sup_{\hat{t}\leq l\leq\hat{t}+\delta}|(\tilde{\cal{L}}{\varphi})(X^{{\hat{\gamma}_{\hat{t}},{u}^{\varepsilon}}}_l,{u}^{\varepsilon}(l))-
                      (\tilde{\cal{L}}{\varphi})({\hat{\gamma}_{\hat{t}},{u}^{\varepsilon}}(l))|=0,
$$
then
\begin{eqnarray}\label{4.18}
\lim_{\delta\rightarrow0}|III| &\leq&\lim_{\delta\rightarrow0}\sup_{\hat{t}\leq l\leq\hat{t}+\delta}\mathbb{E}|(\tilde{\cal{L}}{\varphi})(X^{{\hat{\gamma}_{\hat{t}},{u}^{\varepsilon}}}_l,{u}^{\varepsilon}(l))-
                      (\tilde{\cal{L}}{\varphi})({\hat{\gamma}_{\hat{t}},{u}^{\varepsilon}}(l))|=0;
\end{eqnarray}
and, for some finite constant $C>0$,
\begin{eqnarray}\label{4.19}
|IV|
&\leq&\frac{1}{\delta}\int^{\hat{t}+\delta}_{{\hat{t}}}\mathbb{E}|\Gamma^{\hat{t}}(l)-1||(\tilde{\cal{L}}{\varphi})(X^{\hat{\gamma}_{\hat{t}},{u}^{\varepsilon}}_l,
                               u^{\varepsilon}(l)|
                       dl\nonumber\\
                       &\leq&\frac{1}{\delta}\int^{\hat{t}+\delta}_{{\hat{t}}}(\mathbb{E}(\Gamma^{\hat{t}}(l)-1)^2)^{\frac{1}{2}}(\mathbb{E}((\tilde{\cal{L}}{\varphi})(X^{\hat{\gamma}_{\hat{t}},
                                     {u}^{\varepsilon}}_l,
                               u^{\varepsilon}(l))^2)^{\frac{1}{2}}dl
                    \nonumber\\
                       &\leq& C(1+||\hat{\gamma}_{\hat{t}}||_0)^{\bar{p}}\delta^\frac{1}{2}.
\end{eqnarray}
Substituting  (\ref{4.1611}), (\ref{4.16}) and (\ref{4.19}) into (\ref{4.15}), we have
\begin{eqnarray}\label{4.2000000}
-\varepsilon\leq {\mathcal{L}}{\varphi}(\hat{\gamma}_{\hat{t}})
+III
    +C(1+||\hat{\gamma}_{\hat{t}}||_0)^{\bar{p}}\delta^\frac{1}{2}.
\end{eqnarray}
Sending $\delta$ to $0$, by (\ref{4.18}), we have
$$
-\varepsilon
    \leq {\mathcal{L}}{\varphi}(\hat{\gamma}_{\hat{t}}).
$$
By the arbitrariness of $\varepsilon$, we show
 $V$ is  a viscosity subsolution to (\ref{hjb1}).
 \par
 In a symmetric (even easier) way, we show that  $V$ is also a viscosity supersolution to equation (\ref{hjb1}).
                 This step  completes the proof.\ \ $\Box$

{\bf  Proof (of Theorem \ref{theorem3.2})}. \ \ We prove the subsolution property only.  Assume $v$ is a viscosity subsolution. It is clear that $v(\gamma_T)\leq\phi(\gamma_T)$ for all $\gamma_T\in \Lambda_T$. For any $(t,\gamma_t)\in [0,T)\times \Lambda$, since $v\in C_p^{1,2}({\Lambda})$, by definition of viscosity subsolutions we see that ${\mathcal{L}}v(\gamma_t)\geq0$.
\par
On the other hand, assume $v$ is a  classical subsolution.    Let  $\varphi\in \mathcal{A}^+(\gamma_t,v)$ with $t\in [0,T)$.
 For  every $\alpha\in \mathbb{R}^{d}$ and $\beta\in \mathbb{R}^{d\times n}$,
 let  $$
           X(s)=\gamma_t(t)+\int^{s}_{t}\alpha dl+\int^{s}_{t}\beta dW(l),\ \ s\in [t,T],
 $$
 and $X(s)=\gamma_t(s)$, $s\in [0,t)$. Then $X(\cdot)$ is a continuous semi-martingale on $[t,T]$.
  Applying functional It\^o formula (\ref{statesop}) to $\varphi$ and noticing that $(v-\varphi)(\gamma_t)=0$, we have, for every $0<\delta\leq T-t$,
 \begin{eqnarray}\label{1224}
                            0&\leq&{\mathbb{E}}(\varphi-v)(X_{t+\delta})\nonumber\\
                            &=&{\mathbb{E}}\int^{t+\delta}_{t}[\partial_t(\varphi-v)(X_l)
                 +\langle\partial_x(\varphi-v)(X_l),\alpha\rangle]dl+\frac{1}{2}{\mathbb{E}}\int^{t+\delta}_{t}\mbox{tr}((\partial_{xx}(\varphi-v)(X_l))\beta\beta^\top)dl\nonumber\\
                 &=&
                           {\mathbb{E}}\int^{t+\delta}_{t}\widetilde{{\mathcal{H}}}(X_l)
                              dl,
\end{eqnarray}
where \begin{eqnarray*}
\widetilde{{\mathcal{H}}}(\eta_s)=\partial_t(\varphi-v)(\eta_s)
                 +\langle\partial_x(\varphi-v)(\eta_s),\alpha\rangle
                 +\frac{1}{2}\mbox{tr}((\partial_{xx}(\varphi-v)(\eta_s))\beta\beta^\top), \ \ \ (s,\eta_s)\in [0,T]\times \Lambda.
                 \end{eqnarray*}
                 Letting $\delta\rightarrow0$ in (\ref{1224}),
\begin{eqnarray}\label{040912}
\widetilde{{\mathcal{H}}}(\gamma_t)\geq0.
\end{eqnarray}
Let $\beta=\mathbf{0}$, by the arbitrariness of  $\alpha$,
$$
\partial_t\varphi(\gamma_t)\geq\partial_tv(\gamma_t), \ \ \partial_x\varphi(\gamma_t)=\partial_xv(\gamma_t). 
$$
Then, for every $u\in U$, let $\beta=\sigma(\gamma_t,u)$ in (\ref{040912}).
Noting that $\varphi(\gamma_t)=v(\gamma_t)$, 
we have
\begin{eqnarray*}
 &&\partial_t\varphi(\gamma_t)+\langle\partial_x\varphi(\gamma_t),b(\gamma_t,u)\rangle+\frac{1}{2}\mbox{tr}(\partial_{xx}\varphi(\gamma_t)\sigma(\gamma_t,u)\sigma^\top(\gamma_t,u))\\
                 &&
                 +q(\gamma_t,{\varphi}(\gamma_t),\sigma^\top(\gamma_t,u){\partial_x{\varphi}(\gamma_t)},u)\\
                 &\geq& \partial_tv(\gamma_t)+\langle\partial_xv(\gamma_t),b(\gamma_t,u)\rangle
                 +\frac{1}{2}\mbox{tr}(\partial_{xx}v(\gamma_t)\sigma(\gamma_t,u)\sigma^\top(\gamma_t,u))\\&&+q(\gamma_t,v(\gamma_t),\sigma^\top(\gamma_t,u){\partial_xv(\gamma_t)},u).
\end{eqnarray*}
Note that ${\mathcal{L}}v(\gamma_t)\geq0$,  
                             taking the supremum over $u\in U$, we see that
\begin{eqnarray*}
{\mathcal{L}}\varphi(\gamma_t)\geq{\mathcal{L}}v(\gamma_t)\geq0.
\end{eqnarray*}
 Thus,
   we have that $v$ is a viscosity subsolution of equation (\ref{hjb1}). 
    \ \ $\Box$

\par

\end{document}